\documentclass{article}%
\usepackage{amsfonts}
\usepackage{amsmath}
\usepackage{amssymb}
\usepackage{graphicx}%
\newtheorem{theorem}{Theorem}

\newtheorem{corollary}[theorem]{Corollary}

\newtheorem{definition}[theorem]{Definition}

\newtheorem{lemma}[theorem]{Lemma}

\newtheorem{proposition}[theorem]{Proposition}
\newtheorem{remark}[theorem]{Remark}

\begin{document}

\author{Igor B. Frenkel\\Department of Mathematics, Yale University\\New Haven, CT 06520
\and Andrey N. Todorov\\Department of Mathematics\\University of California\\Santa Cruz, CA 95064\\Bulgarian Academy of Sciences\\Institute of Mathematics,\\Sofia, Bulgaria}
\title{Complex Counterpart of Chern-Simons-Witten Theory and Holomorphic Linking.}
\maketitle

\begin{abstract}
In this paper we are begining to explore the complex counterpart of the
Chern-Simon-Witten theory. We define the complex analogue of \ the Gauss
linking number for complex curves embedded in a Calabi-Yau threefold using the
formal path integral that leads to a rigorous mathematical expression. We give
an analytic and geometric interpretation of our holomorphic linking following
the parallel with the real case. We show in particular that the Green kernel
that appears in the explicit integral for the Gauss linking number is replaced
by the Bochner-Martinelli kernel. We also find canonical expressions of the
holomorphic linking using\ the Grothendieck-Serre duality in local cohomology,
the latter admits a generalization for an arbitrary field.

\end{abstract}
\tableofcontents

\section{Introduction}

The relation between loop groups and their central extensions,
Wess-Zumino-Novikov-Witten (WZNW) two dimensional conformal field theories and
Chern-Simons-Witten (CSW) three dimensional topological theories, emerged as a
unifying principle among various areas in mathematics and theoretical physics.
It has gradually become clear that the three components and the relation
itself admit a remarkable complexification, which combines in a profound way
some further areas in both disciplines. It was shown in \cite{EF} that the
classification of coadjoint orbits for a new class of two dimensional current
groups on Riemann surfaces can be viewed as a classification of stable vector
bundles over these surfaces. In \cite{FK} an analogue of WZNW construction for
two dimensional current groups was obtained by means of Leray's residue theory
in complex analysis cf \cite{LMNS}.

In this paper we begin to explore a complex counterpart of CSW theory. We
consider a generalization of the simplest invariant of two curves in S$^{3}$,
namely the Gauss linking number, which arises in the abelian CSW theory as a
correlation function of holonomy functionals corresponding to two curves. In
the abelian case one can give a precise meaning to the formal path integral
and derive the familiar formula for the Gauss linking number. The
\textquotedblright complexification\textquotedblright\ of loop group theory
discovered in \cite{EF} and extended in \cite{FK} leads to a complex
counterpart of the abelian CSW path integral, which we turn into a rigorous
expression for the holomorphic linking. Then we show that the Green kernel
that appears in the classical integral for the Gauss linking number is
replaced by the Bochner-Martinelli kernel and has deep relation to the theory
of Green currents in Arakelov geometry. The integral formula for the Gauss
linking number leads to its topological realization as an intersection number
and we derive its algebro-geometric analogue in the complex case. It turns out
that the notion of holomorphic linking is related to many structures of
complex and algebraic geometry, which can be viewed as complementary aspects
of one unified picture. At one end the holomorphic linking is presented by the
complex Chern-Simon-Witten path integral, at the other end, it is expressed
via the Grothendieck-Serre duality in local cohomology. The goal of this paper
is to demonstrate different realizations of the holomorphic linking and its
connections with various structures of mathematics and mathematical physics.

Let us now explain the notion of holomorphic linking in some detail. By
definition the holomorphic linking of two complex curves $\Sigma_{1}$ and
$\Sigma_{2}$ in a Calabi-Yau (CY) manifold M is a linear map on the product of
the spaces of holomorphic differentials on $\Sigma_{1}$ and $\Sigma_{2}.$ It
is no longer a topological invariant but depends only on the complex structure
on $\Sigma_{1}$ and $\Sigma_{2}$ and their embedding into a CY manifold M and
not on the metric. To illustrate our notion of holomorphic linking, we
consider an example of the simplest non-compact CY manifold $\mathbb{C}^{3}$
and two affine curves $\Sigma_{1}$ and $\Sigma_{2}$ embedded in it. Then the
analogue of the Gauss formula is the following expression for the holomorphic
linking:
\[
\#((\Sigma_{1},\theta_{1}),(\Sigma_{2},\theta_{2}))=\int_{\Sigma_{1}%
\times\Sigma_{2}}\frac{\varepsilon^{ijk}\overline{(z_{i}-w_{i})}}{|z-w|^{6}%
}\wedge d\overline{z_{j}}\wedge d\overline{w_{k}}\wedge\theta_{1}\wedge
\theta_{2},
\]
where $\varepsilon^{ijk}$ is the sign of the permutation $(i,j,k)$ and
$\theta_{1}$ and $\theta_{2}$ are holomorphic forms on $\Sigma_{1}$ and
$\Sigma_{2},$ respectively. When $\Sigma_{1}$ and $\Sigma_{2}$ are complex
lines we can compute the holomorphic linking and compare it with the
\textquotedblright Gauss linking number\textquotedblright\ of two real lines
in $\mathbb{R}^{3}.$ When two lines are not parallel to each other the
integral formulas both in real and in complex cases can be easily identified
as follows. We pick three vectors $\overrightarrow{e_{1}},\overrightarrow
{e_{2}}$ and $\overrightarrow{e_{3}}$ such that the first two determine the
direction of two lines and the third vector has the initial point on the first
line and the end point on the second line. Then the Gauss integral formula
yields $\frac{1}{2}sign(\det(\overrightarrow{e_{1}},\overrightarrow{e_{2}}$,
$\overrightarrow{e_{3}}))$ and it depends on the choice of the orientation of
the first and the second line, but not on the order of the two lines. A
similar computation in the complex case gives, up to a scalar multiple
$\det(\overrightarrow{e_{1}},\overrightarrow{e_{2}}$, $\overrightarrow{e_{3}%
})^{-1}$ $<\overrightarrow{e_{1}},\theta_{1}><\overrightarrow{e_{2}}%
,\theta_{2}>,$ where $\theta_{i}$ for i=1 and 2 are elements of the dual space
and can be viewed as 1-forms. Thus, the main information of the holomorphic
linking is contained in the determinant. The latter degenerates when two lines
cross each other or become parallel, and is also covariant with respect to the
linear transformation. This elementary example indicates that the complex
linking is the \textquotedblright measure\textquotedblright\ of closeness of
two curves in a three dimensional complex manifold and is the simplest
possible invariant in the complex case. This supports our belief in the
intrinsic nature of the new invariant.

It is a well known that the Gauss linking number of two circles in $S^{3}%
\,\ $has a simple geometric interpretation as an intersection number of one
circle with a disk bounded by the second one. The analytic expression for the
holomorphic linking number can also be written in a similar form. The new
ingredients are the holomorphic forms $\theta_{1}$ and $\theta_{2}$ attached
to the curves $\Sigma_{1}$ and $\Sigma_{2}.$ Instead of a disk with prescribed
boundary circle in the real class one should consider a surface $S_{1}$ with a
prescribed divisor $\Sigma_{1}$ in the complex case. Moreover the holomorphic
form $\theta_{1}$ on $\Sigma_{1}$ is "lifted" to a meromorphic form
$\omega_{1}$ on $S_{1}$ such that $res(\omega_{1})=\theta_{1}.$ See \cite{FK}
and a further generalization in \cite{HR} and \cite{HR1}. Then the formula for
the holomorphic linking in M becomes the following%
\[
\#\left(  \left(  \Sigma_{1},\theta_{1}\right)  ,\left(  \Sigma_{2},\theta
_{2}\right)  \right)  =%
%TCIMACRO{\dsum \limits_{x\in S_{1}\cap\Sigma_{2}}}%
%BeginExpansion
{\displaystyle\sum\limits_{x\in S_{1}\cap\Sigma_{2}}}
%EndExpansion
\frac{\omega_{1}(x)\wedge\theta_{2}(x)}{\eta(x)},
\]
where $\omega_{1}(x)\wedge\theta_{2}(x)\in\wedge^{3}\left(  T_{x,\text{M}%
}^{1,0}\right)  ^{\ast}\approxeq\Omega_{x,\text{M}}^{3}$ and $\eta$ is a
holomorphic volume form on M. As in the real case, the invariant does not
depend on the choice of the complex surface $S_{1}$ and the meromorphic form
$\omega_{1}.$ It is easy to see directly that when $\Sigma_{1}$ and
$\Sigma_{2}$ are complex lines and $S_{1}$ is a complex plane containing
$\Sigma_{1}$ the above geometric formula for the holomorphic linking yields up
to a scalar multiple the expression $\det(\overrightarrow{e_{1}}%
,\overrightarrow{e_{2}}$, $\overrightarrow{e_{3}})^{-1}\left\langle
\overrightarrow{e_{1}},\theta_{1}\right\rangle \left\langle \overrightarrow
{e_{2}},\theta_{2}\right\rangle $ discussed above$.$

As we have mentioned before we have derived the analytic and geometric
formulas for holomorphic linking studying complex counterpart of CSW with a
path integral
\[%
%TCIMACRO{\dint \limits_{\mathcal{A}^{0,1}/\mathcal{G}^{0,1}}}%
%BeginExpansion
{\displaystyle\int\limits_{\mathcal{A}^{0,1}/\mathcal{G}^{0,1}}}
%EndExpansion
\mathcal{D}A\exp\left(  \sqrt{-1}h%
%TCIMACRO{\dint \limits_{\text{M}}}%
%BeginExpansion
{\displaystyle\int\limits_{\text{M}}}
%EndExpansion
A\wedge\overline{\partial}A\wedge\eta\right)  \exp\left(
%TCIMACRO{\dint \limits_{\Sigma_{1}}}%
%BeginExpansion
{\displaystyle\int\limits_{\Sigma_{1}}}
%EndExpansion
A\wedge\theta_{1}\right)  \exp\left(
%TCIMACRO{\dint \limits_{\Sigma_{2}}}%
%BeginExpansion
{\displaystyle\int\limits_{\Sigma_{2}}}
%EndExpansion
A\wedge\theta_{2}\right)  ,
\]
where $h$ is a real parameter, $\mathcal{A}^{0,1}$ is the space of $(0,1)$
forms on M and $\mathcal{G}^{0,1}$ is the subspace of all $\overline{\partial
}$ closed forms. Its rigorous meaning can be expressed in terms of the Green
current for the pair (M$\times$M,$\Delta$), where $\Delta$ is the diagonal.
Its existence follows from the general theory of the Green currents
established by Gillet and Soul\'{e} (see \cite{AB}) in the context of Arakelov
geometry. For any complex subvariety Y in a projective manifold X they defined
a Green current g$_{\text{Y}}$ satisfying the equation%
\[
\partial\overline{\partial}\text{g}_{\text{Y}}+\delta_{\text{Y}}%
=\omega_{\text{Y}},
\]
where $\delta_{\text{Y}}$ is the Dirac delta current corresponding to Y and
$\omega_{\text{Y}}$ is a smooth closed form representing the Poincare dual of
the homology class [Y] of Y. We use the existence of a Green current for the
pair $($X,Y$)=($M$\times$M,$\Delta)$ and then restrict it to $U\times U$,
where $U$ is an affine open set in M. The comparison of the restriction to the
explicit formula for the Bochner-Martinelli kernel on $U\times U$ yields the
above analogue of the Gauss integral formula in the complex case. Similarly a
Green current for the pairs (X,Y)$=($M,$\Sigma_{1})$ or $($M,$\Sigma_{2})$ can
be used for an alternative analytic expression of the holomorphic linking.

The holomorphic linking also admits a certain canonical expression in the
language of homological algebra. For a complex subvariety Y in a projective
manifold X and a top holomorphic form $\theta$ on Y we introduce generalized
notions of Grothendieck and Serre classes of a pair $($Y,$\theta)$ denoted by
$\mu($Y,$\theta)$ and $\lambda($Y,$\theta),$ respectively. By definition
$\mu($Y,$\theta)$ and $\lambda($Y,$\theta)$ belong to $Ext_{\mathcal{O}%
_{\text{X}}}^{d}(\mathcal{O}_{\text{Y}},\Omega_{\text{X}}^{n})$ and
$Ext_{\mathcal{O}_{\text{X}}}^{d-1}(I_{\text{Y}},\Omega_{\text{X}}^{n}%
)\,\ $where $I_{\text{Y}}$ is the ideal sheaf of Y, $\Omega_{\text{X}}^{n}$ is
the locally free sheaf of top holomorphic forms on X, $\mathcal{O}_{\text{X}}$
and $\mathcal{O}_{\text{Y}}$ are the structure sheaves on X and Y and $d$ is
the codimension of Y in X. This is a generalization \ of similar notions
originally studied by Atiyah in \cite{A} in the case of linking of two lines
in a twistor space. From the definition it follows that the Grothendieck class
always exists. The new point of our construction is the dependence of the
Grothendieck and Serre classes on the holomorphic form $\theta,$ which is
crucial for the existence of the Serre classes of pairs (Y,$\theta),$ where Y
is a submanifold in a projective manifold. In particular, we show the
existence of the Serre class of the pair (X,Y) where X=M$\times$M and Y is the
diagonal embedding of a CY manifold M with the holomorphic form $\eta$ on it.
Then the holomorphic linking of two curves in a CY threefold can be expressed
via Yoneda pairing $\left\langle \text{ \ },\text{ }\right\rangle $ of the
corresponding $Ext$ groups as follows:%
\[
\#\left(  \left(  \Sigma_{1},\theta_{1}\right)  ,\left(  \Sigma_{2},\theta
_{2}\right)  \right)  =\left\langle \lambda\left(  \Delta,\eta\right)
|_{\Sigma_{1}\times\Sigma_{2}},\mu\left(  \Sigma_{1}\times\Sigma_{2},\pi
_{1}^{\ast}\left(  \theta\right)  \wedge\pi_{2}^{\ast}\left(  \theta\right)
\right)  \right\rangle ,
\]
where $\pi_{1}$ and $\pi_{2}$ are the projections of M$\times$M to the first
and the second factor, respectively. The advantage of the latter
interpretation of the holomorphic linking is that it makes sense over an
arbitrary field and has universal algebro-geometric meaning. The relation with
the analytic expression follows from the interpretation of the Serre class
$\lambda(\Delta,\eta)$ as an element of $H^{2}($M$\times$M $-$ $\Delta
,\mathcal{O}_{\text{M}\times\text{M}-\Delta}),$ thanks to the Grothendieck
duality, and then from the identification of the local expression \ for
$\lambda(\Delta,\eta)$ as a Bochner-Martinelli kernel. Similarly, we obtain an
alternative expression for the holomorphic linking using the Grothendieck and
Serre classes for the pairs (X,Y$)=($M,$\Sigma_{1}),$ or $($M,$\Sigma_{2}%
\dot{)}.$

We also compare our homological formula for the holomorphic linking with a
similar expressions found by Atiyah in \cite{A} for the linking of the two
lines in a twistor space. In order to obtain the direct relation between our
holomorphic linking and Atiyah's, we extend our construction by considering
complex curves with marked points on them and restrict it to the case of two
spheres with two marked points.

The holomorphic linking studied in this paper has many predecessors. The
Lagrangian of the complex counterpart of CSW theory was first proposed by E.
Witten in \cite{W2}. He derived it as a low energy limit of an open string
theory and argued that the theory is finite in spite of the fact that it is
defined in six-dimensional space. We rediscovered the same theory following
the analogy with the real case, extensively studied in relation to
representation theory of loop groups. Other approaches to holomorphic linking
were previously considered by Atiyah \cite{A}, Penrose \cite{Pen} and
Gerasimov \cite{G}. In particular Gerasimov also suggested a path integral
presentation of the Atiyah linking of two lines in a twistor space, which is
similar to our path integral that arises in complexification of loop group
theory \cite{FK}. We expect that the combination of the different motivations
that lead to the holomorphic linking can be fruitful in the future
developments in this rich new field.

The organization of the paper by sections is as follows:

In \textbf{Section 2 }we recall a path integral derivation of the Gauss
linking number for CSW theory, previously obtained in \cite{P} and \cite{W1}.
The derivation yields, in particular, a well known explicit integral formula
for the Gauss invariant and relates it with the familiar geometric form of the
Gauss linking number.

In \textbf{Section 3} we repeat the path integral derivation in the complex
case and obtain a definition of the holomorphic linking. To make it rigorous,
we use the theory of currents on complex manifolds.

In \textbf{Section 4 }we use the existence of the Green current in Arakelov
geometry \cite{AB} for the pair M$\times$M and its diagonal $\Delta$ to recast
the definition of the holomorphic linking in a more invariant form. As a
result we also obtain a symmetric form of the holomorphic linking.

In \textbf{Section 5} we give explicit formulas for the restriction of the
Green kernel of the diagonal on an open affine set $U\times U$ in M$\times$M
by relating it to the Bochner-Martinelli kernel. This leads to an explicit
analytic formula for the holomorphic linking which is the analogue of the
classical Gauss formula.

In \textbf{Section 6 }we derive a geometric form of holomorphic linking. This
is a direct complex generalization of the usual topological form of the Gauss
linking number.

In \textbf{Section 7 }we introduce generalized notions of Grothendieck and
Serre classes. We prove their existence for any embedded submanifold Y into a
CY manifold M and the diagonal embedded in M$\times$M. Then we recast the
analytic formulas for the holomorphic linking in the universal language of
homological algebra.

In \textbf{Section 8} we give a definition of the holomorphic linking of two
Riemann surfaces with marked points on them. This allows us to define the
holomorphic linking of spheres embedded in a three dimensional complex
manifolds. As a consequence we are able to relate explicitly our notion of
holomorphic linking with the Atiyah holomorphic linking of rational curves in
the twistor space of the four dimensional sphere.

In conclusion we would like to mention some future perspectives that result
from our work. We have seen how different approaches to holomorphic linking
lead to equivalent definitions. However, our approach via complex CSW theory
has one important advantage, namely, it admits a non-abelian generalization at
least at formal level. To make it rigorous, one can try the pertributive
expansion of the complex CSW path integral. S. Donaldson and R. Thomas
obtained the analogue of Casson's invariant for CY manifolds in \cite{DT} ,
\cite{T} and \cite{T1}, which should appear as the first nontrivial term in
the pertributive expansion. In \cite{T} one can find further ideas how to
generalize the invariants obtained by S. Axelrod and I. Singer in \cite{AS}
and \cite{AS1}, coming from the expansion of the Chern-Simon functional in the
complex case. Another possible approach is to find a combinatorial calculus
for the complex CSW path integral and its generalizations. Interpretation of
non-abelian generalizations of the holomorphic linking in the context of local
cohomology can be an equally important application of complex CSW theory.
Finally we would like to note that our notion of holomorphic linking can be
easily generalized to submanifolds N$_{1}$ and N$_{2}$ in a n-dimensional CY
manifold M, where $\dim_{\mathbb{C}}$N$_{1}+\dim_{\mathbb{C}}$N$_{2}=n-1.$ In
this form it might be related to the height pairing for higher dimensional
cycles that was first constructed by Bloch and Beilinson, see \cite{B},
\cite{Be} and \cite{Bo}. Other constructions related to the notion of
holomorphic linking presented in this paper can be found in \cite{HR},
\cite{HR1} and \cite{T}.

\textbf{ACKNOWLEDGMENTS}

We are deeply obliged to A. Gerasimov for pointing out the work of Atiyah and
its relation to holomorphic linking and to S. Shatashvili for giving us
Gerasimov's letter addressed to him. We thank S. Axelrod, M. Khovanov, Yu.
Manin, G. Moore, E. Nisnevich, S. Shatashvili, R. Thomas, B. Wang and G.
Zuckerman for interesting discussions and helpful remarks. We thank B. Khesin
for participation in the beginning of this work and for his critical remarks.
In particular, the explicit geometric formula for the holomorphic linking was
his contribution (see \cite{HR} and \cite{HR1}). We are grateful to the Aspen
Institute for Physics for its hospitality in the summer of 1996 where we wrote
the first draft of this paper which has widely circulated since then in an
unfinished form.. We would like to thank S. Donaldson, M. Kontsevich, Yu.
Manin, E. Witten and S.-T. Yau for their encouragement to finish this paper.
Finally, the authors are grateful to both referees for suggesting some
corrections and for their thoughtful remarks. The research of the authors was
supported by NSF grants.

\section{Abelian Chern-Simons-Witten Theory and the Gauss Linking Number}

Let C$_{1}$ and C$_{2}$ be two non-intersecting knots in a simply connected
three dimensional real manifold M. The simplest nontrivial invariant of the
pair (C$_{1},$C$_{2})$ is the Gauss linking number denoted by \#(C$_{1}%
,$C$_{2})$. The abelian Chern-Simons-Witten theory gives a path integral
presentation of this invariant via a correlation function of holonomy
functionals attached to the knots. See \cite{P} and \cite{W1}. We will use
this presentation as a starting point of our approach and then show that it
leads to familiar analytic and geometric definitions of the Gauss linking
number. This point of view will allow us to produce complex counterparts of
all classical formulas.

Let $\mathcal{A}$ be the space of 1-forms on M and let $\mathcal{G}$ be the
subspace of exact 1-forms on M. The abelian Chern-Simons Lagrangian has the
form
\[
\mathcal{L}(A)=\int_{\text{M}}A\wedge dA,
\]
where $A\in\mathcal{A}$ and it is invariant with respect to translation of an
element from $\mathcal{G}$. One defines a holonomy functional on $\mathcal{A}%
$, called a Wilson loop, by
\[
W_{C}(A)=\exp\left(  \int_{C}A\right)
\]
for any knot $C$ in M. The Witten invariant of a link in M with components
$C_{1},...,C_{n}$ is given by the formal expression
\begin{equation}
\mathcal{Z(}\text{M;}C_{1},...,C_{n})=\int_{\mathcal{A}/\mathcal{G}%
}\mathcal{D}Ae^{\sqrt{-1}h\mathcal{L}(A)}W_{C_{1}}(A)...W_{C_{n}}(A),
\label{0}%
\end{equation}
where $h$ is a real parameter.

Though in the general case the rigorous definition of Witten's path integral
remains a widely open publicized problem, the explicit meaning of the above
definition can be described in the following way: Choose in $\mathcal{A}$ a
linear complement $\mathcal{A}_{0}$ to $\mathcal{G}$. Then for any
$\mathcal{G}$ invariant functional $\mathcal{F}$ on $\mathcal{A}$, one assumes
the following ''gauge fixing'':
\[
\int_{\mathcal{A}/\mathcal{G}}\mathcal{D}Ae^{\sqrt{-1}h\mathcal{L}%
(A)}\mathcal{F}(A)=\int_{\mathcal{A}_{0}}\mathcal{D}Ae^{\sqrt{-1}%
h\mathcal{L}(A)}\mathcal{F(}A)
\]
and the ''quasi-invariance'':
\[
\int_{\mathcal{A}_{0}}\mathcal{D}Ae^{\sqrt{-1}h\mathcal{L}(A+A_{1}%
)}\mathcal{F}(A+A_{1})=\int_{\mathcal{A}_{0}}\mathcal{D}Ae^{\sqrt
{-1}h\mathcal{L}(A)}\mathcal{F(}A),
\]
where $A_{1}\in\mathcal{A}_{0}.$

The Gauss linking number of $C_{1}$ and $C_{2}$ appears if we ''factor out''
the self-linking. Namely, one has the following identity, which we will use as
a definition:

\begin{definition}
\label{igor1}The Gauss linking number $\#(C_{1},C_{2})$ is defined by the
following formula:
\begin{equation}
\exp\left(  \frac{\sqrt{-1}}{2h}\#\left(  C_{1},C_{2}\right)  \right)
=\frac{\mathcal{Z(}\text{M;}C_{1},C_{2})\mathcal{Z(}\text{M)}}{\mathcal{Z(}%
\text{M;}C_{1})\mathcal{Z(}\text{M;}C_{2})}. \label{W0}%
\end{equation}

\end{definition}

We will show by formal application of the quasi-invariance that Definition
\ref{igor1} yields a certain well-defined functional on C$_{1}$ and C$_{2}$
which does not depend on the choice of the \textquotedblright gauge
fixing\textquotedblright. Let us choose two forms $\omega_{1}$ and $\omega
_{2}$ such that
\[
\int_{\text{M}}A\wedge\omega_{i}=\int_{C_{i}}A\text{ for }i=1,2.
\]
In fact, the forms $\omega_{1}$ and $\omega_{2}$ are currents, i.e. linear
functionals on the space of smooth 1-forms. However, by duality one can define
a differential on the space of currents. (see \cite{GH} and section
\ref{curr}.) Then it follows that $\omega_{i}$ are exact for $i=1,$ $2.$ We
also obtain:
\[
\mathcal{L}\left(  A+\frac{\sqrt{-1}}{2h}A_{1}+\frac{\sqrt{-1}}{2h}%
A_{2}\right)  =
\]%
\[
\mathcal{L}(A)+\frac{-1}{4h^{2}}\mathcal{L(}A_{1})+\frac{-1}{4h^{2}%
}\mathcal{L}(A_{2})+
\]%
\[
\frac{\sqrt{-1}}{h}\left(  \int_{\text{M}}A\wedge dA_{1}+\int_{\text{M}%
}A\wedge dA_{2}\right)  -\frac{1}{4h^{2}}\int_{\text{M}}A_{1}\wedge dA_{2}.
\]
Let us choose the real currents $A_{1}$ and $A_{2}$ satisfying the following
relations
\[
dA_{i}=\omega_{i}\text{ for }i=1,2.
\]
The currents $A_{i}$ should be viewed as generalized one forms and we choose
them in a unique way by requiring that they belong to the completion of
$\mathcal{A}_{0}$. We will also denote by
\[
d^{-1}\omega_{i}=A_{i}%
\]
for $i=1,$ 2. With the above choice of $A_{1}$ and $A_{2}$ the integrand in
$\mathcal{Z}($M;C$_{1},C_{2})$ becomes constant and the right hand side of
\ the expression in Definition \ref{igor1} transforms to the following form:
\[
\exp\left(  \frac{\sqrt{-1}}{2h}\int_{\text{M}}A_{1}\wedge dA_{2}\right)
=\exp\left(  \frac{\sqrt{-1}}{2h}\int_{\text{M}}d^{-1}\omega_{1}\wedge
\omega_{2}\right)  .
\]
It is easy to check that the integral in the exponent does not depend on the
choice of the complement $\mathcal{A}_{0}$ and is well defined. Thus,
Definition \ref{igor1} assumes the following form:
\[
\#(C_{1},C_{2})=\int_{\text{M}}A_{1}\wedge dA_{2}.
\]
In order to relate it with the standard expression of the classical Gauss
linking number we will rewrite the above integral in a different form. Let us
denote by $\delta_{\Delta}$ a Dirac current concentrated on the diagonal
$\Delta\subset$M$\times$M, namely for any smooth 3-form $\xi$ on M one has%
\[
\int_{\text{M}\times\text{M}}\pi_{i}^{\ast}(\xi)\wedge\delta_{\Delta}%
=\int_{\Delta}\pi_{i}^{\ast}(\xi)=\int_{\text{M}}\xi,\text{ }i=1,2
\]
where $\pi_{1}$ and $\pi_{2}$ are the projections on the first and the second
factor, respectively.

It is well known (see e.g. \cite{GH}, Chapter 3) that on a compact manifold M
the cohomology theory defined by the currents coincides with the De Rham
cohomology. We will recall the basic definitions related to currents in
Section 3.2 below. Let us choose a closed smooth form $\omega_{\Delta}$ which
represents the same cohomology class as the Dirac current $\delta_{\Delta}.$
Then there exists a current S such that
\begin{equation}
dS+\delta_{\Delta}=\omega_{\Delta}. \label{2}%
\end{equation}
The equation $\left(  \ref{2}\right)  $ implies that the singular support of
the current S is exactly the diagonal $\Delta,$ so the restriction of S to
M$\times$M $-$ $\Delta$ is a smooth form, which we denote by S$_{\Delta}.$

Since the cohomology class of $\delta_{\Delta},$ and therefore of
$\omega_{\Delta},$ is the Poincare dual class of $\Delta,$ the restriction of
the closed smooth form $\omega_{\Delta}$ to M$\times$M $-$ $\Delta$ is an
exact smooth form, i.e.
\[
\omega_{\Delta}|_{\text{M}\times\text{M }-\text{ }\Delta}=d\psi_{\Delta}%
\]
for some smooth form $\psi_{\Delta}.$ Then we introduce a smooth two-form on
M$\times$M $-$ $\Delta$%
\begin{equation}
\phi_{\Delta}=\psi_{\Delta}-\text{S}_{\Delta}. \label{3}%
\end{equation}
The cohomology class of $\phi_{\Delta}$ on M$\times$M $-$ $\Delta$ does not
depend on the choices involved in the above construction.

Now we can rewrite%
\[%
%TCIMACRO{\dint \limits_{\text{M}}}%
%BeginExpansion
{\displaystyle\int\limits_{\text{M}}}
%EndExpansion
d^{-1}\omega_{1}\wedge\omega_{2}=%
%TCIMACRO{\dint \limits_{\text{M}\times\text{M}}}%
%BeginExpansion
{\displaystyle\int\limits_{\text{M}\times\text{M}}}
%EndExpansion
\pi_{1}^{\ast}\left(  d^{-1}\omega_{1}\right)  \wedge\pi_{2}^{\ast}\left(
\omega_{2}\right)  \wedge\delta_{\Delta}=
\]%
\[%
%TCIMACRO{\dint \limits_{\text{M}\times\text{M}}}%
%BeginExpansion
{\displaystyle\int\limits_{\text{M}\times\text{M}}}
%EndExpansion
\pi_{1}^{\ast}\left(  d^{-1}\omega_{1}\right)  \wedge\pi_{2}^{\ast}\left(
\omega_{2}\right)  \wedge\left(  \omega_{\Delta}-d\text{S}\right)  =
\]%
\[%
%TCIMACRO{\dint \limits_{\text{M}\times\text{M$-$}\Delta}}%
%BeginExpansion
{\displaystyle\int\limits_{\text{M}\times\text{M$-$}\Delta}}
%EndExpansion
\pi_{1}^{\ast}\left(  d^{-1}\omega_{1}\right)  \wedge\pi_{2}^{\ast}\left(
\omega_{2}\right)  \wedge d\psi_{\Delta}-%
%TCIMACRO{\dint \limits_{\text{M}\times\text{M}}}%
%BeginExpansion
{\displaystyle\int\limits_{\text{M}\times\text{M}}}
%EndExpansion
\pi_{1}^{\ast}\left(  d^{-1}\omega_{1}\right)  \wedge\pi_{2}^{\ast}\left(
\omega_{2}\right)  \wedge d\text{S}=
\]%
\[%
%TCIMACRO{\dint \limits_{\text{M}\times\text{M $-$ }\Delta}}%
%BeginExpansion
{\displaystyle\int\limits_{\text{M}\times\text{M $-$ }\Delta}}
%EndExpansion
\pi_{1}^{\ast}\left(  \omega_{1}\right)  \wedge\pi_{2}^{\ast}\left(
\omega_{2}\right)  \wedge\psi_{\Delta}-%
%TCIMACRO{\dint \limits_{\text{M}\times\text{M $-$ }\Delta}}%
%BeginExpansion
{\displaystyle\int\limits_{\text{M}\times\text{M $-$ }\Delta}}
%EndExpansion
\pi_{1}^{\ast}\left(  \omega_{1}\right)  \wedge\pi_{2}^{\ast}\left(
\omega_{2}\right)  \wedge\text{S}_{\Delta}=
\]%
\begin{equation}%
%TCIMACRO{\dint \limits_{\text{M}\times\text{M $-$ }\Delta}}%
%BeginExpansion
{\displaystyle\int\limits_{\text{M}\times\text{M $-$ }\Delta}}
%EndExpansion
\pi_{1}^{\ast}\left(  \omega_{1}\right)  \wedge\pi_{2}^{\ast}\left(
\omega_{2}\right)  \wedge\phi_{\Delta}=%
%TCIMACRO{\dint \limits_{\text{C}_{1}\times\text{C}_{2}}}%
%BeginExpansion
{\displaystyle\int\limits_{\text{C}_{1}\times\text{C}_{2}}}
%EndExpansion
\phi_{\Delta}. \label{4}%
\end{equation}
Let M$=S^{3}=\mathbb{R}^{3}\cup\infty.$ Standard facts about Green currents
imply that the form $\phi_{\Delta}$ is given on $\mathbb{R}^{3}\times
\mathbb{R}^{3}-\Delta$ by
\begin{equation}
\phi_{\Delta}=\frac{1}{4\pi}\text{ }\frac{\varepsilon^{ijk}(x_{i}-y_{i}%
)}{|x-y|^{3}}dx_{j}\wedge dy_{k}, \label{5}%
\end{equation}
where $x,y\in\mathbb{R}^{3},$ $\varepsilon^{ijk}$ is the sign of permutation
$(i,j,k)$ and we assume the summation over all three indexes. Thus $\left(
\ref{4}\right)  $ and $\left(  \ref{5}\right)  $ imply the classical Gauss
formula for the linking number%
\begin{equation}
\#\left(  \text{C}_{1},\text{C}_{2}\right)  =\frac{1}{4\pi}\text{ }%
%TCIMACRO{\dint \limits_{\text{C}_{1}\times\text{C}_{2}}}%
%BeginExpansion
{\displaystyle\int\limits_{\text{C}_{1}\times\text{C}_{2}}}
%EndExpansion
\frac{\varepsilon^{ijk}(x_{i}-y_{i})}{|x-y|^{3}}dx_{j}\wedge dy_{k}. \label{6}%
\end{equation}

One can also derive a topological interpretation of the Gauss linking number
by considering a disk $D_{1}$ in M with a boundary $\partial D_{1}=C_{1}$ and
"reversing" the steps in $\left(  \ref{4}\right)  $:%
\[%
%TCIMACRO{\dint \limits_{\text{C}_{1}\times\text{C}_{2}}}%
%BeginExpansion
{\displaystyle\int\limits_{\text{C}_{1}\times\text{C}_{2}}}
%EndExpansion
\phi_{\Delta}=%
%TCIMACRO{\dint \limits_{\text{C}_{1}\times\text{C}_{2}}}%
%BeginExpansion
{\displaystyle\int\limits_{\text{C}_{1}\times\text{C}_{2}}}
%EndExpansion
\left(  \psi_{\Delta}-S_{\Delta}\right)  =
\]%
\[%
%TCIMACRO{\dint \limits_{\text{D}_{1}\times\text{C}_{2}\text{ $-$ }\Delta}}%
%BeginExpansion
{\displaystyle\int\limits_{\text{D}_{1}\times\text{C}_{2}\text{ $-$ }\Delta}}
%EndExpansion
d\psi_{\Delta}-%
%TCIMACRO{\dint \limits_{\text{D}_{1}\times\text{C}_{2}}}%
%BeginExpansion
{\displaystyle\int\limits_{\text{D}_{1}\times\text{C}_{2}}}
%EndExpansion
dS=
\]%
\begin{equation}%
%TCIMACRO{\dint \limits_{\text{D}_{1}\times\text{C}_{2}\text{ $-$ }\Delta}}%
%BeginExpansion
{\displaystyle\int\limits_{\text{D}_{1}\times\text{C}_{2}\text{ $-$ }\Delta}}
%EndExpansion
\omega_{\Delta}-%
%TCIMACRO{\dint \limits_{\text{D}_{1}\times\text{C}_{2}}}%
%BeginExpansion
{\displaystyle\int\limits_{\text{D}_{1}\times\text{C}_{2}}}
%EndExpansion
dS=%
%TCIMACRO{\dint \limits_{\text{D}_{1}\times\text{C}_{2}\text{ }}}%
%BeginExpansion
{\displaystyle\int\limits_{\text{D}_{1}\times\text{C}_{2}\text{ }}}
%EndExpansion
\delta_{\Delta}. \label{7}%
\end{equation}
The latter integral should be viewed as a pairing between the cohomology class
represented by the current $\delta_{\Delta}$ and the cycle D$_{1}\times C_{2}$
(See Section 3.2. for details.) It counts the number of the intersection
points of $D_{1}$ and $C_{2}$ with a sign that depends on the orientation.
This presentation provides a simple topological meaning of the Gauss linking
number and implies its integrality properties.

Finally we note that various forms and integrals that appeared in the above
discussion of Gauss linking number admit a natural cohomological
interpretation. The Gauss linking number can be obtained as a pairing of
cohomology classes with the appropriate cycles.

In the rest of the paper we will produce complex counterparts of various
manifestation of the Gauss linking number.

\section{Complex Counterpart of the Abelian CSW Theory and Holomorphic
Linking}

The complexification of loop group theory developed in \cite{EF} and \cite{FK}
yields in particular a complex analogue of the abelian Chern-Simons-Witten
Lagrangian and Wilson loop functional previously considered in the context of
string theory. \cite{W2}. This leads to a path integral definition of the
holomorphic linking of Riemann surfaces in a CY threefold. We show that the
formal calculation of the path integral gives rise to a rigorous mathematical
notion of holomorphic linking which can be viewed as a complex counterpart of
the Gauss linking number.

\subsection{Path Integral Definition of the Holomorphic Linking}

Let M be a K\"{a}hler manifold of complex dimension $n.$ We recall that M is a
CY manifold if it admits a metric with holonomy group $SU(n).$ It is a well
known fact that this definition of a CY manifold M implies that there exists a
unique up to constant holomorphic n-form $\eta$ without zeroes. See
\cite{besse}. Let M be a CY threefold and $\eta$ be a holomorphic three form
on M. Until Section \ref{At} we assume that M is a compact variety.

The Lagrangian of the complex abelian CSW theory is defined as follows. Let
$\mathcal{A}^{0,1}$ denote the complex space of (0,1) forms on M and let
$\mathcal{G}^{0,1}$ be the subspace of $\overline{\partial}$-exact forms. Then
we set%
\[
\mathcal{L}(A):=\int_{\text{M}}A\wedge\overline{\partial}A\wedge\eta
,\ A\in\mathcal{A}^{0,1}.
\]
It follows immediately from Stokes' Theorem that the Lagrangian is invariant
with respect to $\mathcal{G}^{0,1},$ i.e.
\[
\mathcal{L}(A+\overline{\partial}\phi)=\mathcal{L}(A)
\]
for any function $\phi.$

For a Riemann surface $\Sigma$, together with a holomorphic 1-form $\theta$ we
define an analogue of the Wilson loop functional
\begin{equation}
W_{(\Sigma,\theta)}(A):=\exp\left(  \int_{\Sigma}A\wedge\theta\right)  ,\text{
}A\in\mathcal{A}^{0,1}. \label{W}%
\end{equation}
The expression in $\left(  \ref{W}\right)  $ is also invariant with respect to
$\mathcal{G}^{0,1}.$ We will define the complex version of $\left(
\ref{0}\right)  $ as follows:
\[
\mathcal{Z(}M;(\Sigma_{1},\theta_{1}),...,(\Sigma_{n},\theta_{n}%
)):=\int_{\mathcal{A}^{0,1}/\mathcal{G}^{0,1}}\mathcal{D}A\exp\left(
\sqrt{-1}h\mathcal{L}(A)\right)  \prod_{i=1}^{n}W_{(\Sigma_{i},\theta_{i}%
)}(A),
\]
where as in the real case, the path integral should satisfy the gauge fixing
and quasi-invariance. Following the analogy with the real case we now define a
complex counterpart of the Gauss linking number.

\begin{definition}
\label{igor2}The holomorphic linking of \ two Riemann surfaces $\Sigma_{1}$
and $\Sigma_{2}$ with chosen holomorphic forms $\theta_{1}$ and $\theta_{2},$
respectively, is defined by the following formula:
\begin{equation}
\exp\left(  \frac{\sqrt{-1}}{2h}\#\left(  (\Sigma_{1},\theta_{1}),(\Sigma
_{2},\theta_{2})\right)  \right)  =\frac{\mathcal{Z(}\text{M};\text{ }%
(\Sigma_{1},\theta_{1}),(\Sigma_{2},\theta_{2}))\mathcal{Z(}\text{M}%
\mathcal{)}}{\mathcal{Z(}\text{M};\text{ }(\Sigma_{1},\theta_{1}%
))\mathcal{Z(}\text{M};\text{ }(\Sigma_{2},\theta_{2}))}. \label{W1}%
\end{equation}

\end{definition}

Let $\mathcal{H}^{1}(\Sigma)$ denote the space of holomorphic differentials on
$\Sigma.$ We will show that for any two Riemann surfaces $\Sigma_{1}$ and
$\Sigma_{2}$ the formal expression in Definition \ref{igor2} gives rise to a
well defined linear map $\#$, which we will call the holomorphic linking:
\[
\#\left(  (\Sigma_{1},\circ),(\Sigma_{2},\circ)\right)  :\mathcal{H}%
^{1}(\Sigma_{1})\times\mathcal{H}^{1}(\Sigma_{2})\rightarrow\mathbb{C}.
\]

\subsection{\label{curr}Currents on Complex Manifolds}

In order to transform the formal definition of the holomorphic linking to a
rigorous form, we will need to recall some basic facts about currents on a
complex manifold X. For more details see \cite{GH}.

Let X be a compact K\"{a}hler complex manifold of complex dimension n. We will
denote by $A^{p,q}(X)$ the vector space of $C^{\infty}$ complex valued forms
of type (p,q). The space $A^{m}(X)$ of complex valued $C^{\infty}$ m- forms is
given by
\[
A^{m}(X)=\underset{p+q=m}{\oplus}A^{p,q}(X).
\]
We denote by
\[
\partial:A^{p,q}(X)\rightarrow A^{p+1,q}(X)
\]
and
\[
\overline{\partial}:A^{p,q}(X)\rightarrow A^{p,q+1}(X)
\]
the Dolbeault differentials.

Let $D_{m}$ denote the dual space of $A^{m}(X)$ with respect to the standard
Frechet topology. We also denote by $D_{p,q}(X)$ the dual space of
$A^{p,q}(X).$ We define $D^{p,q}(X)$ as $D_{n-p,n-q}(X).$ For any X there is a
natural inclusion
\[
A^{p,q}(X)\subset D^{p,q}(X).
\]
In fact, to each form $\omega\in A^{p,q}(X)$ we can associate a continuous
functional on $A^{n-p,n-q}(X)$ by the formula:
\[
\left\langle \omega,\alpha\right\rangle =\omega(\alpha):=\int_{X}\omega
\wedge\alpha
\]
for any $\alpha\in A^{n-p,n-q}(X).$

The differentials $d,$ $\partial$ and $\overline{\partial}$ act on the spaces
$D^{m}(X)$ and $D^{p,q}(X)$ in a natural way by extending the action on the
subspaces $A^{m}(X)$ and $A^{p,q}(X),$ respectively. From Stokes' Theorem we
have:
\[
\left\langle \omega,d\alpha\right\rangle =\int_{X}\omega\wedge d\alpha
=(-1)^{m-1}\int_{X}d\omega\wedge\alpha
\]
for any $\alpha\in A^{2n-m}(X)$ and $\omega\in A^{m}(X).$ So we define
$d\omega$ for any $\omega\in D^{m}(X)$ as follows
\[
\left\langle d\omega,\alpha\right\rangle :=(-1)^{m-1}\left\langle
\omega,d\alpha\right\rangle
\]
for any $\alpha\in A^{2n-m-1}.$ In the same way we start with the actions of
$\partial$ and $\overline{\partial}$ on the everywhere dense subspaces
\[
A^{p,q}(X)\subset D^{p,q}(X)
\]
and then we define $\partial\omega$ and $\overline{\partial}\omega$
respectively for any $\omega\in D^{p-1,q}(X)$ or $\omega\in D^{p,q-1}(X)$
\[
\left\langle \partial\omega,\alpha\right\rangle :=(-1)^{p+q-1}\left\langle
\omega,\partial\alpha\right\rangle
\]
and
\[
\left\langle \overline{\partial}\omega,\alpha\right\rangle :=(-1)^{p+q-1}%
\left\langle \omega,\overline{\partial}\alpha\right\rangle
\]
for any $\alpha\in A^{n-p,n-q}.$ A current $\omega\in D^{p,q}(X)$ is
$\overline{\partial}$ closed if and only if for any $\overline{\partial}$
exact form $\overline{\partial}\alpha\in$ $A^{n-p,n-q}(X)$ one has
\[
\left\langle \omega,\overline{\partial}\alpha\right\rangle =0.
\]
Similar characterization is valid for $\partial$-closed currents.

\begin{definition}
\label{andrey48} \textbf{i. }Let Y be a complex subvariety in a projective
variety X$.$ We define a current $\delta_{\text{Y}}$ corresponding to Y via
the integration pairing:
\[
\left\langle \delta_{\text{Y}},\omega\right\rangle =\int_{\text{Y}}%
\omega|_{\text{Y}},
\]
where $\omega$ is a smooth form on X. In particular, \textit{the Dirac kernel
}$\delta_{\Delta}$\textit{ of the diagonal embedding of M is a current of type
}$(n,n)$ \textit{on M}$\times$\textit{M, }$\dim_{\mathbb{C}}$M$=n,$\textit{
such that }$\delta_{\Delta}$\textit{ is zero on M}$\times$\textit{M $-$
}$\Delta$\textit{ and }
\[
\left\langle \delta_{\Delta},\alpha\right\rangle =\int_{\Delta\subset
\text{M}\times\text{M}}\alpha|_{\Delta}%
\]
\textit{for any smooth form }$\alpha$\textit{\ of type }$(n,n)$\textit{\ on
}M$\times$M.\textit{\ }
\end{definition}

\textbf{ii. }\textit{Let X and Y be as above, such that X is CY, }%
$\dim_{\mathbb{C}}$X$=n,\dim_{\mathbb{C}}$Y$=n-m.$\textit{ We define the Dirac
antiholomorphic current }$\overline{\theta}_{\text{Y}}$ \textit{of type
}$(0,m)$ \textit{on X corresponding to} \textit{the holomorphic }$n-m$
\textit{form} $\theta_{\text{Y}}$ \textit{on the subvariety Y as follows: Let}
$\beta$ \textit{any smooth form of type }$(n,n-m)$ \textit{on X, then }%
\[
\left\langle \overline{\theta}_{\text{Y}},\beta\right\rangle =%
%TCIMACRO{\dint \limits_{\text{Y}}}%
%BeginExpansion
{\displaystyle\int\limits_{\text{Y}}}
%EndExpansion
\theta\wedge\frac{\beta}{\eta}.
\]
\textit{The holomorphic Dirac kernel }$\theta_{\text{Y}}$ \textit{is defined
in a similar manner. In particular the antiholomorphic Dirac kernel
}$\overline{\eta}_{\Delta}\mathit{\ }$\textit{is a current of type }%
$(0,n)$\textit{ on M}$\times$\textit{M such that}:
\begin{equation}
\left\langle \overline{\eta}_{\Delta},\beta\right\rangle =\int_{\Delta}%
\frac{\beta}{\pi_{1}^{\ast}(\eta)}=\int_{\Delta}\frac{\beta}{\pi_{2}^{\ast
}(\eta)} \label{LP2}%
\end{equation}
\textit{for any C}$^{\infty}$\textit{ form }$\beta$ \textit{of type }$(2n,n)$
\textit{on M}$\times$\textit{M, where }$\pi_{1}$\textit{\ and }$\pi_{2}$
\textit{are the\ projections\ of }M$\times$M\textit{\ on the first and the
second factor, respectively. Since the current }$\overline{\eta}_{\Delta}%
$\textit{ is supported on the diagonal its definition does not depend on the
choice of the projection. }

For any $p$ dimensional complex subvariety Y$\subset$X we define the Green
current $g_{\text{Y}}\in D^{p-1,p-1}($X$)$ such that%

\begin{equation}
\partial\overline{\partial}g_{Y}+\delta_{\text{Y}}=\omega_{\text{Y}},
\label{LP}%
\end{equation}
where $\omega_{Y}$ is some closed $C^{\infty}$ $(p,p)$ form representing the
Poincare dual class of [Y]$.$ The existence of the Green current with only
logarithmic singularities is one of the key result in Arakelov geometry (see
\cite{AB}).

We will use two basic facts from the theory of currents. We can define both
Dolbeault's and De Rham cohomology using currents\ instead of $C^{\infty}$
forms. The basic theorem, due to Whitney, asserts that on compact smooth
manifolds De Rham or, in the complex category, Dolbeault's cohomology are
isomorphic to De Rham or Dolbeault's cohomology obtained from currents. (See
Chapter 5 in \cite{GH}.) Also according to H\"{o}rmander \cite{Hor} we can
define the exterior product $\alpha\wedge\beta$ of two currents as another
current, if the singular supports of the currents $\alpha$ and $\beta$ are
disjoint as sets.

Since currents represent De Rham and Dolbeault cohomology classes thanks to
the theorem of Whitney we will sometimes denote their pairing with homology
classes by the integral, e.g.
\[%
%TCIMACRO{\dint \limits_{\text{Y}^{^{\prime}}}}%
%BeginExpansion
{\displaystyle\int\limits_{\text{Y}^{^{\prime}}}}
%EndExpansion
\delta_{\text{Y}}%
\]
where Y$^{^{\prime}}\subset$X is another subvariety in X of complimentary
complex dimension to Y. In particular, that is the precise meaning of the
formula $\left(  \ref{7}\right)  $ for the Gauss linking number.

\begin{proposition}
\label{igor7}For a pair $(\Sigma,\theta),$ where $\Sigma$ is a Riemannian
surface in the three dimensional CY manifolds M and $\theta$ is a holomorphic
form on $\Sigma,$ the antiholomorphic Dirac current $\overline{\theta}%
_{\Sigma}$ is $\overline{\partial}$ exact current of type $(0,2)$ and
\begin{equation}
\overline{\theta}_{\Sigma}=\overline{\partial}A, \label{3.2.II}%
\end{equation}
for some current $A$ of type (0,1) with a singular support on $\Sigma.$
\end{proposition}

\textbf{Proof:} Since M is a CY manifold then any closed $(0,1)$ form $\beta$
is an exact form. Let $\beta=$ $\overline{\partial}f\wedge\eta$ on M$.$ Then
Stokes' theorem implies
\[
\left\langle \overline{\theta}_{\Sigma},\overline{\partial}f\wedge
\eta\right\rangle =\int_{\Sigma}\theta\wedge\overline{\partial}f=\int_{\Sigma
}d\left(  f\text{ }\theta\right)  =0.
\]
Therefore the antiholomorphic Dirac current $\overline{\theta}_{\Sigma}$ is
$d$ and $\overline{\partial}$ closed. From the theorem of Whitney, and \ since
$H^{2}($M$,\Omega^{3})=0$ for a CY manifold, we obtain that the
antiholomorphic Dirac current $\overline{\theta}_{\Sigma}$ is a $\overline
{\partial}$ exact current of type $(0,2)$ on M. Thus $\left(  \ref{3.2.II}%
\right)  $ holds for some current $A$ of type (0,1). It follows from
Definition \ref{andrey48} of the antiholomorphic Dirac current $\overline
{\theta}_{\Sigma}$ and $\left(  \ref{3.2.II}\right)  $ that the singular
supports of the currents $A$ and $\overline{\partial}A$ are on $\Sigma.$
Proposition \ref{igor7} is proved. $\blacksquare$

\subsection{Definition and Properties of Holomorphic Linking}

We will use some formal path integral computations to find explicit formula
for the expression in Definition \ref{igor2} in the complex case following the
exposition in Section 2. We define $\mathcal{A}^{0,1}$ as the space of (0,1)
forms on M and $\mathcal{G}^{0,1}$ as the space of $\overline{\partial}$
closed (0,1) forms on M. As in the real case we choose a complement to
$\mathcal{G}^{0,1}$ which we call $\mathcal{A}_{0}^{0,1}.$ For any invariant
functional $F$ on $\mathcal{A}^{0,1}$ with respect to $\mathcal{G}^{0,1}$ we
assume that%

\[
\mathbf{a.}\int_{\mathcal{A}^{0,1}/\mathcal{G}^{0,1}}\mathcal{D}A\exp\left(
\sqrt{-1}h\mathcal{L(}A)\right)  F(A)=\int_{\mathcal{A}_{0}^{0,1}}%
\mathcal{D}A\exp\left(  \sqrt{-1}h\mathcal{L(}A)\right)  F(A).
\]
and for $A_{1}\in\mathcal{A}_{0}^{0,1}$
\[
\mathbf{\ b.}\int_{\mathcal{A}_{0}^{0,1}}\mathcal{D}A\exp\left(  \sqrt
{-1}h\mathcal{L(}A+A_{1})\right)  F(A+A_{1})=\int_{\mathcal{A}_{0}^{0,1}%
}\mathcal{D}A\exp\left(  \sqrt{-1}h\mathcal{L(}A)\right)  F(A).
\]

The conditions \textbf{a} and \textbf{b} imply that we have the following
expression for the formula $\left(  \ref{W1}\right)  $:
\[
\frac{\mathcal{Z(}\text{M,}\left(  \Sigma_{1},\theta_{1}\right)  ,\left(
\Sigma_{2},\theta_{2}\right)  )\mathcal{Z(}\text{M)}}{\mathcal{Z(}%
\text{M,}\left(  \Sigma_{1},\theta_{1}\right)  )\mathcal{Z(}\text{M,}\left(
\Sigma_{2},\theta_{2}\right)  )}=\exp\left(  \frac{\sqrt{-1}}{2h}%
\int_{\text{M}}A_{1}\wedge\overline{\partial}A_{2}\wedge\eta\right)  ,
\]
where $A_{1}$ and $A_{2}$ are forms in $\mathcal{A}_{0}^{0,1}$defined in
Proposition \ref{igor7} corresponding to $\Sigma_{1}$ and $\Sigma_{2},$
respectively. Now we can give a rigorous definition of the holomorphic linking.

\begin{definition}
\label{igor8} The holomorphic linking between two curves $\Sigma_{1}$ and
$\Sigma_{2}$ with two holomorphic 1-forms $\theta_{1}$ and $\theta_{2}$ on
them is defined by the formula:
\begin{equation}
\#\left(  \left(  \Sigma_{1},\theta_{1}\right)  ,\left(  \Sigma_{2},\theta
_{2}\right)  \right)  =\int_{\text{M}}A_{1}\wedge\overline{\partial}%
A_{2}\wedge\eta, \label{W3}%
\end{equation}
\textit{where }$A_{1}$\textit{\ and }$A_{2}$\textit{\ are currents in
$\mathcal{A}_{0}^{0,1}$ }defined by $\left(  \ref{3.2.II}\right)  $ \textit{in
Proposition }\ref{igor7} for the pairs $\left(  \Sigma_{1},\theta_{1}\right)
$ and $\left(  \Sigma_{2},\theta_{2}\right)  $, respectively.
\end{definition}

The integral in $\left(  \ref{W3}\right)  $ makes sense since $A_{1}$
and\textit{\ }$\overline{\partial}A_{2}$ have disjoint supports. This follows
from the definitions of $A_{1}$ and $A_{2}$ and since\textit{\ }$\Sigma
_{1}\cap\Sigma_{2}=\emptyset$ as in the real case$.$ Substituting the
expression for $A_{1}$ and $\overline{\partial}A_{2}$ from Proposition
\ref{igor7} and the Definition \ref{andrey48} of antiholomorphic currents
$\left(  \overline{\theta}_{i}\right)  _{\Sigma_{i}}$ we obtain the first
explicit expression for the holomorphic linking%
\[
\#\left(  \left(  \Sigma_{1},\theta_{1}\right)  ,\left(  \Sigma_{2},\theta
_{2}\right)  \right)  =
\]%
\begin{equation}
\int_{\text{M}}\overline{\partial}^{-1}\left(  \left(  \overline{\theta}%
_{1}\right)  _{\Sigma_{1}}\right)  \wedge\left(  \overline{\theta}_{2}\right)
_{\Sigma_{2}}\wedge\eta=\int_{\Sigma_{2}}\left(  \overline{\partial}%
^{-1}\left(  \left(  \overline{\theta}_{1}\right)  _{\Sigma_{1}}\right)
\right)  |_{\Sigma_{2}}\wedge\theta_{2}. \label{W3a}%
\end{equation}
It is easy to check the following properties:

\begin{proposition}
\label{gr}\textbf{i. }The holomorphic linking does not depend on the choice of
the complement $\mathcal{A}_{0}^{0,1}.$
\end{proposition}

\textbf{ii.} \textit{The holomorphic linking is symmetric:}\textbf{\ }
\[
\#\left(  \left(  \Sigma_{1},\theta_{1}\right)  ,\left(  \Sigma_{2},\theta
_{2}\right)  \right)  =\#\left(  \left(  \Sigma_{2},\theta_{2}\right)
,\left(  \Sigma_{1},\theta_{1}\right)  \right)  .
\]

\textbf{iii. }$\#\left(  \left(  \Sigma_{1},\theta_{1}\right)  ,\left(
\Sigma_{2},\theta_{2}\right)  \right)  $\textit{\ is linear in } $\theta_{1}$
\textit{and }$\theta_{2}.$

\section{Green Kernel and Symmetric Form of Holomorphic Linking}

Now that we have a rigorous definition of the holomorphic linking, we might
try to express it in a more invariant form following the analogy with the real
case. In Section 2 we derived the Gauss kernel for the linking number using
the Green kernel for the exterior derivative operator. In order to obtain a
similar formula for the holomorphic linking, we first show the existence of
the analogue of the Green kernel of the operator $\overline{\partial}$, which
can be viewed as a complexification of the exterior derivative in the real case.

\subsection{The Cohomology Groups of Zariski Open Varieties}

Let $X_{0}=%
%TCIMACRO{\dbigcup \limits_{i=1}^{N}}%
%BeginExpansion
{\displaystyle\bigcup\limits_{i=1}^{N}}
%EndExpansion
C_{i}$ be a divisor with normal crossings in a projectve variety $X.$The
groups of the cohomology of the variety $X-X_{0}$ can be computed as the
cohomology groups of the de Rham log complex $\mathcal{A}^{\ast}%
(X,\log\left\langle X_{0}\right\rangle )$.

\begin{definition}
\label{d1}\textbf{i. }We will say that a form $\omega$ on one of the
components $C_{i_{0}}$ of $X_{0}$ has a logarithmic singularities if for each
point
\[
z\in C_{i_{0}}\cap...\cap C_{i_{k}}%
\]
and some open neighborhood $U\subset X$ of the point $z$ we have%
\begin{equation}
\omega|_{U}=\alpha\wedge\frac{dz_{i_{1}}}{z_{i_{1}}}\wedge...\wedge
\frac{dz_{i_{k}}}{z_{i_{k}}}, \label{log}%
\end{equation}
\textit{where } $\alpha$ \textit{is a }C$^{\infty}\,$\ \textit{form in }$U$
and $z_{i_{0}}\times...\times z_{i_{k}}\omega$ and $z_{i_{0}}\times...\times
z_{i_{k}}d\omega$ are smooth forms in $U$.\textit{ \ }$X_{0}\cap U$ \textit{is
defined by the equations }
\[
z_{i_{1}}\times...\times z_{i_{k}}=0
\]
in $U.$
\end{definition}

\textbf{ii}\textit{. We define the de Rham log complex with the standard
differential as follows:}%
\[
\mathcal{A}^{\ast}(X_{0},\log\left\langle X_{0}\right\rangle )
\]%
\[
=\{\omega\in C^{\infty}\left(  X-\text{ }X_{0},\Omega^{\ast}\right)
|\omega\text{ \textit{and}}\mathit{\ }d\omega\text{ \textit{are} }C^{\infty
}\text{\textit{forms}}\mathit{\ }%
\]%
\begin{equation}
\text{\textit{on}}\mathit{\ }X-\text{ }X_{0}\text{ \textit{which have}%
}\mathit{\ \log\ }\text{\textit{singularities along} }X_{0}\}. \label{LOG}%
\end{equation}

By using the Poincare residue map Deligne proved the following:

\begin{theorem}
\label{D0}The cohomology of $X$ $-$ $X_{0}$ are equal to the cohomology of the
De Rham log complex $\mathcal{A}^{\ast}(X_{0},\log\left\langle X_{0}%
\right\rangle ).$
\end{theorem}

For the proof of Theorem \ref{D0} see \cite{GS}.

\subsection{Basic Definitions and The Existence Theorem}

Let M be a CY of complex dimension n. We will normalize the holomorphic n-form
$\eta$ as follows:
\begin{equation}
(-1)^{\frac{n(n-1)}{2}}(-i)^{n}\int_{\text{M}}\eta\wedge\overline{\eta}=1.
\label{14}%
\end{equation}
According to the Theorem of Gillet and Soule stated in Section\textbf{
}3.2,\textbf{ }for each closed form of type $(n,n)$ that represents the
Poincare dual of the homology class of the diagonal $\Delta\subset$M$\times$M
in H$_{2n}($M$\times$M,$\mathbb{Z})$ there exists a Green current $g_{\Delta}$
of type $(n-1,n-1)$ with a logarithmic growth along the diagonal $\Delta$ in
M$\times$M such that%

\begin{equation}
\partial\overline{\partial}g_{\Delta}+\delta_{\Delta}=\omega_{\Delta},
\label{LP1}%
\end{equation}
where $\delta_{\Delta}$ is a Dirac current of $\Delta$ and $\omega_{\Delta}$
is a C$^{\infty}$ form on M$\times$M which represents the same cohomology
class as the current $\delta_{\Delta}.$ From the fact that the singular
support of the Dirac kernel is $\Delta$ and the equation $\left(
\ref{LP1}\right)  $ we can conclude that the restriction of the current
$g_{\Delta}$ on M$\times$M $-$ $\Delta$ is represented by a smooth C$^{\infty
}$ closed form of type $(n-1,n-1)$.

\begin{proposition}
\label{LOG0}Let $\mathcal{U}$ be an affine open set in M. Then we have the
following presentation for the restriction of $\omega_{\Delta}$ on
$\mathcal{U}\times\mathcal{U}$
\begin{equation}
\omega_{\Delta}|_{\mathcal{U}\times\mathcal{U}}=\overline{\partial}%
\psi_{\Delta}, \label{LP1a}%
\end{equation}
where $\psi_{\Delta}$ is a smooth form of type $(n,n-1)$ on $\mathcal{U}%
\times\mathcal{U}.$ Moreover the form $\psi_{\Delta}$ has singularities of the
type described by $\left(  \ref{log}\right)  $ along the divisor which is the
complement of $\mathcal{U}\times\mathcal{U}$ in M$\times$M$.$
\end{proposition}

\textbf{Proof: }It is a standard fact that for any coherent sheaf
$\mathfrak{F}$ on the affine set $\mathcal{U}\times\mathcal{U}$ we have
$H^{k}(\mathcal{U}\times\mathcal{U},\mathfrak{F})=0$ for $k>0.$ Since
$\omega_{\Delta}$ is a closed form of type $(n,n)$ representing the Poincare
dual class of the diagonal $\Delta$ in $H^{n}($M$\times$M,$\mathbb{Z})$ its
restriction on $\mathcal{U}\times\mathcal{U}$ is a non zero element of
$H^{n}(\mathcal{U}\times\mathcal{U},\mathbb{Z}).$ On the other hand
$\omega_{\Delta}$ is a closed form of type $(n,n)$ and thus it is
$\overline{\partial}$ closed. Thus its restriction on $\mathcal{U}%
\times\mathcal{U}$ is an element of $H^{n}(\mathcal{U}\times\mathcal{U}%
,\Omega^{n}).$Since $\Omega^{n}$ is a coherent sheaf the restriction of
$\omega_{\Delta}$ on $\mathcal{U}\times\mathcal{U}$ is a $\overline{\partial}$
exact smooth form of type ($n,n$). Thus $\omega_{\Delta}|_{\mathcal{U}%
\times\mathcal{U}}=\overline{\partial}\psi_{\Delta},$ where $\psi_{\Delta}$ is
a smooth form of type $(n,n-1)$ on $\mathcal{U}\times\mathcal{U}$. Since
$\omega_{\Delta}$ is a smooth form on M$\times$M then according to Theorem
\ref{D0} the restriction of $\omega_{\Delta}$ on \ $\mathcal{U}\times
\mathcal{U}$ can be represented as a cohomology class of the de Rham log
complex $\left(  \ref{LOG}\right)  $. Thus we can choose $\psi_{\Delta}$ to
have singularities of the type described by $\left(  \ref{log}\right)  $ along
the divisor which is the complement of $\mathcal{U}\times\mathcal{U}$ in
M$\times$M$.$ Proposition \ref{LOG0} is proved. $\blacksquare$

According to Proposition \ref{LOG0}, the C$^{\infty}$ form $\left(  \frac
{\psi_{\Delta}-\partial g_{\Delta}}{\pi_{1}^{\ast}(\eta)}\right)  $ is of type
$(0,n-1).$ It is non zero and it is well defined on $\mathcal{U}%
\times\mathcal{U}.$ We will show that the integration of any smooth $(2n,n+1)$
form $\beta$ on M$\times$M with this form defines a current on M$\times$M,
which we denote by $\overline{\partial}^{-1}(\overline{\eta}_{\Delta}),$
namely
\begin{equation}
\left\langle \overline{\partial}^{-1}(\overline{\eta}_{\Delta}),\beta
\right\rangle =%
%TCIMACRO{\dint \limits_{\mathcal{U}\times\mathcal{U}-\Delta}}%
%BeginExpansion
{\displaystyle\int\limits_{\mathcal{U}\times\mathcal{U}-\Delta}}
%EndExpansion
\left(  \frac{\psi_{\Delta}-\partial g_{\Delta}}{\pi_{1}^{\ast}(\eta)}\right)
\wedge\beta. \label{ig00}%
\end{equation}

\begin{theorem}
\label{igor9}The current $\overline{\partial}^{-1}(\overline{\eta}_{\Delta})$
is well defined and it does not depend on the choice of $\psi_{\Delta}$ and
$\mathcal{U}.$
\end{theorem}

\textbf{Proof: }Let $\beta$ be any $(2n,n+1)$ smooth form on M$\times$M. First
we will show that the integral $\left(  \ref{ig00}\right)  $ does not depend
on the choice of $\psi_{\Delta}$ when we fix the affine open set $\mathcal{U}%
$. To prove that the current $\overline{\partial}^{-1}(\overline{\eta}%
_{\Delta})$ is well defined it is enough to show that
\begin{equation}%
%TCIMACRO{\dint \limits_{\mathcal{U}\times\mathcal{U}-\Delta}}%
%BeginExpansion
{\displaystyle\int\limits_{\mathcal{U}\times\mathcal{U}-\Delta}}
%EndExpansion
\left(  \frac{\psi_{\Delta}-\partial g_{\Delta}}{\pi_{1}^{\ast}(\eta)}\right)
\wedge\beta\label{Int0}%
\end{equation}
converges. So we need to prove that the integral $\left(  \ref{Int0}\right)  $
converges. In fact the integral%

\begin{equation}%
%TCIMACRO{\dint \limits_{\text{M}\times\text{M}-\Delta}}%
%BeginExpansion
{\displaystyle\int\limits_{\text{M}\times\text{M}-\Delta}}
%EndExpansion
\frac{\partial g_{\Delta}}{\pi_{1}^{\ast}(\eta)}\wedge\beta\label{Int1}%
\end{equation}
converges since the current $g_{\Delta}$ has a logarithmic growth along
$\Delta$. See \cite{AB}. The integral
\[%
%TCIMACRO{\dint \limits_{\text{M}\times\text{M}-\Delta}}%
%BeginExpansion
{\displaystyle\int\limits_{\text{M}\times\text{M}-\Delta}}
%EndExpansion
\frac{\psi_{\Delta}}{\pi_{1}^{\ast}(\eta)}\wedge\beta
\]
converges since according to Proposition \ref{LOG0} $\psi_{\Delta}$ is a
smooth form on $\mathcal{U}\times\mathcal{U}$ and it has singularities along
the complement of $\mathcal{U}\times\mathcal{U}$ in M$\times$M described by
$\left(  \ref{log}\right)  .$ Thus we proved that $\overline{\partial}%
^{-1}(\overline{\eta}_{\Delta})$ defines a current on M$\times$M.

Next we will verify that the current $\overline{\partial}^{-1}(\overline{\eta
}_{\Delta})$ does not depend on the choice of \ $\psi_{\Delta}.$ Indeed if we
choose $\psi_{\Delta}^{^{\prime}}$ in $\mathcal{U}\times\mathcal{U}$ with
logarithmic singularities along the complement of $\mathcal{U}\times
\mathcal{U}$ in M$\times$M such that%
\[
\overline{\partial}\psi_{\Delta}^{^{\prime}}=\overline{\partial}\psi_{\Delta
}=\omega|_{\mathcal{U}\times\mathcal{U}}%
\]
then
\[
\overline{\partial}\left(  \psi_{\Delta}^{^{\prime}}-\psi_{\Delta}\right)
=0.
\]
Thus the class of cohomology of $\left(  \psi_{\Delta}^{^{\prime}}%
-\psi_{\Delta}\right)  $ on $\mathcal{U}\times\mathcal{U}$ is zero. In
particular the singularities will cancel in $\left(  \psi_{\Delta}^{^{\prime}%
}-\psi_{\Delta}\right)  $. Thus $\left(  \psi_{\Delta}^{^{\prime}}%
-\psi_{\Delta}\right)  =\overline{\partial}\phi,$ where $\phi$ is a smooth
form defined on M$\times$M. So
\begin{equation}%
%TCIMACRO{\dint \limits_{\mathcal{U}\times\mathcal{U}}}%
%BeginExpansion
{\displaystyle\int\limits_{\mathcal{U}\times\mathcal{U}}}
%EndExpansion
\left(  \frac{\psi_{\Delta}}{\pi_{1}^{\ast}(\eta)}-\frac{\psi_{\Delta
}^{^{\prime}}}{\pi_{1}^{\ast}(\eta)}\right)  \wedge\beta=%
%TCIMACRO{\dint \limits_{\text{M}\times\text{M}}}%
%BeginExpansion
{\displaystyle\int\limits_{\text{M}\times\text{M}}}
%EndExpansion
d\left(  \frac{\phi}{\pi_{1}^{\ast}(\eta)}\wedge\beta\right)  =0. \label{Int2}%
\end{equation}
Suppose that $\mathcal{U}_{1}$ is another affine open set in M. Since
$\mathcal{U}_{1}\cap\mathcal{U}$ $\ $is an affine set and the complement of
$\mathcal{U}_{1}\cap\mathcal{U}$ in both of the affine open sets has measure
zero then by repeating the arguments that we used to prove $\left(
\ref{Int2}\right)  $ will show that
\[%
%TCIMACRO{\dint \limits_{\mathcal{U}\times\mathcal{U}-\Delta}}%
%BeginExpansion
{\displaystyle\int\limits_{\mathcal{U}\times\mathcal{U}-\Delta}}
%EndExpansion
\left(  \frac{\psi_{\Delta}-\partial g_{\Delta}}{\pi_{1}^{\ast}(\eta)}\right)
\wedge\beta=%
%TCIMACRO{\dint \limits_{\mathcal{U}_{1}\times\mathcal{U}_{1}-\Delta}}%
%BeginExpansion
{\displaystyle\int\limits_{\mathcal{U}_{1}\times\mathcal{U}_{1}-\Delta}}
%EndExpansion
\left(  \frac{\psi_{1,\Delta}-\partial g_{\Delta}}{\pi_{1}^{\ast}(\eta
)}\right)  \wedge\beta
\]
Theorem \ref{igor9} is proved. $\blacksquare$

\begin{corollary}
The form $\frac{\psi_{\Delta}-\partial g_{\Delta}}{\pi_{1}^{\ast}(\eta)}$ is a
$\overline{\partial}$ closed form of type $(0,n-1)$ on M$\times$M$-\Delta$ and
it defines a class of cohomology in $H_{\overline{\partial}}^{n-1}($M$\times
$M$-\Delta,\mathcal{O}_{\text{M}\times\text{M}-\Delta}).$
\end{corollary}

\begin{remark}
We would like to note that since the cohomology class of $\overline{\eta
}_{\Delta}$ is non zero one can not define $\overline{\partial}^{-1}%
(\overline{\eta}_{\Delta})$ in a straightforward way. We hope that our
notation would not lead to a possible confusion.
\end{remark}

We will call the current $\overline{\partial}^{-1}(\overline{\eta}_{\Delta})$
the Green kernel of the operator $\overline{\partial},$ and we will often use
the same notation for its defining $C^{\infty}$ form in $\left(
\ref{ig00}\right)  $.

\subsection{Symmetric Form of Holomorphic Linking and Complex Linking Number}

From now on we will suppose that $\Sigma_{1}$ and $\Sigma_{2}$ are Riemann
surfaces embedded in a CY threefold and $\Sigma_{1}\cap\Sigma_{2}=\emptyset.$

\begin{theorem}
\label{igor14}Let us consider the embedding $\Sigma_{1}\times\Sigma_{2}%
\subset$M$\times$M$.$ Then the following formula holds:
\begin{equation}
\#\left(  \left(  \Sigma_{1},\theta_{1}\right)  ,\left(  \Sigma_{2},\theta
_{2}\right)  \right)  =\int_{\Sigma_{1}\times\Sigma_{2}}\left(  \overline
{\partial}^{-1}\left(  \overline{\eta}_{\Delta}\right)  \right)  |_{\Sigma
_{1}\times\Sigma_{2}}\wedge\pi_{1}^{\ast}(\theta_{1})\wedge\pi_{2}^{\ast
}(\theta_{2}). \label{Ig}%
\end{equation}

\end{theorem}

\textbf{Proof:} We can rewrite the holomorphic linking in the following form:%
\[
\#\left(  \left(  \Sigma_{1},\theta_{1}\right)  ,\left(  \Sigma_{2},\theta
_{2}\right)  \right)  =\int_{\text{M}}\overline{\partial}^{-1}\left(  \left(
\overline{\theta}_{1}\right)  _{\Sigma_{1}}\right)  \wedge\left(
\overline{\theta}_{2}\right)  _{\Sigma_{2}}\wedge\eta=
\]%
\begin{equation}
\int_{\text{M$\times$M}}\delta_{\Delta}\wedge\pi_{1}^{\ast}\left(
\overline{\partial}^{-1}\left(  \left(  \overline{\theta}_{1}\right)
_{\Sigma_{1}}\right)  \right)  \wedge\pi_{2}^{\ast}\left(  \left(
\overline{\theta}_{2}\right)  _{\Sigma_{2}}\wedge\eta\right)  . \label{Li2}%
\end{equation}
Substituting in $\left(  \ref{Li2}\right)  $ the expression for the current
$\delta_{\Delta}$ stated in $\left(  \ref{LP1}\right)  :$
\[
\delta_{\Delta}=\pi_{1}^{\ast}(\eta)\wedge\left(  \frac{\omega_{\Delta
}-\overline{\partial}\partial g_{\Delta}}{\pi_{1}^{\ast}(\eta)}\right)  ,
\]
where
\[
\omega_{\Delta}|_{\text{M$\times$M}-\Delta}=\overline{\partial}\psi_{\Delta}%
\]
and the $\omega_{\Delta}$ is a form representing the Poincare dual class of
$\Delta,$ we obtain $\,$that%
\[
\#\left(  \left(  \Sigma_{1},\theta_{1}\right)  ,\left(  \Sigma_{2},\theta
_{2}\right)  \right)  =
\]%
\[
\int_{\text{M$\times$M}}\left(  \left(  \frac{\omega_{\Delta}-\overline
{\partial}\partial g_{\Delta}}{\pi_{1}^{\ast}(\eta)}\right)  \right)
\wedge\left(  \overline{\partial}^{-1}\left(  \pi_{1}^{\ast}\left(  \left(
\overline{\theta}_{1}\right)  _{\Sigma_{1}}\wedge\eta\right)  \right)
\right)  \wedge\pi_{2}^{\ast}\left(  \left(  \overline{\theta}_{2}\right)
_{\Sigma_{2}}\wedge\eta\right)  =
\]%
\[
\int_{\text{M$\times$M}-\Delta}\left(  \overline{\partial}\left(  \frac
{\psi_{\Delta}}{\pi_{1}^{\ast}(\eta)}\right)  \right)  \wedge\left(
\overline{\partial}^{-1}\left(  \pi_{1}^{\ast}\left(  \left(  \overline
{\theta}_{1}\right)  _{\Sigma_{1}}\wedge\eta\right)  \right)  \right)
\wedge\pi_{2}^{\ast}\left(  \left(  \overline{\theta}_{2}\right)  _{\Sigma
_{2}}\right)  -
\]%
\begin{equation}
\int_{\text{M$\times$M}}\left(  \overline{\partial}\left(  \frac{\partial
g_{\Delta}}{\pi_{1}^{\ast}(\eta)}\right)  \right)  \wedge\left(
\overline{\partial}^{-1}\left(  \pi_{1}^{\ast}\left(  \left(  \overline
{\theta}_{1}\right)  _{\Sigma_{1}}\wedge\eta\right)  \right)  \right)
\wedge\pi_{2}^{\ast}\left(  \left(  \overline{\theta}_{2}\right)  _{\Sigma
_{2}}\right)  . \label{Li3}%
\end{equation}
From Stokes' Theorem we obtain
\[
\#\left(  \left(  \Sigma_{1},\theta_{1}\right)  ,\left(  \Sigma_{2},\theta
_{2}\right)  \right)  =
\]%
\[
\int_{\text{M$\times$M$-$}\Delta}\frac{\psi_{\Delta}}{\pi_{1}^{\ast}(\eta
)}\wedge\left(  \pi_{1}^{\ast}\left(  \left(  \overline{\theta}_{1}\right)
_{\Sigma_{1}}\wedge\eta\right)  \right)  \wedge\pi_{2}^{\ast}\left(  \left(
\overline{\theta}_{2}\right)  _{\Sigma_{2}}\wedge\eta\right)  -
\]%
\[
\int_{\text{M$\times$M$-$}\Delta}\frac{\partial g_{\Delta}}{\pi_{1}^{\ast
}(\eta)}\wedge\pi_{1}^{\ast}\left(  \left(  \overline{\theta}_{1}\right)
_{\Sigma_{1}}\wedge\eta\right)  \wedge\pi_{2}^{\ast}\left(  \left(
\overline{\theta}_{2}\right)  _{\Sigma_{2}}\wedge\eta\right)  =
\]%
\begin{equation}
\int_{\text{M$\times$M$-$}\Delta}\left(  \frac{\psi_{\Delta}-\partial
g_{\Delta}}{\pi_{1}^{\ast}(\eta)}\right)  \wedge\pi_{1}^{\ast}\left(  \left(
\overline{\theta}_{1}\right)  _{\Sigma_{1}}\wedge\eta\right)  \wedge\pi
_{2}^{\ast}\left(  \left(  \overline{\theta}_{2}\right)  _{\Sigma_{2}}%
\wedge\eta\right)  . \label{Li5}%
\end{equation}
The antiholomorphic Dirac currents $\left(  \overline{\theta}_{1}\right)
_{\Sigma_{1}},$ $\left(  \overline{\theta}_{2}\right)  _{\Sigma_{2}}$ and
$\overline{\eta}_{\Delta}$ have disjoint singularities sets. This fact
guarantees that the current $\left(  \pi_{1}^{\ast}\left(  \left(
\overline{\theta}_{1}\right)  _{\Sigma_{1}}\right)  \right)  \wedge\pi
_{2}^{\ast}\left(  \left(  \overline{\theta}_{2}\right)  _{\Sigma_{2}}\right)
$ is well defined. Thus the convergence of all integrals is established.
Finally using the formula $\left(  \ref{ig00}\right)  $ for $\overline
{\partial}^{-1}\left(  \overline{\eta}_{\Delta}\right)  $ and the definition
of the antiholomorphic Dirac currents $\pi_{1}^{\ast}\left(  \left(
\overline{\theta}_{1}\right)  _{\Sigma_{1}}\right)  $ and $\pi_{2}^{\ast
}\left(  \left(  \overline{\theta}_{2}\right)  _{\Sigma_{2}}\right)  $ we
deduce from $\left(  \ref{Li5}\right)  $ that
\[
\#\left(  \left(  \Sigma_{1},\theta_{1}\right)  ,\left(  \Sigma_{2},\theta
_{2}\right)  \right)  =
\]%
\[
\int_{\text{M$\times$M $-$ }\Delta}\left(  \overline{\partial}^{-1}\left(
\overline{\eta}_{\Delta}\right)  \right)  \wedge\left(  \pi_{1}^{\ast}\left(
\left(  \overline{\theta}_{1}\right)  _{\Sigma_{1}}\wedge\eta\right)  \right)
\wedge\pi_{2}^{\ast}\left(  \left(  \overline{\theta}_{2}\right)  _{\Sigma
_{2}}\wedge\eta\right)  =
\]%
\begin{equation}
\int_{\Sigma_{1}\times\Sigma_{2}}\left(  \overline{\partial}^{-1}\left(
\overline{\eta}_{\Delta}\right)  \right)  |_{\Sigma_{1}\times\Sigma_{2}}%
\wedge\pi_{1}^{\ast}(\theta_{1})\wedge\pi_{2}^{\ast}(\theta_{2}). \label{Ig2}%
\end{equation}
Thus we established formula $\left(  \ref{Ig}\right)  $. Theorem \ref{igor14}
is proved. $\blacksquare$

\begin{corollary}
The holomorphic linking does not depend on the choice of the representative of
the current $\overline{\partial}^{-1}(\overline{\eta}_{\Delta}).$
\end{corollary}

\textbf{Proof: }The corollary follows directly from Stokes' Theorem.
$\blacksquare$

The expression of the holomorphic linking given by Definition \ref{igor8} and
Theorem \ref{igor14} can be viewed as complex counterparts of the
corresponding expressions for the Gauss linking number reviewed in Section 2.
The essential difference between the two notions is that in the real case we
obtain a topological invariant while in the complex case our linking map
depends on the way the Riemannian surfaces are embedded in the CY manifold M
and the choice of holomorphic forms on them but does not depend on the choice
of the metric on CY manifold M.

The comparison of our holomorphic linking and the Gauss linking number
suggests an alternative definition:

\begin{definition}
\label{igor16} The complex linking number of two Riemann surfaces $\Sigma_{1}$
and $\Sigma_{2}$ embedded in a CY threefold such that $\Sigma_{1}\cap
\Sigma_{2}=\emptyset$ is defined by the formula:
\[
\#\left(  \Sigma_{1},\Sigma_{2}\right)  =\int_{\Sigma_{1}\times\Sigma_{2}%
}\left(  \overline{\partial}^{-1}\left(  \overline{\eta}_{\Delta}\right)
\right)  \wedge\left(  \partial^{-1}\left(  \eta_{\Delta}\right)  \right)  ,
\]
where $\eta$ is normalized by condition $\left(  \ref{14}\right)  $.
\end{definition}

Simple examples show that the complex linking number contains less information
than the holomorphic linking. It is an interesting question to find a relation
between the two invariants.

\section{Analytic Expression of Holomorphic Linking}

In this section we will give an explicit integral formula for the holomorphic
linking. It can be viewed as a complex counterpart of the Gauss formula for
the linking number. The key ingredient of our formula is an expression for the
Green kernel $\overline{\partial}^{-1}(\overline{\eta}_{\Delta})$ in terms of
the Bochner-Martinelli form of type (n,n-1). The letter form written in affine
coordinates is precisely the classical Bochner-Martinelli form that yields a
generalization of the Cauchy integral formula. It turns out that it finds
another application in the integral formula for the holomorphic linking.

\subsection{Bochner-Martinelli Kernel\label{BMK}}

Let $\mathcal{U}$ be an affine open set in Zariski topology of the CY manifold
M. Let $z^{1},...,z^{n}$ be local coordinates\ in $\mathcal{U}$ such that
restriction of the holomorphic n-form $\eta$ on $\mathcal{U}$ is expressed
as:
\[
\eta\left\vert _{\mathcal{U}}\right.  =dz^{1}\wedge...\wedge dz^{n}.
\]
Let
\[
\mathit{\ }\Phi_{j}(z):=(-1)^{j-1}z^{j}dz^{1}\wedge...\wedge dz^{j-1}\wedge
dz^{j+1}\wedge...\wedge dz^{n}.
\]
Following \cite{GH} Chapter 3 we will define the Bochner-Martinelli kernel on
\[
\mathcal{U}\times\mathcal{U}\text{ }-\text{ }\Delta_{\mathcal{U}},\text{
}\Delta_{\mathcal{U}}:=\mathcal{U\times U}\cap\Delta,
\]
by the formula: \textit{ }%
\begin{equation}
\mathcal{K}_{\mathcal{U}\times\mathcal{U}}^{n,n-1}=\frac{C_{n}\overline
{\left(  \sum_{r=1}^{n}\Phi_{r}(w-z\right)  }\wedge\left(  dz^{1}%
\wedge...\wedge dz^{n}\right)  }{\left\Vert w-z\right\Vert ^{2n}},
\label{Boch}%
\end{equation}
$C_{n}$ is the volume of the unit $2n-1$ dimensional sphere, $\{w^{k}\}$ and
$\{z^{k}\}$ are local coordinates in $\mathcal{U}\times\mathcal{U}$%
\textit{\ }such that\textit{ }
\begin{equation}
\pi_{1}^{\ast}(\eta)|_{\mathcal{U}\times\mathcal{U}}=dw^{1}\wedge...\wedge
dw^{n},\text{ }\pi_{2}^{\ast}(\eta)|_{\mathcal{U}\times\mathcal{U}}%
=dz^{1}\wedge...\wedge dz^{n}, \label{COOR}%
\end{equation}
and%
\[
\left\Vert w-z\right\Vert ^{2n}=\left(  \sum_{k=1}^{n}\left\vert w^{k}%
-z^{k}\right\vert ^{2}\right)  ^{n}.
\]

It is proved in Chapter 3 of \cite{GH} that on $\mathcal{U}\times\mathcal{U}$
$-$ $\Delta_{\mathcal{U}}$ the forms $\mathcal{K}_{\mathcal{U}\times
\mathcal{U}}^{n,n-1}$ of type $(n,n-1)$ given by the expression $\left(
\ref{Boch}\right)  $ are $d$ and so $\overline{\partial}$ closed, i.e.
\begin{equation}
\overline{\partial}\left(  \mathcal{K}_{\mathcal{U}\times\mathcal{U}}%
^{n,n-1}\right)  |_{\mathcal{U}\times\mathcal{U}-\Delta_{\mathcal{U}}%
}=d\left(  \mathcal{K}_{\mathcal{U}\times\mathcal{U}}^{n,n-1}\right)
|_{\mathcal{U}\times\mathcal{U}-\Delta_{\mathcal{U}}}=0. \label{Boch2}%
\end{equation}

Let $T_{\varepsilon}(\Delta)$ be a tubular neighborhood of the diagonal
$\Delta$ in M$\times$M. In \cite{GH} Chapter 5 it is proved that the limit%
\[
\underset{\varepsilon\rightarrow0}{\lim}%
%TCIMACRO{\dint \limits_{\partial(T_{\varepsilon}(\Delta)\cap\mathcal{U}%
%\times\mathcal{U}}}%
%BeginExpansion
{\displaystyle\int\limits_{\partial(T_{\varepsilon}(\Delta)\cap\mathcal{U}%
\times\mathcal{U}}}
%EndExpansion
\mathcal{K}_{\mathcal{U}\times\mathcal{U}}^{n,n-1}\wedge\omega
\]
exists. The form $\mathcal{K}_{\mathcal{U}\times\mathcal{U}}^{n,n-1}$ defines
a current on M$\times$M as follows%
\begin{equation}
\left\langle \mathcal{K}_{\mathcal{U}\times\mathcal{U}}^{n,n-1},\beta
\right\rangle :=\underset{\varepsilon\rightarrow0}{\lim}%
%TCIMACRO{\dint \limits_{\mathcal{U}\times\mathcal{U}\text{ $-$ }\left(
%T_{\varepsilon}(\Delta)\cap\mathcal{U}\times\mathcal{U}\right)  }}%
%BeginExpansion
{\displaystyle\int\limits_{\mathcal{U}\times\mathcal{U}\text{ $-$ }\left(
T_{\varepsilon}(\Delta)\cap\mathcal{U}\times\mathcal{U}\right)  }}
%EndExpansion
\mathcal{K}_{\mathcal{U}\times\mathcal{U}}^{n,n-1}\wedge\beta, \label{Boch2a}%
\end{equation}
where $\beta$ is any smooth form of type $(n,n+1)$ on M$\times$M. Based on
Stokes' theorem and $\left(  \ref{Boch2a}\right)  $ we can compute\ the
current
\[
\overline{\partial}\left(  \mathcal{K}_{\mathcal{U}\times\mathcal{U}}%
^{n,n-1}\right)  =d\left(  \mathcal{K}_{\mathcal{U}\times\mathcal{U}}%
^{n,n-1}\right)
\]
on M$\times$M as follows
\begin{equation}
\left\langle \overline{\partial}\mathcal{K}_{\mathcal{U}\times\mathcal{U}%
}^{n,n-1},\omega\right\rangle =\left\langle d\mathcal{K}_{\mathcal{U}%
\times\mathcal{U}}^{n,n-1},\omega\right\rangle =\underset{\varepsilon
\rightarrow0}{\lim}%
%TCIMACRO{\dint \limits_{\partial(T_{\varepsilon}(\Delta)\cap\mathcal{U}%
%\times\mathcal{U}}}%
%BeginExpansion
{\displaystyle\int\limits_{\partial(T_{\varepsilon}(\Delta)\cap\mathcal{U}%
\times\mathcal{U}}}
%EndExpansion
\mathcal{K}_{\mathcal{U}\times\mathcal{U}}^{n,n-1}\wedge\omega\label{29}%
\end{equation}
for any smooth form $\omega$ of type $(n,n)$ on M$\times$M.

\begin{theorem}
\label{igor2001}Let M be a CY manifold. \textit{Then we have the following
equality of currents: }%
\begin{equation}
\overline{\partial}(\mathcal{K}_{\mathcal{U}\times\mathcal{U}}^{n,n-1}%
)=\delta_{\Delta}|_{\mathcal{U}\times\mathcal{U}} \label{Boch1}%
\end{equation}
for any affine open set $\mathcal{U}$ $\ $in M.
\end{theorem}

\textbf{Proof:} In order to prove Theorem \ref{igor2001} we need to prove that
for any smooth form $\omega$ of type $(n,n)$ on M we have%
\begin{equation}
\underset{\varepsilon\rightarrow0}{\lim}%
%TCIMACRO{\dint \limits_{\partial(T_{\varepsilon}(\Delta)\cap\mathcal{U}%
%\times\mathcal{U}}}%
%BeginExpansion
{\displaystyle\int\limits_{\partial(T_{\varepsilon}(\Delta)\cap\mathcal{U}%
\times\mathcal{U}}}
%EndExpansion
\mathcal{K}_{\mathcal{U}\times\mathcal{U}}^{n,n-1}\wedge\omega=%
%TCIMACRO{\dint \limits_{\mathcal{U}_{\Delta}}}%
%BeginExpansion
{\displaystyle\int\limits_{\mathcal{U}_{\Delta}}}
%EndExpansion
\omega=%
%TCIMACRO{\dint \limits_{\Delta}}%
%BeginExpansion
{\displaystyle\int\limits_{\Delta}}
%EndExpansion
\omega. \label{Boch3}%
\end{equation}
Since $d\left(  \mathcal{K}_{\mathcal{U}\times\mathcal{U}}^{n,n-1}\right)  =0$
on $\mathcal{U}\times\mathcal{U}-\left(  \mathcal{U}_{\Delta}\right)  $ we
have the following equation of currents:
\[
d\left(  \mathcal{K}_{\mathcal{U}\times\mathcal{U}}^{n,n-1}\wedge
\omega\right)  =d\left(  \mathcal{K}_{\mathcal{U}\times\mathcal{U}}%
^{n,n-1}\right)  \wedge\omega+(-1)^{2n-1}\mathcal{K}_{\mathcal{U}%
\times\mathcal{U}}^{n,n-1}\wedge\overline{\partial}\omega=
\]%
\begin{equation}
-\mathcal{K}_{\mathcal{U}\times\mathcal{U}}^{n,n-1}\wedge\overline{\partial
}\omega. \label{Boch4}%
\end{equation}
Let us denote $\left(  \mathcal{U}\times\mathcal{U}\right)  \mathcal{\cap
}T_{\varepsilon}(\Delta)=T_{\epsilon}(\Delta_{\mathcal{U}}).$ Formula $\left(
\ref{Boch4}\right)  $ and the Stokes' Theorem imply%
\[%
%TCIMACRO{\dint \limits_{\mathcal{U\times U}\text{ }-\text{ }T_{\epsilon
%}(\Delta_{\mathcal{U}})}}%
%BeginExpansion
{\displaystyle\int\limits_{\mathcal{U\times U}\text{ }-\text{ }T_{\epsilon
}(\Delta_{\mathcal{U}})}}
%EndExpansion
d\left(  \mathcal{K}_{\mathcal{U}\times\mathcal{U}}^{n,n-1}\wedge
\omega\right)  =
\]%
\begin{equation}
-%
%TCIMACRO{\dint \limits_{\mathcal{U\times U}\text{ }-\text{ }T_{\epsilon
%}(\Delta_{\mathcal{U}})}}%
%BeginExpansion
{\displaystyle\int\limits_{\mathcal{U\times U}\text{ }-\text{ }T_{\epsilon
}(\Delta_{\mathcal{U}})}}
%EndExpansion
\mathcal{K}_{\mathcal{U}\times\mathcal{U}}^{n,n-1}\wedge\overline{\partial
}\omega=-%
%TCIMACRO{\dint \limits_{\partial\left(  T_{\epsilon}(\Delta_{\mathcal{U}%
%})\right)  }}%
%BeginExpansion
{\displaystyle\int\limits_{\partial\left(  T_{\epsilon}(\Delta_{\mathcal{U}%
})\right)  }}
%EndExpansion
\mathcal{K}_{\mathcal{U}\times\mathcal{U}}^{n,n-1}\wedge\omega. \label{Boch5}%
\end{equation}
We notice that $\partial(T_{\varepsilon}(\Delta)\cap\mathcal{U}\times
\mathcal{U}$ is a fibration of $2n-1$ spheres over $\mathcal{U}_{\Delta}.$
Since $\omega$ is a form of type $(n,n),$ then when we restrict it to the
diagonal we will get%
\[
\omega|_{\Delta}=\phi(w)\left(  \eta\wedge\overline{\eta}\right)
\]
for some C$^{\infty}$ function $\phi(w).$ Now we may apply the Fubini theorem
and rewrite the integral in $\left(  \ref{Boch5}\right)  $ as follows:
\begin{equation}%
%TCIMACRO{\dint \limits_{\partial(T_{\varepsilon}(\Delta)\cap\mathcal{U}%
%\times\mathcal{U}}}%
%BeginExpansion
{\displaystyle\int\limits_{\partial(T_{\varepsilon}(\Delta)\cap\mathcal{U}%
\times\mathcal{U}}}
%EndExpansion
\mathcal{K}_{\mathcal{U}\times\mathcal{U}}^{n,n-1}\wedge\omega=%
%TCIMACRO{\dint \limits_{\mathcal{U}_{\Delta}}}%
%BeginExpansion
{\displaystyle\int\limits_{\mathcal{U}_{\Delta}}}
%EndExpansion
\text{ }\left(
%TCIMACRO{\dint \limits_{\left\Vert z-w\right\Vert =\varepsilon,\omega\in
%\Delta}}%
%BeginExpansion
{\displaystyle\int\limits_{\left\Vert z-w\right\Vert =\varepsilon,\omega
\in\Delta}}
%EndExpansion
\phi(z)\mathcal{K}_{\mathcal{U}\times\mathcal{U}}^{n,n-1}\right)  \eta
\wedge\overline{\eta}. \label{Boch6}%
\end{equation}
The expression in $\left(  \ref{Boch6}\right)  $ makes sense since in
\cite{GH} Chapter 3 paragraph 1 it is proved that if in the expression of the
Bochner-Martinelli kernel given by $\left(  \ref{Boch}\right)  $ we
substituted $z=w_{0}$ then
\begin{equation}%
%TCIMACRO{\dint \limits_{||z-w_{0}||=\varepsilon}}%
%BeginExpansion
{\displaystyle\int\limits_{||z-w_{0}||=\varepsilon}}
%EndExpansion
\phi(z)\mathcal{K}_{\mathcal{U}\times\mathcal{U}}^{n,n-1}=\phi(w_{0})
\label{Boch7}%
\end{equation}
for any continuous function $\phi.$ Using the representation of the
restriction of $\omega$ on $\Delta$ together with $\left(  \ref{Boch6}\right)
$ and $\left(  \ref{Boch7}\right)  $ we get formula $\left(  \ref{Boch3}%
\right)  $. Theorem \ref{igor2001} is proved. $\blacksquare$

We define \textit{ }%
\begin{equation}
\mathcal{K}_{\mathcal{U}\times\mathcal{U}}^{0,n-1}=\frac{C_{n}\overline
{\left(  \sum_{r=1}^{n}\Phi_{r}(w-z\right)  }}{\left\Vert w-z\right\Vert
^{2n}} \label{Bocha}%
\end{equation}
and we call $\mathcal{K}_{\mathcal{U}\times\mathcal{U}}^{0,n-1}$ the
holomorphic Bochner-Martinelli kernel. We have the following formula
\begin{equation}
\mathcal{K}_{\mathcal{U}\times\mathcal{U}}^{0,n-1}\wedge\pi_{1}^{\ast}%
(\eta)=\mathcal{K}_{\mathcal{U}\times\mathcal{U}}^{n,n-1} \label{Bochb}%
\end{equation}
The form $\mathcal{K}_{\mathcal{U}\times\mathcal{U}}^{0,n-1}$ defines a
current of type $(0,n-1)$ on M$\times$M as follows%
\begin{equation}
\left\langle \mathcal{K}_{\mathcal{U}\times\mathcal{U}}^{0,n-1},\beta
\right\rangle :=\underset{\varepsilon\rightarrow0}{\lim}%
%TCIMACRO{\dint \limits_{\mathcal{U}\times\mathcal{U}\text{ $-$ }\left(
%T_{\varepsilon}(\Delta)\cap\mathcal{U}\times\mathcal{U}\right)  }}%
%BeginExpansion
{\displaystyle\int\limits_{\mathcal{U}\times\mathcal{U}\text{ $-$ }\left(
T_{\varepsilon}(\Delta)\cap\mathcal{U}\times\mathcal{U}\right)  }}
%EndExpansion
\mathcal{K}_{\mathcal{U}\times\mathcal{U}}^{0,n-1}\wedge\beta, \label{Bochc}%
\end{equation}
where $\beta$ is any smooth form of type $(2n,n+1)$ on M$\times$M.

By using formulas $\left(  \ref{Bochb}\right)  ,$ $\left(  \ref{Bochc}\right)
$ and by repeating the arguments of the proof of Theorem \ref{igor2001} we obtain

\begin{theorem}
\label{igor2001a}Let M be a CY manifold. \textit{Then we have the following
equality of currents: }%
\begin{equation}
\overline{\partial}(\mathcal{K}_{\mathcal{U}\times\mathcal{U}}^{0,n-1}%
)=\overline{\eta}_{\Delta}|_{\mathcal{U}\times\mathcal{U}}, \label{Boch1a}%
\end{equation}
where $\overline{\eta}_{\Delta}$ is \textit{the antiholomorphic Dirac kernel
}defined in Definition\textbf{ }\ref{andrey48}.
\end{theorem}

\begin{remark}
Instead of the equation $\left(  \ref{LP}\right)  $ we could consider the
holomorphic analogue of (3):
\begin{equation}
\overline{\partial}\text{S}+\overline{\eta}_{\Delta}=\omega_{\Delta},
\label{hbl}%
\end{equation}
where $\omega_{\Delta}$ is a smooth form on M$\times$M which realizes the
class of the current $\overline{\eta}_{\Delta}.$ Then the holomorphic
Bochner-Martinelli kernel represents a solution S of $\left(  \ref{hbl}%
\right)  $ on M$\times$M$-\Delta.$
\end{remark}

\subsection{The Holomorphic Analogue of the Gauss Formula}

Now we are ready to reexpress the formula for the holomorphic linking
established in Theorem \ref{igor14} for the complex linking number introduced
in Definition \ref{igor16} using the Bochner-Martinelli kernel. This gives us
the complex analogue of the Gauss integral formula for linking number.

Let $\Sigma_{1}$ and $\Sigma_{2}$ be two Riemann surfaces embedded in a CY
threefold M such that
\[
\Sigma_{1}\cap\Sigma_{2}=\emptyset.
\]
Let $\mathcal{U}$ be an affine open set in $\ $M$,$ then $\mathcal{U}%
\cap\Sigma_{i}$ are affine open sets in $\Sigma_{i},$ i.e. $\mathcal{U}%
\cap\Sigma_{i}$ are the Riemann surfaces $\Sigma_{i}$ minus finite number of
points. Suppose that $\theta_{1}$ and $\theta_{2}$ are two non-zero
holomorphic forms on $\Sigma_{1}$ and $\Sigma_{2}.$ We also may assume that
\begin{equation}
\theta_{1}\left\vert _{\mathcal{U\cap}\Sigma_{1}}\right.  =f_{1}%
(w)dw\ \text{\ and }\theta_{2}\left\vert _{\mathcal{U\cap}\Sigma_{2}}\right.
=f_{2}(z)dz. \label{theta}%
\end{equation}
Next we will derive an explicit integral form of the holomorphic linking:

\begin{theorem}
\label{b1}The following formula holds:
\[
\#((\Sigma_{1},\theta_{1}),(\Sigma_{2},\theta_{2}))=C_{3}\int_{(\mathcal{U}%
\cap\Sigma_{1})\times(\mathcal{U}\cap\Sigma_{2})}\mathcal{K}_{\mathcal{U}%
\times\mathcal{U}}^{0,2}\left(  f_{1}(w)dw\right)  \wedge\left(
f_{2}(z)dz\right)
\]%
\begin{equation}
C_{3}\int_{(\mathcal{U}\cap\Sigma_{1})\times(\mathcal{U}\cap\Sigma_{2})}%
\frac{\left(  \overline{\sum_{i=1}^{3}\Phi_{i}(z-w)}\right)  \wedge\left(
f_{1}(w)dw\right)  \wedge\left(  f_{2}(z)dz\right)  }{\left(  \sum_{j=1}%
^{3}\left\vert z^{j}-w^{j}\right\vert ^{2}\right)  ^{3}}, \label{Igo}%
\end{equation}
where $\mathcal{U}$ is any open affine subset in M, C$_{3}$ is the volume of
the unit sphere in $\mathbb{C}^{3}$ and the coordinates in $\mathcal{U}%
\times\mathcal{U}$ are chosen as in $\left(  \ref{COOR}\right)  $.
\end{theorem}

\textbf{Proof:} According to formula $\left(  \ref{Li5}\right)  $ we have:%
\[
\#\left(  \left(  \Sigma_{1},\theta_{1}\right)  ,\left(  \Sigma_{2},\theta
_{2}\right)  \right)  =
\]%
\begin{equation}
\int_{\text{M$\times$M$-$}\Delta}\left(  \frac{\psi_{\Delta}-\partial
g_{\Delta}}{\pi_{1}^{\ast}(\eta)}\right)  \wedge\pi_{1}^{\ast}\left(  \left(
\overline{\theta}_{1}\right)  _{\Sigma_{1}}\wedge\eta\right)  \wedge\pi
_{2}^{\ast}\left(  \left(  \overline{\theta}_{2}\right)  _{\Sigma_{2}}%
\wedge\eta\right)  , \label{Ig3}%
\end{equation}
where $\psi_{\Delta}$and $g_{\Delta}$ are defined as in Section 4.1. Since the
complement to the affine set $\mathcal{U}$ in M has measure zero, we can
rewrite the formula $\left(  \ref{Ig3}\right)  $ as follows:%
\[
\#\left(  \left(  \Sigma_{1},\theta_{1}\right)  ,\left(  \Sigma_{2},\theta
_{2}\right)  \right)  =
\]%
\begin{equation}
\int_{\mathcal{U}\text{$\times$}\mathcal{U}\text{ $-$ }\Delta_{\mathcal{U}}%
}\left(  \frac{\psi_{\Delta}-\partial g_{\Delta}}{\pi_{1}^{\ast}(\eta
)}\right)  \wedge\pi_{1}^{\ast}\left(  \left(  \overline{\theta}_{1}\right)
_{\Sigma_{1}}\wedge\eta\right)  \wedge\pi_{2}^{\ast}\left(  \left(
\overline{\theta}_{2}\right)  _{\Sigma_{2}}\wedge\eta\right)  . \label{Ig4}%
\end{equation}
Comparing the equality $\left(  \ref{Ig4}\right)  $ with the equality $\left(
\ref{Boch1}\right)  $ of Theorem \ref{igor2001} we can substitute the current
$\psi_{\Delta}-\partial g_{\Delta}$ with the Bochner-Martinelli kernel
$\mathcal{K}_{\mathcal{U\times U}}^{3,2}$ since their derivatives restricted
on $\mathcal{U\times U}$ give the Dirac current $\delta_{\Delta}$ of the
diagonal restricted on $\mathcal{U\times U}$ and we get%
\[
\#\left(  \left(  \Sigma_{1},\theta_{1}\right)  ,\left(  \Sigma_{2},\theta
_{2}\right)  \right)  =
\]%
\[
\int_{\mathcal{U}\text{$\times$}\mathcal{U}\text{ $-$ }\Delta_{\mathcal{U}}%
}\left(  \frac{\mathcal{K}_{\mathcal{U}\times\mathcal{U}}^{3,2}}{\pi_{1}%
^{\ast}(\eta)}\right)  \wedge\pi_{1}^{\ast}\left(  \left(  \overline{\theta
}_{1}\right)  _{\Sigma_{1}}\wedge\eta\right)  \wedge\pi_{2}^{\ast}\left(
\left(  \overline{\theta}_{2}\right)  _{\Sigma_{2}}\wedge\eta\right)  =
\]%
\[
\int_{\mathcal{U}\text{$\times$}\mathcal{U}\text{ $-$ }\Delta_{\mathcal{U}}%
}\mathcal{K}_{\mathcal{U}\times\mathcal{U}}^{0,2}\wedge\pi_{1}^{\ast}\left(
\left(  \overline{\theta}_{1}\right)  _{\Sigma_{1}}\wedge\eta\right)
\wedge\pi_{2}^{\ast}\left(  \left(  \overline{\theta}_{2}\right)  _{\Sigma
_{2}}\wedge\eta\right)  =
\]%
\begin{equation}
\int_{\left(  \mathcal{U}\cap\Sigma_{1}\right)  \text{$\times$}\left(
\mathcal{U}\cap\Sigma_{2}\right)  \text{ }}\mathcal{K}_{\mathcal{U}%
\times\mathcal{U}}^{0,2}\wedge\theta_{1}\wedge\theta_{2}. \label{Ig6}%
\end{equation}
Substitituing the expression for the holomorphic Bochner-Martinelli kernel
$\left(  \ref{Bocha}\right)  $ and the local expressions for $\theta_{1}$ and
$\theta_{2}$ in $\left(  \ref{theta}\right)  $ we obtain
\[
\#\left(  \left(  \Sigma_{1},\theta_{1}\right)  ,\left(  \Sigma_{2},\theta
_{2}\right)  \right)  =
\]%
\[
C_{3}\int_{(\mathcal{U}\cap\Sigma_{1})\times(\mathcal{U}\cap\Sigma_{2})}%
\frac{\left(  \overline{\sum_{i=1}^{3}\Phi_{i}(z-w)}\right)  \wedge\left(
f_{1}(w)dw\right)  \wedge\left(  f_{2}(z)dz\right)  }{\left(  \sum_{j=1}%
^{3}\left\vert z^{j}-w^{j}\right\vert ^{2}\right)  ^{3}}.
\]
Theorem \ref{b1} is proved. $\blacksquare$

In the same way we prove.

\begin{corollary}
\label{b2}The following formula holds:
\[
\#(\Sigma_{1},\Sigma_{2})=C_{3}\int_{(\mathcal{U}\cap\Sigma_{1})\times
(\mathcal{U}\cap\Sigma_{2})}\frac{\left(
%TCIMACRO{\dsum \limits_{j=1}^{3}}%
%BeginExpansion
{\displaystyle\sum\limits_{j=1}^{3}}
%EndExpansion
\overline{\Phi_{j}(z-w)}\right)  \wedge\left(
%TCIMACRO{\dsum \limits_{j=1}^{3}}%
%BeginExpansion
{\displaystyle\sum\limits_{j=1}^{3}}
%EndExpansion
\Phi_{j}(z-w)\right)  }{\left(
%TCIMACRO{\dsum \limits_{j=1}^{3}}%
%BeginExpansion
{\displaystyle\sum\limits_{j=1}^{3}}
%EndExpansion
\left\vert z^{j}-w^{j}\right\vert ^{2}\right)  ^{3}}.
\]

\end{corollary}

The formulas of Theorem \ref{b1} and Corollary \ref{b2} suggest that we can
define holomorphic linking for non compact CY manifolds$.$ In the case of
\ $\mathbb{C}^{3}$ we can choose the holomorphic form $\eta$ to be
$dz^{1}\wedge dz^{2}\wedge dz^{3}.$ In this case the formula for the
holomorphic linking becomes the complex version of the Gauss formula in
$\mathbb{R}^{3}$ as it appears in the Introduction.

\section{Geometric Interpretation of Holomorphic Linking}

In this section we will give a geometric interpretation of the holomorphic
linking, which is a direct complex generalization of the usual topological
definition of the Gauss linking number of two knots in $\mathbb{R}^{3}.$ We
will derive the geometric formula using the integral form of the holomorphic
linking via the Green kernel. The derivation is parallel to the real case but
it also uses the Leray residue theory as in \cite{FK}. The explicit formula of
\ Theorem \ref{igor19} was communicated to us by B. Khesin.

\subsection{Meromorphic forms with Prescribed Residues on Riemann Surfaces}

We will suppose that the Riemann surfaces $\Sigma_{1}$ and $\Sigma_{2}$ are of
genus $\geq1$ embedded in a CY threefold M and
\[
\Sigma_{1}\cap\Sigma_{2}=\emptyset.
\]
Let us fix two non-zero holomorphic forms $\theta_{i}$ on each of $\Sigma_{i}$
for $i=1,\ 2.$ According to \cite{AK} we can always find a very ample non
singular divisor $S_{1}$ containing $\Sigma_{1}$ since%
\[
\dim_{\mathbb{C}}\Sigma_{1}<\frac{1}{2}\dim_{\mathbb{C}}\text{M.}%
\]

\begin{proposition}
\label{igor17}Let $S_{1}$ be a non-singular hypersurface section of
M$\subset\mathbb{CP}^{m}$ which contains $\Sigma_{1}$, then there exist a
holomorphic two form $\omega_{1}$ on $S_{1}$ with a pole of order one along
$\Sigma_{1}$ such that Leray's residues of $\omega_{1}$ on $\Sigma_{1}$ is
equal to $\theta_{1}.$
\end{proposition}

\textbf{Proof:} We have the following exact sequences:
\[
0\rightarrow\Omega_{S_{1}}^{2}\rightarrow\Omega_{S_{1}}^{2}\left(  \Sigma
_{1}\right)  \overset{res}{\rightarrow}\Omega_{\Sigma_{1}}^{1}\rightarrow0
\]
and
\[
0\rightarrow H^{0}\left(  S_{1,}\Omega_{S_{1}}^{2}\right)  \rightarrow
H^{0}\left(  S_{1,}\Omega_{S_{1}}^{2}\left(  \Sigma_{1}\right)  \right)
\overset{res}{\rightarrow}%
\]%
\begin{equation}
\rightarrow H^{0}\left(  \Sigma_{1},\Omega_{\Sigma_{1}}^{1}\right)
\rightarrow H^{1}\left(  S_{1,}\Omega_{S_{1}}^{2}\right)  \rightarrow...
\label{es3}%
\end{equation}
where $\Omega_{S_{1}}^{2}\left(  \Sigma_{1}\right)  $ is the locally free
sheaf of holomorphic two forms on S$_{1}$ with a pole of order 1 on
$\Sigma_{1}.$ In order to deduce Proposition \ref{igor17}, we need to prove
that $H^{1}\left(  S_{1,}\Omega_{S_{1}}^{2}\right)  =0.$ The definition of a
CY manifold implies that $H_{1}($M,$\mathbb{C})=0$ and so by the Lefschetz
Theorem (see \cite{GH}) we conclude that $H_{1}(S_{1}$,$\mathbb{C})=0$. Then
the Poincare duality implies that $H^{3}(S_{1},\mathbb{C})=0.$ Hodge theory
implies that $H^{1}\left(  S_{1,}\Omega_{S_{1}}^{2}\right)  =0.$ So the map
\[
H^{0}\left(  S_{1,}\Omega_{S_{1}}^{2}\left(  \Sigma_{1}\right)  \right)
\overset{res}{\rightarrow}H^{0}\left(  \Sigma_{1},\Omega_{\Sigma_{1}}%
^{1}\right)  \rightarrow0
\]
is surjective. Proposition \ref{igor17} is proved. $\blacksquare$

\begin{proposition}
\label{igor18}Suppose that we can represent $\Sigma_{1}$ as
\begin{equation}
S_{1}\cap S_{2}=\Sigma_{1}, \label{Cu1}%
\end{equation}
where $S_{1}$ and $S_{2}$ are non-singular hypersurface sections on M $.$
\textit{Then there exist a meromorphic three form }$\eta_{1}$ on M
\textit{with poles of order one along }$S_{1}$ and $S_{2}$\textit{\ such that
the double Leray residue of \ }$\eta_{1}$ is equal to $\theta_{1}.$
\end{proposition}

\textbf{Proof:} We have the following exact sequences:
\[
0\rightarrow\Omega_{\text{M}}^{3}(S_{2})\rightarrow\Omega_{\text{M}}%
^{3}\otimes\mathcal{O}_{\text{M}}\left(  S_{1}\right)  \otimes\mathcal{O}%
_{\text{M}}\left(  S_{2}\right)  \overset{res}{\rightarrow}\Omega_{S_{1}}%
^{2}(S_{2})\rightarrow0
\]
and
\[
0\rightarrow H^{0}\left(  \text{M},\Omega_{\text{M}}^{3}(S_{2})\right)
\rightarrow H^{0}\left(  M,\Omega_{\text{M}}^{3}\otimes\mathcal{O}_{\text{M}%
}\left(  S_{1}\right)  \otimes\mathcal{O}_{\text{M}}\left(  S_{2}\right)
\right)  \overset{res}{\rightarrow}%
\]%
\begin{equation}
\overset{res}{\rightarrow}H^{0}\left(  S_{1},\Omega_{S_{1}}^{2}(S_{2})\right)
\rightarrow H^{1}\left(  \text{M},\Omega_{\text{M}}^{3}(S_{2})\right)
\rightarrow... \label{es4}%
\end{equation}
Since $S_{2}$ is a hypersurface section of CY manifold M, we can conclude from
the Kodaira vanishing Theorem that:
\begin{equation}
H^{1}\left(  \text{M},\Omega_{\text{M}}^{3}(S_{2})\right)  =0. \label{es5}%
\end{equation}
Combining $\left(  \ref{es4}\right)  $ and $\left(  \ref{es5}\right)  $, we
deduce that the Leray's residue map:
\[
H^{0}\left(  \text{M},\Omega_{\text{M}}^{3}\otimes\mathcal{O}_{\text{M}%
}\left(  S_{1}\right)  \otimes\mathcal{O}_{\text{M}}\left(  S_{2}\right)
\right)  \overset{res}{\rightarrow}H^{0}\left(  S_{1},\Omega_{S_{1}}^{2}%
(S_{2})\right)
\]
is surjective. This fact combined with Proposition \ref{igor17} implies
Proposition \ref{igor18}. $\blacksquare$

From now on we will consider curves on CY threefold M which are represented by
$\left(  \ref{Cu1}\right)  $.

\subsection{Geometric Formula for Holomorphic Linking}

Now we will express the holomorphic linking of $(\Sigma_{1},\theta_{1})$ and
$(\Sigma_{2},\theta_{2})$ as a sum of residues of certain meromorphic one form
over the intersection points of $S_{1}$ and $\Sigma_{2}.$ In \cite{HR} and
\cite{HR1} the expression $\left(  \ref{BK96}\right)  $ is interpreted via
polar homologies.

\begin{theorem}
\label{igor19}Let $\omega_{1}$ be a meromorphic two form defined as in
Proposition \ref{igor17} such that
\[
res_{\text{S}_{1}}\omega_{1}=\theta_{1}%
\]
Then the following formula is true:
\begin{equation}
\#\left(  (\Sigma_{1},\theta_{1}),(\Sigma_{2},\theta_{2})\right)
=\underset{x\in S_{1}\cap\Sigma_{2}}{\sum}\frac{\omega_{1}(x)\wedge\theta
_{2}(x)}{\eta(x)}, \label{BK96}%
\end{equation}
where $\omega_{1}(x)\wedge\theta_{2}(x)\in\wedge^{3}\left(  T_{x,\text{M}%
}^{1,0}\right)  ^{\ast}\approxeq\Omega_{x,\text{M}}^{3}.$
\end{theorem}

\textbf{Proof:} Let T$_{\varepsilon}(\Sigma_{1})$ be a tubular neighborhood of
$\ \Sigma_{1}\times\Sigma_{2}$ in $S_{1}\times\Sigma_{2}.$ From the definition
of Leray's residue and formula $\left(  \ref{Ig}\right)  $ we deduce that
\[
\#((\Sigma_{1},\theta_{1}),(\Sigma_{2},\theta_{2}))=\int_{\Sigma_{1}%
\times\Sigma_{2}}\left(  \overline{\partial}^{-1}\left(  \overline{\eta
}_{\Delta}\right)  \right)  \wedge\pi_{1}^{\ast}(\theta_{1})\wedge\pi
_{2}^{\ast}(\theta_{2})=
\]%
\begin{equation}
\underset{\varepsilon\rightarrow0}{\lim}\int_{\partial T_{\varepsilon}%
(\Sigma_{1})\times\Sigma_{2}}\left(  \overline{\partial}^{-1}(\overline{\eta
}_{\Delta})\right)  \wedge\pi_{1}^{\ast}(\omega_{1})\wedge\pi_{2}^{\ast
}(\theta_{2}), \label{Dia}%
\end{equation}
where $\partial T_{\varepsilon}(\Sigma_{1})$ is the boundary of
$T_{\varepsilon}(\Sigma_{1})$ in S$_{1}.$ Stokes' Theorem implies that
\[
\underset{\varepsilon\rightarrow0}{\lim}\int_{\partial T_{\varepsilon}%
(\Sigma_{1})\times\Sigma_{2}}\left(  \overline{\partial}^{-1}(\overline{\eta
}_{\Delta})\right)  \wedge\pi_{1}^{\ast}(\omega_{1})\wedge\pi_{2}^{\ast
}(\theta_{2})=
\]%
\[
\int_{\text{S}_{1}\text{ }\times\Sigma_{2}\text{$-$ }\Delta}d\left(
\frac{\psi_{\Delta}}{\pi_{1}^{\ast}(\eta)}\right)  \wedge\pi_{1}^{\ast}%
(\omega_{1})\wedge\pi_{2}^{\ast}(\theta_{2})-
\]%
\[
\int_{\text{S}_{1}\text{ }\times\Sigma_{2}}d\left(  \frac{\partial g_{\Delta}%
}{\pi_{1}^{\ast}(\eta)}\right)  \wedge\pi_{1}^{\ast}(\omega_{1})\wedge\pi
_{2}^{\ast}(\theta_{2})=
\]%
\[
\int_{\text{S}_{1}\text{ }\times\Sigma_{2}}\delta_{\Delta}\wedge\frac{\pi
_{1}^{\ast}(\omega_{1})\wedge\pi_{2}^{\ast}(\theta_{2})}{\pi_{1}^{\ast}(\eta
)}=\underset{x\in S_{1}\cap\Sigma_{2}}{\sum}\frac{\omega_{1}(x)\wedge
\theta_{2}(x)}{\eta(x)}.
\]
Theorem \ref{igor19} is proved. $\blacksquare$

One can also deduce an alternative version of the formula $\left(
\ref{BK96}\right)  $ using meromorphic three forms $\eta_{1}$ and $\eta
_{2}\,\ $with poles along $S_{1}$ and $S_{2}$ constructed in Proposition
\ref{igor18}. Then the ratio $\frac{\eta_{1}}{\eta}$ is a meromorphic section
on M of the line bundle $\mathcal{O}_{\text{M}}(S_{1})$ which we restrict to
$\Sigma_{2}$ and multiply by the holomorphic form $\theta_{2}.$ The result is
a meromorphic one form
\[
\left(  \frac{\eta_{1}}{\eta}\left\vert _{\Sigma_{2}}\right.  \right)
\theta_{2}%
\]
whose residues at $x\in S_{1}\cap\Sigma_{2}$ are precisely the values%
\[
\frac{\omega_{1}(x)\wedge\theta_{2}(x)}{\eta(x)}.
\]
This implies

\begin{corollary}
\label{igor19a}Let $\eta_{1}$ be a meromorphic two form defined as in
Proposition \ref{igor18}. Then the following formula holds%
\begin{equation}
\#\left(  (\Sigma_{1},\theta_{1}),(\Sigma_{2},\theta_{2})\right)  =%
%TCIMACRO{\dsum \limits_{x\in S_{1}\cap\Sigma_{2}}}%
%BeginExpansion
{\displaystyle\sum\limits_{x\in S_{1}\cap\Sigma_{2}}}
%EndExpansion
\text{res}_{x}\left(  \frac{\eta_{1}}{\eta}\left\vert _{\Sigma_{2}}\right.
\right)  \theta_{2}. \label{IGOR}%
\end{equation}

\end{corollary}

We can explore the symmetry of the holomorphic linking and represent the other
Riemann surface $\Sigma_{2}$ as an intersection of two surfaces. Then we
obtain another version of formulas $\left(  \ref{BK96}\right)  $ and $\left(
\ref{IGOR}\right)  $ using the meromorphic "lifts" of $\theta_{2}$ instead of
$\theta_{1}.$

Theorem \ref{igor19} and Corollary \ref{igor19a} are the complex analogues of
the geometric form of the linking number of two knots in $\mathbb{R}^{3}.$ In
the real case we took a disk whose boundary is one of the knots. We counted
the number of points of intersection of the disk with the other knot taking
into account the orientation. The linking number of two knots is the
intersection number just described. In the complex case, instead of a disk, we
took a complex surface which contains one of the original curves. They summed
again over points of intersection of the surface with the other curve weighted
with the natural ratios of the forms associated to the geometric objects
involved in the construction. Since we always recover the same holomorphic
linking, the geometric realization in the complex case also does not depend on
the choice of the complex surface and the meromorphic two form with prescribed residue.

\section{Cohomological Interpretation of Holomorphic Linking}

In this section we will reformulate the analytic expressions for the
holomorphic linking in the language of homological algebra. The homological
interpretation of the analytic formulas for the holomorphic linking is based
on the notions of generalized Grothendieck and Serre classes $\mu($Y,$\theta$)
and $\lambda($Y,$\theta$) of subvariety Y in a projective variety X together
with top holomorphic form $\theta$ on Y. Our definitions of Grothendieck and
Serre classes are generalizations of the similar notions introduced by Atiyah
in \cite{A}. The Grothendieck class $\mu($Y,$\theta$) always exists by
definition. We will prove the existence of generalized Serre classes for the
embedding of a submanifold Y in a CY manifold M and the diagonal embedding of
M into M$\times$M. Then we derive homological expressions for the holomorphic
linking using the Yoneda pairing. These expressions make sense over arbitrary
field. They should be related to the height pairing in \cite{B}, \cite{Be} and
\cite{Bo}. We will also illustrate why the original Atiyah's notion of
Grothendieck and Serre classes can not be used in our setting. We will start
with the calculation of the local cohomology of the diagonal embedding and we
will recall some basic facts from the general theory of local cohomology; a
more detailed account of the theory can be found in \cite{H}, \cite{H1} and
\cite{Hart};

\subsection{Local Cohomology of the Diagonal}

The key property of the local cohomology is the existence of the analogue of
the Mayer-Vietoris exact sequence, namely, let $\mathcal{F}$ be a sheaf on a
complex manifold X and Y be a closed submanifold, then the local cohomology
groups $H_{\text{Y}}^{k}(X,\mathcal{F})$ will satisfy the following exact
sequence:
\[
\rightarrow H^{k-1}(X,\mathcal{F)}\rightarrow H^{k-1}(X\text{ }-\text{
}Y,\mathcal{F})\overset{\delta_{k-1}}{\rightarrow}H_{\text{Y}}^{k}%
(X,\mathcal{F)}\rightarrow H^{k}(X,\mathcal{F)\rightarrow}%
\]
The boundary map $\delta$ in the exact sequence can be interpreted as an
\textquotedblright explicit\textquotedblright\ construction and a
generalization of distributions of \textquotedblright boundary
values\textquotedblright\ of sections of the coherent sheaf $\mathcal{F}.$

One of the main results that was proved in \cite{H1} is that the local
cohomology groups are related to the functor $Ext.$ The computation of the
local cohomology group $H_{Y}^{k-1}(X,\mathcal{F)}$ is done by using the
following spectral sequence:

First we define the sheaf of extensions in a standard way for any two coherent
sheaves $\mathcal{F}$ and $\mathcal{G}$ on X \cite{SGA2}. We will denote this
sheaf by \underline{$Ext$}$_{\mathcal{O}_{X}}^{i}(\mathcal{G},\mathcal{F)}$
and its global sections by $Ext_{\mathcal{O}_{X}}^{i}(\mathcal{G}%
,\mathcal{F)}.$ Then we define the sheaf $\underline{\mathcal{H}}_{Y}%
^{i}(\mathcal{F)}$ as the following projective limit:
\begin{equation}
\underline{\mathcal{H}}_{\text{Y}}^{i}(\mathcal{F)=}\underset{n}%
{\underrightarrow{\lim}}\underline{Ext}_{\mathcal{O}_{X}}^{i}(\mathcal{O}%
_{n,Y},\mathcal{F)} \label{Grot}%
\end{equation}
where $\mathcal{O}_{n,\text{Y}}=\mathcal{O}_{\text{X}}/I_{\text{Y}}^{n}$ and
$I_{\text{Y}}$ is the ideal sheaf in $\mathcal{O}_{\text{X}}$ that defines Y.
According to \cite{SGA2} Exp. I page 9 and Exp. II page 2 spectral sequence
with the initial term
\[
E_{2}^{p,q}=H^{p}(X,\underline{\mathcal{H}}_{\text{Y}}^{q}(\mathcal{F)})
\]
converges to $H_{\text{Y}}^{p+q}(X,\mathcal{F}).$

The following theorem can be found in \cite{Ma}:

\begin{theorem}
\label{ma} Let M be a non singular algebraic manifold. Let $\Delta\subset
$M$\times$M be the diagonal. Let $I_{\Delta}$ be the ideal sheaf of the
diagonal, then $I_{\Delta}/I_{\Delta}^{2}$ \ is isomorphic to the cotangent
sheaf $\Omega_{\text{M}}^{1}$ and, moreover,
\[
I_{\Delta}^{k}/I_{\Delta}^{k+1}\backsimeq S^{k}(I_{\Delta}/I_{\Delta}^{2})
\]
is a locally free $\mathcal{O}_{\text{M}\times\text{M}}/I_{\Delta}$ module.
\end{theorem}

We will use this fact to establish the following:

\begin{theorem}
$\label{ma1}$ Let M be a Calabi Yau manifold and let $T_{\text{M}}$ be the
tangent bundle of M, then
\[
H^{0}(\text{M},S^{k}(T_{\text{M}}))=H^{0}(\text{M,}S^{k}(\Omega_{\text{M}}%
^{1}))=0
\]
for all $k>0.$
\end{theorem}

\textbf{Proof:} When $k=1$ Theorem \ref{ma1} follows from the isomorphism
$T_{\text{M}}\overset{\lrcorner\eta_{\text{M}}}{\approxeq}\Omega_{\text{M}%
}^{n-1}.$ Thus we have
\[
H^{0}(\text{M},T_{\text{M}})=H^{0}(\text{M},\Omega_{\text{M}}^{n-1}).
\]
The definition of a CY manifold implies that $H^{0}(M,\Omega_{\text{M}}%
^{n-1})=0.$ Thus $H^{0}($M$,T_{\text{M}})=H^{0}($M,$\Omega_{\text{M}}^{1})=0.$

Bochner's principle for the Ricci flat K\"{a}hler metric implies that if
$\phi$ is any holomorphic tensor on a CY manifold, then it is parallel with
respect to the Levi Cevita connection with respect to Yau's metric. $\left(
\text{See }\cite{besse}.\right)  $ We\ also know from \cite{besse} that for a
simply connected CY manifold the holonomy group of the CY metric is $SU(n).$
These two facts imply that the globally defined holomorphic symmetric one
forms are obtained from the $SU(n)$ invariant k symmetric tensors at one point
by parallel transportation. So we have the following equality:
\[
S^{k}(\mathbb{C}^{n})^{SU(n)}=H^{0}(M,S^{k}(T_{\text{M}}))=H^{0}%
(\text{M,}S^{k}(\Omega_{\text{M}}^{1})).
\]
Since $S^{k}(\mathbb{C}^{n})^{SU(n)}=0,$ we conclude from Theorem \ref{ma}
that
\[
H^{0}(\text{M},S^{k}(T_{\text{M}}))=H^{0}(\text{M,}S^{k}(\Omega_{\text{M}}%
^{1}))=H^{0}(\text{M}\times\text{M},I_{\Delta}^{k}/I_{\Delta}^{k+1})=0
\]
for $k\geq1.$Theorem \ref{ma1} is proved. $\blacksquare$

Next we will compute some local cohomology groups that will be needed in the
construction and the cohomological interpretation of the Green current
associated with subvarieties. By the definition $\left(  \ref{Grot}\right)
$:
\[
\underline{\mathcal{H}}_{\Delta}^{i}(\mathcal{O}_{\text{M}\times\text{M}%
}\mathcal{)=}\underset{k}{\underrightarrow{\lim}}\underline{Ext}%
_{\mathcal{O}_{\text{M}\times\text{M}}}^{i}(\mathcal{O}_{k,\Delta}%
,\mathcal{O}_{\text{M}\times\text{M}}\mathcal{)},
\]
where\textit{\ }$\mathcal{O}_{k,\Delta}:=\mathcal{O}_{\text{M}\times\text{M}%
}/I_{\Delta}^{k}$ and $I_{\Delta}$ is the ideal sheaf of $\Delta\subset
$M$\times$M$.$

\begin{theorem}
\label{lc1}$\underline{\mathcal{H}}_{\Delta}^{i}(\mathcal{O}_{\text{M}%
\times\text{M}}\mathcal{)}=0$ for $i\neq n$ and $\underline{\mathcal{H}%
}_{\Delta}^{n}(\mathcal{O}_{\text{M}\times\text{M}}\mathcal{)}\neq0.$
\end{theorem}

\textbf{Proof:} Let $\mathcal{U}$ be an open affine set in M$\times$M. We will
prove by induction on $k$ that we have:
\begin{equation}
\underline{Ext}_{\mathcal{O}_{\text{M}\times\text{M}}}^{i}(\mathcal{O}%
_{\text{M}\times\text{M}}/I_{\Delta}^{k},\mathcal{O}_{\text{M}\times\text{M}%
})=0\text{ } \label{lcs1}%
\end{equation}
for $i\neq n$\ and
\begin{equation}
\underline{Ext}_{\mathcal{O}_{\text{M}\times\text{M}}}^{n}(\mathcal{O}%
_{\text{M}\times\text{M}}/I_{\Delta}^{k},\mathcal{O}_{\text{M}\times\text{M}%
})\neq0. \label{lcs2}%
\end{equation}
The diagonal $\Delta\subset$M$\times$M is a smooth algebraic variety.
Therefore it is a local complete intersection in M$\times$M. This we have%
\[
\mathcal{O}_{\Delta}=\mathcal{O}_{\text{M}\times\text{M}}/(f_{1},...,f_{n})
\]
where $\Delta$ is locally defined by the regular sequence $f_{1},...,f_{n}$ of
analytic functions in $\mathcal{O}_{\text{M}\times\text{M}}$ by
\[
f_{1}=...=f_{n}=0.
\]
In Chapter 5 "Residues", paragraph "Kozul Complex and Its Applications" of
\cite{GH} it is proved that for any regular sequence $(f_{1},...,f_{p})$ in
the local ring $\mathcal{O}_{\text{N}}$ of any complex manifold N we have:%
\begin{equation}
\underline{Ext}_{\mathcal{O}_{\text{N}}}^{j}(\mathcal{O}_{\text{N}}%
/(f_{1},...,f_{p}),\mathcal{O}_{\text{N}})=\left\{
\begin{array}
[c]{cc}%
\mathcal{O}_{\text{N}}/(f_{1},...,f_{p}) & j=p\\
0 & j\neq p
\end{array}
\right.  . \label{lcs2a}%
\end{equation}
Thus we proved $\left(  \ref{lcs1}\right)  $ and $\left(  \ref{lcs2}\right)  $
for $k=1.$

In order to proceed with the induction on $k$, we assume that $\left(
\ref{lcs1}\right)  $ and $\left(  \ref{lcs2}\right)  $ are true for $k=m.$ The
following exact sequence:
\[
0\rightarrow I_{\Delta}^{m}/I_{\Delta}^{m+1}\rightarrow\mathcal{O}%
_{\text{M}\times\text{M}}/I_{\Delta}^{m+1}\rightarrow\mathcal{O}%
_{\text{M}\times\text{M}}/I_{\Delta}^{m}\rightarrow0
\]
implies the long exact sequence:%
\[
...\rightarrow\underline{Ext}_{O_{\text{M}\times\text{M}}}^{i}(\mathcal{O}%
_{\text{M}\times\text{M}}/I_{\Delta}^{m},\mathcal{O}_{\text{M}\times\text{M}%
})\text{ }\rightarrow\underline{Ext}_{\mathcal{O}_{\text{M}\times\text{M}}%
}^{i}(\mathcal{O}_{\text{M}\times\text{M}}/I_{\Delta}^{m+1},\mathcal{O}%
_{\text{M}\times\text{M}})\rightarrow
\]%
\begin{equation}
\rightarrow\underline{Ext}_{\mathcal{O}_{\text{M}\times\text{M}}}%
^{i}(I_{\Delta}^{m}/I_{\Delta}^{m+1},\mathcal{O}_{\text{M}\times\text{M}%
})\rightarrow... \label{LC3}%
\end{equation}
Since by Theorem \ref{ma} $I_{\Delta}^{m}/I_{\Delta}^{m+1}$ is a free
$\mathcal{O}_{\Delta}$ module, we obtain from
\begin{equation}
\underline{Ext}_{\mathcal{O}_{\text{M}\times\text{M}}}^{j}(I_{\Delta}%
^{m}/I_{\Delta}^{m+1},\mathcal{O}_{\text{M}\times\text{M}})=\left\{
\begin{array}
[c]{cc}%
I_{\Delta}^{m}/I_{\Delta}^{m+1} & j=n\\
0 & j\neq n
\end{array}
\right.  . \label{LC1}%
\end{equation}
From the induction hypothesis, $\left(  \ref{LC1}\right)  $ and the long exact
sequence $\left(  \ref{LC3}\right)  $, we can deduce $\left(  \ref{lcs1}%
\right)  $ and $\left(  \ref{lcs2}\right)  $ for any $k.$ Thus Theorem
\ref{lc1} follows from the definition of the sheaves $\underline{\mathcal{H}%
}_{\Delta}^{i}(\mathcal{O}_{\text{M}\times\text{M}}\mathcal{)}$ given by
$\left(  \ref{Grot}\right)  .\blacksquare$

\begin{corollary}
\label{zero}The following formula%
\[
H^{0}(\text{M}\times\text{M},\underline{Ext}_{\mathcal{O}_{\text{M}%
\times\text{M}}}^{k}(I_{\Delta}^{m}/I_{\Delta}^{m+1},\mathcal{O}%
_{\text{M}\times\text{M}}))=
\]%
\begin{equation}
Ext_{\mathcal{O}_{\text{M}\times\text{M}}}^{k}(I_{\Delta}^{m}/I_{\Delta}%
^{m+1},\mathcal{O}_{\text{M}\times\text{M}})=0 \label{Z}%
\end{equation}
holds for $k\leq n$ and $m>0$.
\end{corollary}

\textbf{Proof:} Formula $\left(  \ref{LC1}\right)  $ implies Corollary
\ref{zero} for $k\neq n=\dim_{\mathbb{C}}$M. Thus $\left(  \ref{Z}\right)  $
is proved for any $k\neq n.$

Suppose that $k=n.$ Formula $\left(  \ref{LC1}\right)  ,$ Theorems \ref{ma}
and \ref{ma1} imply that for $m>0$%
\[
Ext_{\mathcal{O}_{\text{M}\times\text{M}}}^{n}(I_{\Delta}^{m}/I_{\Delta}%
^{m+1},\mathcal{O}_{\text{M}\times\text{M}})=H^{0}\left(  \text{M}%
\times\text{M},\underline{Ext}_{\mathcal{O}_{\text{M}\times\text{M}}}%
^{n}(I_{\Delta}^{m}/I_{\Delta}^{m+1},\mathcal{O}_{\text{M}\times\text{M}%
})\right)  =
\]%
\begin{equation}
H^{0}(\text{M}\times\text{M},I_{\Delta}^{m}/I_{\Delta}^{m+1})=0. \label{Z2}%
\end{equation}
So Corollary \ref{zero} is proved. $\blacksquare$

\begin{theorem}
\label{igor11} $H_{\Delta}^{k}($M$\times$M,$\mathcal{O}_{\text{M}%
\times\text{M}})=0$ for $k\neq n$ and dim$_{\mathbb{C}}H_{\Delta}^{n}%
($M$\times$M,$\mathcal{O}_{\text{M}\times\text{M}})=1.$
\end{theorem}

\textbf{Proof:} By Theorem \ref{lc1} $\underline{\mathcal{H}}_{\Delta}%
^{q}(\mathcal{O}_{\text{M}\times\text{M}}\mathcal{)}\mathbf{=}0$ for $q\neq0$.
Therefore $H_{\Delta}^{k}($M$\times$M,$\mathcal{O}_{\text{M}\times\text{M}%
})=0$ for $k\neq n.$

According to \cite{SGA2} $H_{\Delta}^{n}($M$\times$M,$\mathcal{O}%
_{\text{M}\times\text{M}})$ is obtained by the spectral sequence with a first
term $E_{2}^{p,q}=H^{p}($M$\times$M$,\underline{\mathcal{H}}_{\Delta}%
^{q}(\mathcal{O}_{\text{M}\times\text{M}}\mathcal{)}).$ Thus we have for
$p+q=n$
\[
E_{2}^{p,q}=H^{p}(\text{M}\times\text{M},\underline{\mathcal{H}}_{\Delta}%
^{q}(\mathcal{O}_{\text{M}\times\text{M}}\mathcal{)})\Longrightarrow
H_{\Delta}^{n}(\text{M}\times\text{M},\mathcal{O}_{\text{M}\times\text{M}}).
\]
From Theorem \ref{lc1} and the definition of the cohomology group $H_{\Delta
}^{n}($M$\times$M$,\mathcal{O}_{\text{M}\times\text{M}})$ by the spectral
sequence we get that
\[
H_{\Delta}^{n}(\text{M}\times\text{M},\mathcal{O}_{\text{M}\times\text{M}%
})\backsimeq H^{0}(\text{M}\times\text{M},\underline{\mathcal{H}}_{\Delta}%
^{n}(\mathcal{O}_{\text{M}\times\text{M}}\mathcal{)})\backsimeq
\]%
\[
\backsimeq H^{0}\left(  \text{M}\times\text{M},\underset{k}{\underrightarrow
{\lim}}\underline{Ext}_{\mathcal{O}_{\text{M}\times\text{M}}}^{n}%
(\mathcal{O}_{\text{M}\times\text{M}}/\left(  I_{\Delta}\right)
^{k},\mathcal{O}_{\text{M}\times\text{M}}\mathcal{)}\right)  .
\]
Thus the $n^{th}$ local cohomology $H_{\Delta}^{n}($M$\times$M$,\mathcal{O}%
_{\text{M}\times\text{M}})$ is the inductive limit of Cech cohomologies
\[
\underset{k}{\underrightarrow{\lim}}H^{0}\left(  \text{M}\times\text{M}%
,\underline{Ext}_{\mathcal{O}_{\text{M}\times\text{M}}}^{n}\left(
\mathcal{O}_{\text{M}\times\text{M}}/\left(  I_{\Delta}\right)  ^{k}%
,\mathcal{O}_{\text{M}\times\text{M}}\right)  \right)  .
\]
In \cite{SGA2} in Exp. II it is proved that
\[
H^{0}\left(  \text{M}\times\text{M},\underset{k}{\underrightarrow{\lim}%
}\underline{Ext}_{\mathcal{O}_{\text{M}\times\text{M}}}^{n}\left(
\mathcal{O}_{\text{M}\times\text{M}}/\left(  I_{\Delta}\right)  ^{k}%
,\mathcal{O}_{\text{M}\times\text{M}}\right)  \right)  =
\]%
\[
\underset{k}{\underrightarrow{\lim}}H^{0}\left(  \text{M}\times\text{M}%
,\underline{Ext}_{\mathcal{O}_{\text{M}\times\text{M}}}^{n}\left(
\mathcal{O}_{\text{M}\times\text{M}}/\left(  I_{\Delta}\right)  ^{k}%
,\mathcal{O}_{\text{M}\times\text{M}}\right)  \right)  .
\]
Thus we have the following isomorphism:
\begin{equation}
H_{\Delta}^{n}(\text{M}\times\text{M},\mathcal{O}_{\text{M}\times\text{M}%
})=\underset{k}{\underrightarrow{\lim}}H^{0}\left(  \text{M}\times
\text{M},\underline{Ext}_{\mathcal{O}_{\text{M}\times\text{M}}}^{n}%
(\mathcal{O}_{\text{M}\times\text{M}}/\left(  I_{\Delta}\right)
^{k},\mathcal{O}_{\text{M}\times\text{M}}\mathcal{)}\right)  . \label{W4}%
\end{equation}

\begin{lemma}
\label{mar2} The natural restriction maps%
\[
H^{0}(\text{M}\times\text{M},\underline{Ext}_{\mathcal{O}_{\text{M}%
\times\text{M}}}^{n}(\mathcal{O}_{\text{M}\times\text{M}}/\left(  I_{\Delta
}\right)  ^{k+1},\mathcal{O}_{\text{M}\times\text{M}}))\rightarrow
H^{0}(\text{M}\times\text{M},\underline{Ext}_{\mathcal{O}_{\text{M}%
\times\text{M}}}^{n}(\mathcal{O}_{\text{M}\times\text{M}}/\left(  I_{\Delta
}\right)  ^{k},\mathcal{O}_{\text{M}\times\text{M}}))
\]
are isomorphisms. Thus the following formula is true:
\[
\dim_{\mathbb{C}}H^{0}(\text{M}\times\text{M},\underline{Ext}_{\mathcal{O}%
_{\text{M}\times\text{M}}}^{n}(\mathcal{O}_{\text{M}\times\text{M}}/\left(
I_{\Delta}\right)  ^{k},\mathcal{O}_{\text{M}\times\text{M}}))=1\text{
}for\text{ }k\geq0.
\]

\end{lemma}

\textbf{Proof:} The proof of Lemma \ref{mar2} is by induction. It is based on
the following long exact sequence of sheaves:
\[
0\rightarrow\underline{Ext}_{\mathcal{O}_{\text{M}\times\text{M}}}%
^{n}(\mathcal{O}_{\text{M}\times\text{M}}/I_{\Delta}^{k},\mathcal{O}%
_{\text{M}\times\text{M}})\rightarrow\underline{Ext}_{\mathcal{O}%
_{\text{M}\times\text{M}}}^{n}(\mathcal{O}_{\text{M}\times\text{M}}/I_{\Delta
}^{k+1},\mathcal{O}_{\text{M}\times\text{M}})\rightarrow
\]%
\begin{equation}
\rightarrow\underline{Ext}_{\mathcal{O}_{\text{M}\times\text{M}}}%
^{n}(I_{\Delta}^{k}/I_{\Delta}^{k+1},\mathcal{O}_{\text{M}\times\text{M}%
})\rightarrow0. \label{lces}%
\end{equation}
From Theorem \ref{ma} states that
\[
\underline{Ext}_{\mathcal{O}_{\text{M}\times\text{M}}}^{n}(I_{\Delta}%
^{k}/I_{\Delta}^{k+1},\mathcal{O}_{\text{M}\times\text{M}})\backsimeq
I_{\Delta}^{k}/I_{\Delta}^{k+1}.
\]
From $\left(  \ref{LC1}\right)  $ we have:%
\[
H^{0}(\text{M}\times\text{M},\underline{Ext}_{\mathcal{O}_{\text{M}%
\times\text{M}}}^{n}(I_{\Delta}^{k}/I_{\Delta}^{k+1},\mathcal{O}%
_{\text{M}\times\text{M}}))=
\]%
\[
H^{0}(\text{M}\times\text{M},I_{\Delta}^{k}/I_{\Delta}^{k+1})=H^{0}%
(\Delta,S^{k}(\Omega_{\Delta}^{1}))=H^{0}(\text{M},S^{k}(\Omega_{\text{M}}%
^{1}))=0.
\]
Combining this fact with the exact sequence $\left(  \ref{lces}\right)  $ we
obtain that:%
\[
H^{0}(\text{M}\times\text{M},\underline{Ext}_{\mathcal{O}_{\text{M}%
\times\text{M}}}^{n}(\mathcal{O}_{\text{M}\times\text{M}}/I_{\Delta}%
^{k},\mathcal{O}_{\text{M}\times\text{M}})=
\]%
\begin{equation}
H^{0}(\text{M}\times\text{M},\underline{Ext}_{\mathcal{O}_{\text{M}%
\times\text{M}}}^{n}(\mathcal{O}_{\text{M}\times\text{M}}/I_{\Delta}%
^{k+1},\mathcal{O}_{\text{M}\times\text{M}})) \label{isomk}%
\end{equation}
for $k\geq1.$ The isomorphism
\[
\underline{Ext}_{\mathcal{O}_{\text{M}\times\text{M}}}^{n}(\mathcal{O}%
_{\text{M}\times\text{M}}(\mathcal{U})/I_{\Delta},\mathcal{O}_{\text{M}%
\times\text{M}}(\mathcal{U}))\backsimeq\mathcal{O}_{\text{M}\times\text{M}%
}(\mathcal{U})/I_{\Delta},
\]
implies that%
\[
\dim_{\mathbb{C}}H^{0}\left(  \text{M}\times\text{M},\underline{Ext}%
_{\mathcal{O}_{\text{M}\times\text{M}}}^{n}\left(  \mathcal{O}_{\text{M}%
\times\text{M}}/I_{\Delta},\mathcal{O}_{\text{M}\times\text{M}}\right)
\right)  =
\]%
\begin{equation}
\dim_{\mathbb{C}}H^{0}\left(  \text{M}\times\text{M},\mathcal{O}%
_{\text{M}\times\text{M}}/I_{\Delta}\right)  =\dim_{\mathbb{C}}H^{0}\left(
\Delta,\mathcal{O}_{\Delta}\right)  =\dim_{\mathbb{C}}H^{0}\left(
\text{M},\mathcal{O}_{\text{M}}\right)  =1. \label{isom}%
\end{equation}
Lemma \ref{mar2} follows directly from $\left(  \ref{isomk}\right)  $
with\ $\left(  \ref{isom}\right)  $. $\blacksquare$

Theorem \ref{igor11} follows from the isomorphism $\left(  \ref{W4}\right)  $
and Lemma \ref{mar2}. $\blacksquare$

\begin{corollary}
\label{igor11a}$H_{\Delta}^{n}($M$\times$M,$\mathcal{O}_{\text{M}%
\times\text{M}})\approxeq Ext_{\mathcal{O}_{\text{M}\times\text{M}}}%
^{n}\left(  \mathcal{O}_{\text{M}\times\text{M}}/I_{\Delta},\mathcal{O}%
_{\text{M}\times\text{M}}\right)  \approxeq\mathbb{C}$ $.$
\end{corollary}

\textbf{Proof:} Corollary \ref{igor11a} follows from the proof of Lemma
\ref{mar2} and the definition of $H_{\Delta}^{n}($M$\times$M,$\mathcal{O}%
_{\text{M}\times\text{M}}).$ $\blacksquare$

From the Grothendieck duality it follows that $Ext_{\mathcal{O}_{\text{M}%
\times\text{M}}}^{n}\left(  \mathcal{O}_{\text{M}\times\text{M}}/I_{\Delta
},\mathcal{O}_{\text{M}\times\text{M}}\right)  $ can be identified with the
restriction of the space of antiholomorphic $n$ forms on M$\times$M on the
diagonal $\Delta.$ Thus from the definition of the local cohomology it follows
that $H_{\Delta}^{n}($M$\times$M,$\mathcal{O}_{\text{M}\times\text{M}})$ is
generated by $\overline{\eta}_{\Delta}.$

\subsection{Definition of the Grothendieck and Serre Classes}

Let $E\rightarrow$X be a vector bundle over a complex manifold X of a complex
dimension n. Serre showed that the pairing
\[
H^{p}(\text{X},E)\times H^{n-p}(\text{X},E^{\ast}\otimes\Omega_{X}%
^{n})\rightarrow\mathbb{C}%
\]
given by the integration over X of the corresponding pointwise pairing of the
cohomology classes is nondegenerate.

Let $\mathcal{F},$ $\mathcal{H}$ and $\mathcal{G}$ be three coherent sheaves
on X. The Yoneda product
\[
Ext_{\mathcal{O}_{\text{X}}}^{p}(\mathcal{F},\mathcal{G})\times
Ext_{\mathcal{O}_{\text{X}}}^{q}(\mathcal{G},\mathcal{H})\rightarrow
Ext_{\mathcal{O}_{\text{X}}}^{p+q}(\mathcal{F},\mathcal{H})
\]
is defined in a natural way by the composition of long exact sequences. It is
a well known fact that if $\mathcal{F}$ is the sheaf of holomorphic function
on X denoted by $\mathcal{O}_{\text{X}}$, then
\[
Ext^{p}(\mathcal{O}_{\text{X}},\mathcal{F})\mathcal{\approxeq}H^{p}%
(\text{X},\mathcal{F}).
\]
Grothendieck proved that the map
\begin{equation}
H^{p}(\text{X},\mathcal{F})\times Ext_{\mathcal{O}_{\text{X}}}^{n-p}%
(\mathcal{F},\Omega_{\text{X}}^{n})\rightarrow H^{n}(\text{X},\Omega
_{\text{X}}^{n})\mathcal{\approxeq}\mathbb{C} \label{W7}%
\end{equation}
given by the Yoneda pairing is non-degenerate. This pairing is called the
Grothendieck duality.

If $\mathcal{E}$ is a locally free sheaf, i.e. $\mathcal{E}$ is the sheaf of
sections of a vector bundle $E\rightarrow$X$,$ then Grothendieck's duality
implies Serre's duality by using the following isomorphism
\begin{equation}
Ext_{\mathcal{O}_{\text{X}}}^{n-p}(\mathcal{E},\Omega_{\text{X}}^{n})\approxeq
H^{n-p}(\text{X},\mathcal{E}^{\ast}\otimes\Omega_{\text{X}}^{n}). \label{W6}%
\end{equation}
Indeed by $\left(  \ref{W7}\right)  $
\begin{equation}
Ext_{\mathcal{O}_{\text{X}}}^{n-p}(\mathcal{E},\Omega_{\text{X}}^{n})\approxeq
H^{p}(\text{X},\mathcal{E})^{\ast}. \label{W8}%
\end{equation}
On the other hand we know that
\begin{equation}
H^{p}(\text{X},\mathcal{E})^{\ast}\approxeq H^{n-p}(\text{X},\mathcal{E}%
^{\ast}\otimes\Omega_{\text{X}}^{n}). \label{W9}%
\end{equation}
Combining $\left(  \ref{W8}\right)  $ and $\left(  \ref{W9}\right)  $ we get
$\left(  \ref{W6}\right)  .$

For any submanifold Y$\subset$X of codimension $m$, we will denote by
$I_{\text{Y }}$ the ideal sheaf consisting of functions vanishing on Y. Then
the definition of the sheaf $I_{\text{Y}}$ implies that the quotient sheaf
$\mathcal{O}_{\text{X}}$/$I_{\text{Y}}$ is naturally identified with the
structure sheaf $\mathcal{O}_{\text{Y}}$ extended by zero on X $-$Y.

Let $($Y,$\theta)$ be a pair, where Y is a submanifold in X of codimension $m$
and $\theta$ is a holomorphic $n-m=\dim_{\mathbb{C}}$Y form on Y. By applying
the Grothendieck duality twice one can deduce that
\begin{equation}
H^{0}(\text{Y},\Omega_{\text{Y}}^{n-m})\approxeq\left(  Ext_{\mathcal{O}%
_{\text{X}}}^{n-m}(\mathcal{O}_{\text{Y}},\mathcal{O}_{\text{Y}})\right)
^{\ast}\mathcal{\approxeq}Ext_{\mathcal{O}_{\text{X}}}^{m}(\mathcal{O}%
_{\text{Y}},\Omega_{\text{X}}^{n}). \label{W5}%
\end{equation}
Thus $\left(  \ref{W5}\right)  $ defines a canonical isomorphism
\begin{equation}
\mu:H^{0}(\text{Y},\Omega_{\text{Y}}^{n-m})\mathcal{\approxeq}Ext_{\mathcal{O}%
_{\text{X}}}^{m}(\mathcal{O}_{\text{Y}},\Omega_{\text{X}}^{n}). \label{W5a}%
\end{equation}

\begin{definition}
\label{igor25}We define the Grothendieck class $\mu$(Y,$\theta$) of the pair
(Y,$\theta)$ in X as the image of
\[
\theta\in H^{0}(\text{Y},\Omega_{\text{Y}}^{n-m})
\]
under the canonical map $\mu.$
\end{definition}

\begin{proposition}
\label{GR}The Grothendieck class $\mu(\Delta,\pi_{1}^{\ast}(\eta)|_{\Delta})$
can be canonically identified with the class of the antiholomorphic Dirac
current $\overline{\eta}_{\Delta}$ given by Definition \ref{andrey48}.
\end{proposition}

\textbf{Proof: }From the proof of Lemma \ref{mar2} it follows that there exist
canonical identifications:
\[
Ext_{\mathcal{O}_{\text{M}\times\text{M}}}^{n}(\mathcal{O}_{\text{M}%
\times\text{M}}/I_{\Delta}^{k},\mathcal{O}_{\text{M}\times\text{M}})\approxeq
Ext_{\mathcal{O}_{\text{M}\times\text{M}}}^{n}(\mathcal{O}_{\text{M}%
\times\text{M}}/I_{\Delta},\mathcal{O}_{\text{M}\times\text{M}})\approxeq
Ext_{\mathcal{O}_{\text{M}\times\text{M}}}^{n}(\mathcal{O}_{\Delta
},\mathcal{O}_{\text{M}\times\text{M}})
\]
for $k>0.$ According to Corollary \ref{igor11a} the following canonical
identification
\[
Ext_{\mathcal{O}_{\text{M}\times\text{M}}}^{n}(\mathcal{O}_{\Delta
},\mathcal{O}_{\text{M}\times\text{M}})=H_{\Delta}^{n}(\text{M}\times
\text{M},\mathcal{O}_{\text{M}\times\text{M}})
\]
exists and by Theorem \ref{igor11}
\[
\dim_{\mathbb{C}}Ext_{\mathcal{O}_{\text{M}\times\text{M}}}^{n}(\mathcal{O}%
_{\Delta},\mathcal{O}_{\text{M}\times\text{M}})=\dim_{\mathbb{C}}H_{\Delta
}^{n}(\text{M}\times\text{M},\mathcal{O}_{\text{M}\times\text{M}})=1.
\]
At the end of Section\textbf{ }7.1\textbf{. }we identified the generator of
$H_{\Delta}^{n}($M$\times$M,$\mathcal{O}_{\text{M}\times\text{M}})$ with the
class of the antiholomorphic Dirac current $\overline{\eta}_{\Delta}.$ Thus
the Grothendieck class $\mu(\Delta,\pi_{1}^{\ast}(\eta))$ can be interpreted
as a local cohomology class and in particular it can be identified with the
class of the antiholomorphic Dirac current $\overline{\eta}_{\Delta}.$
$\blacksquare$

Since the Grothendieck class $\mu(\Delta,\pi_{1}^{\ast}(\eta)|_{\Delta})$ does
not depends of the pullback of the holomorphic form $\eta$ on M$\times$M we
will denoted it by $\mu(\Delta,\eta).$

\begin{remark}
\label{At0}The definition of the Grothendieck class $\mu$(Y) of a subvariety Y
in an algebraic variety X given by Atiyah in \cite{A} differs from ours. When
we replace the holomorphic form $\theta\in H^{0}($Y$,\Omega_{\text{Y}}^{n-m})$
by $1\in H^{0}($Y,$\mathcal{O}_{\text{Y}})$ we will get the Atiyah definition
of the Grothendieck class
\[
\mu(\text{Y})\in Ext_{\mathcal{O}_{\text{X}}}^{m}(\Omega_{\text{Y}}%
^{n-m},\Omega_{\text{X}}^{n}).
\]
In particular $\mu(\Delta)$ can be identified with the Dirac current
$\delta_{\Delta}.$
\end{remark}

Now we will define the Serre class $\lambda($Y,$\theta)$ of a pair
$($Y,$\theta)$ in X. The definition of the Serre class $\lambda($Y,$\theta)$
is possible under the assumption that the coboundary map
\begin{equation}
d_{m-1}:Ext_{\mathcal{O}_{\text{X}}}^{m-1}(I_{\text{Y}},\Omega_{\text{X}}%
^{n})\rightarrow Ext_{\mathcal{O}_{\text{X}}}^{m}(\mathcal{O}_{\text{Y}%
},\Omega_{\text{X}}^{n}) \label{LP4}%
\end{equation}
resulting from the exact sequence%
\[
0\rightarrow I_{\text{Y}}\rightarrow\mathcal{O}_{\text{X}}\rightarrow
\mathcal{O}_{\text{Y}}\rightarrow0
\]
is an isomorphism .

\begin{definition}
\label{GSC}Suppose that Y is a non singular subvariety in the projective
smooth variety X and suppose that the map $\delta_{m-1}$ of $\left(
\ref{LP4}\right)  $ is an isomorphism. Then%
\[
\lambda(\text{Y,}\theta\text{):}=d_{m-1}^{-1}\left(  \mu(\text{Y,}%
\theta\text{)}\right)  \in Ext^{m-1}(I_{\text{Y}},\Omega_{\text{X}}^{n}).
\]
is uniquely defined. We will call $\lambda($Y,$\theta$) the Serre class of the
pair $($Y,$\theta)$ in X$.$
\end{definition}

We will prove the existence of the Serre class $\lambda$($\Delta,\eta)$ of the
diagonal in M$\times$M for the CY manifold M as an element of
$Ext_{\mathcal{O}_{\text{M}\times\text{M}}}^{n-1}(I_{\Delta},\mathcal{O}%
_{\text{M}\times\text{M}})$ by establishing the isomorphism $\left(
\ref{LP4}\right)  $ for $\Delta\subset$M$\times$M for the sheaf $\Omega
_{\text{M}\times\text{M}}^{2n}\approxeq\mathcal{O}_{\text{M}\times\text{M}}.$

\begin{theorem}
\label{GSCY}Let M be a CY manifold of dimension n. Then we have the following
canonical isomorphism:%
\begin{equation}
\delta_{n-1}:Ext_{\mathcal{O}_{\text{M}\times\text{M}}}^{n-1}(I_{\Delta
},\mathcal{O}_{\text{M}\times\text{M}})\rightarrow Ext_{\mathcal{O}%
_{\text{M}\times\text{M}}}^{n}(\mathcal{O}_{\Delta},\mathcal{O}_{\text{M}%
\times\text{M}}). \label{LP3}%
\end{equation}

\end{theorem}

\textbf{Proof:} The proof is based on the following long exact sequence:%
\[
...\rightarrow Ext_{\mathcal{O}_{\text{M}\times\text{M}}}^{n-1}(\mathcal{O}%
_{\text{M}\times\text{M}},\mathcal{O}_{\text{M}\times\text{M}})\rightarrow
Ext_{\mathcal{O}_{\text{M}\times\text{M}}}^{n-1}(I_{\Delta},\mathcal{O}%
_{\text{M}\times\text{M}})\overset{d_{n-1}}{\rightarrow}Ext_{\mathcal{O}%
_{\text{M}\times\text{M}}}^{n}(\mathcal{O}_{\Delta},\mathcal{O}_{\text{M}%
\times\text{M}})\rightarrow
\]%
\begin{equation}
\rightarrow Ext_{\mathcal{O}_{\text{M}\times\text{M}}}^{n}(\mathcal{O}%
_{\text{M}\times\text{M}},\mathcal{O}_{\text{M}\times\text{M}})\rightarrow
Ext_{\mathcal{O}_{\text{M}\times\text{M}}}^{n}(I_{\Delta},\mathcal{O}%
_{\text{M}\times\text{M}})\rightarrow Ext_{\mathcal{O}_{\text{M}\times
\text{M}}}^{n+1}(\mathcal{O}_{\Delta},\mathcal{O}_{\text{M}\times\text{M}%
})\rightarrow.... \label{lp1}%
\end{equation}
From the Grothendieck duality we obtain
\[
Ext_{\mathcal{O}_{\text{M}\times\text{M}}}^{k}(\mathcal{O}_{\text{M}%
\times\text{M}},\mathcal{O}_{\text{M}\times\text{M}})=H^{k}\left(
\text{M}\times\text{M,}\mathcal{O}_{\text{M}\times\text{M}}\right)  ,
\]
for $k\geq0.$ The definition of the CY manifold and the Kunneth formula imply
that for $0<k\neq n$ we have%
\[
Ext_{\mathcal{O}_{\text{M}\times\text{M}}}^{k}(\mathcal{O}_{\text{M}%
\times\text{M}},\mathcal{O}_{\text{M}\times\text{M}})=H^{k}\left(
\text{M}\times\text{M,}\mathcal{O}_{\text{M}\times\text{M}}\right)  =
\]%
\[
=%
%TCIMACRO{\dbigoplus \limits_{p+q=k}}%
%BeginExpansion
{\displaystyle\bigoplus\limits_{p+q=k}}
%EndExpansion
H^{p}\left(  \text{M,}\mathcal{O}_{\text{M}}\right)  \otimes H^{q}\left(
\text{M,}\mathcal{O}_{\text{M}}\right)  =0.
\]
Thus from $\left(  \ref{lp1}\right)  $ we obtain:%
\begin{equation}
0\rightarrow Ext_{\mathcal{O}_{\text{M}\times\text{M}}}^{n-1}(I_{\Delta
},\mathcal{O}_{\text{M}\times\text{M}})\overset{d_{n-1}}{\rightarrow
}Ext_{\mathcal{O}_{\text{M}\times\text{M}}}^{n}(\mathcal{O}_{\Delta
},\mathcal{O}_{\text{M}\times\text{M}})\overset{i_{n}}{\rightarrow
}Ext_{\mathcal{O}_{\text{M}\times\text{M}}}^{n}(\mathcal{O}_{\text{M}%
\times\text{M}},\mathcal{O}_{\text{M}\times\text{M}})\rightarrow...
\label{lp2}%
\end{equation}
We will prove in Proposition \ref{RT} bellow that the map $i_{n}$ in $\left(
\ref{lp2}\right)  $ is zero. Then Theorem \ref{GSCY} follows immediately.
$\blacksquare$

\begin{proposition}
\label{RT}The map
\[
i_{n}:Ext_{\mathcal{O}_{\text{M}\times\text{M}}}^{n}(\mathcal{O}_{\Delta
},\mathcal{O}_{\text{M}\times\text{M}})\rightarrow Ext_{\mathcal{O}%
_{\text{M}\times\text{M}}}^{n}(\mathcal{O}_{\text{M}\times\text{M}%
},\mathcal{O}_{\text{M}\times\text{M}})
\]
in the long exact sequence $\left(  \ref{lp2}\right)  $ is the zero map.
\end{proposition}

\textbf{Proof:} Let us consider the long exact sequence%
\[
...\rightarrow H^{n-1}(\text{M}\times\text{M}-\Delta,\mathcal{O}%
_{\text{M}\times\text{M$-$}\Delta})\overset{\delta_{n-1}}{\rightarrow
}H_{\Delta}^{n}(\text{M}\times\text{M},\mathcal{O}_{\text{M}\times\text{M}%
})\overset{\iota_{n}}{\rightarrow}%
\]%
\begin{equation}
\overset{\iota_{n}}{\rightarrow}H^{n}(\text{M}\times\text{M},\mathcal{O}%
_{\text{M}\times\text{M}})\overset{r_{n}}{\rightarrow}H^{n}(\text{M}%
\times\text{M}-\Delta,\mathcal{O}_{\text{M}\times\text{M$-$}\Delta
})\rightarrow... \label{lp2a}%
\end{equation}
From $\left(  \ref{lp2a}\right)  $ we can conclude that we have%
\[
0\rightarrow H^{n-1}(\text{M}\times\text{M}-\Delta,\mathcal{O}_{\text{M}%
\times\text{M$-$}\Delta})\overset{\delta_{n-1}}{\rightarrow}H_{\Delta}%
^{n}(\text{M}\times\text{M},\mathcal{O}_{\text{M}\times\text{M}}%
)\overset{\iota_{n}}{\rightarrow}%
\]%
\begin{equation}
\overset{\iota_{n}}{\rightarrow}H^{n}(\text{M}\times\text{M},\mathcal{O}%
_{\text{M}\times\text{M}})\overset{r_{n}}{\rightarrow}H^{n}(\text{M}%
\times\text{M}-\Delta,\mathcal{O}_{\text{M}\times\text{M$-$}\Delta
})\rightarrow0 \label{LP3m}%
\end{equation}
We will need the following Lemma:

\begin{lemma}
\label{lp3}The map $\iota_{n}$ in the long exact sequence $\left(
\ref{LP3m}\right)  $ is the zero map.
\end{lemma}

\textbf{Proof: }It is easy to see that $\iota_{n}$~is the zero map if and only
if $r_{n}$ is an isomorphism. We will prove that $r_{n}$ is an isomorphism by contradiction.

The map $r_{n}$ is induced by the restriction map:%
\[
r:\text{M}\times\text{M}\rightarrow\text{M}\times\text{M $-$ }\Delta.
\]
According to Theorem \ref{igor11}%
\begin{equation}
\dim_{\mathbb{C}}H_{\Delta}^{n}(\text{M}\times\text{M},\mathcal{O}%
_{\text{M}\times\text{M}})=1. \label{ig}%
\end{equation}
From the definition of CY manifold and the Kunneth formula we derive that
\begin{equation}
\dim_{\mathbb{C}}H^{n}(\text{M}\times\text{M},\mathcal{O}_{\text{M}%
\times\text{M}})=2 \label{ig1}%
\end{equation}
and $H^{n}($M$\times$M$,\mathcal{O}_{\text{M}\times\text{M}})$ has a basis
$\pi_{1}^{\ast}(\overline{\eta})$ and $\pi_{2}^{\ast}(\overline{\eta}).$ If we
assume that $r_{n}$ is not an isomorphism, then $\left(  \ref{ig}\right)  $,
$\left(  \ref{ig1}\right)  $ and $\left(  \ref{LP3m}\right)  $ will imply
that
\begin{equation}
\dim_{\mathbb{C}}H^{n}(\text{M}\times\text{M}-\Delta,\mathcal{O}%
_{\text{M}\times\text{M$-$}\Delta})\leq1 \label{ig2}%
\end{equation}
Formula $\left(  \ref{ig2}\right)  $ means that on M$\times$M $-$ $\Delta$
some linear combination $\alpha\pi_{1}^{\ast}(\overline{\eta})+\beta\pi
_{2}^{\ast}(\overline{\eta})$ for $\alpha,\beta\in\mathbb{C}$ will represent
zero. So the holomorphic $n$ form $\alpha\pi_{1}^{\ast}(\eta)+\beta\pi
_{2}^{\ast}(\eta)$ will be an exact form on M$\times$M $-$ $\Delta$. So on
M$\times$M$-\Delta$ there exists a holomorphic $n-1$ form $\omega(n-1,0)$ such
that
\begin{equation}
\alpha\pi_{1}^{\ast}(\overline{\eta})+\beta\pi_{2}^{\ast}(\overline{\eta
})=\overline{\partial}\left(  \overline{\omega(n-1,0)}\right)  . \label{ex}%
\end{equation}
Let $\mathcal{U}$ be an affine open set in M. Thus the restriction of the
holomorphic $n-1$ form $\ \omega(n-1,0)$ on $\mathcal{U}\times\mathcal{U}$ $-$
$\left(  \mathcal{U}\times\mathcal{U}\cap\Delta\right)  $ can be represented
as follows:
\[
\omega(n-1,0)=%
%TCIMACRO{\dsum \limits_{I_{p},J_{q};p+q=n-1}}%
%BeginExpansion
{\displaystyle\sum\limits_{I_{p},J_{q};p+q=n-1}}
%EndExpansion
f_{I_{p},J_{q}}(z,w)dz^{I_{[}}\wedge dw^{J_{q}},
\]
where $(z^{1},...,z^{n},w^{1},...,w^{n})$ are local coordinates in
$\mathcal{U}\times\mathcal{U},$
\[
I_{p}=(i_{1},...,i_{p})\text{ for }0<i_{1}<...<i_{p}\leq n,
\]%
\[
J_{q}=(j_{1},...,j_{q})\text{ for }0<j_{1}<...<j_{q}\leq n
\]
and $f_{I_{p},J_{q}}(z,w)$ are holomorphic functions in $\mathcal{U}%
\times\mathcal{U}-\left(  \mathcal{U}\times\mathcal{U}\cap\Delta\right)  .$
Since $\Delta$ has a codimension $n\geq3$ in M$\times$M, the Hartogs principle
implies that the holomorphic functions $f_{I_{p},J_{q}}(z,w)$ are well defined
on $\mathcal{U}\times\mathcal{U}.$ Thus $\omega(n-1,0)$ is a well defined
holomorphic $n-1$ form on M$\times$M. This implies that $H^{n-1}($M$\times
$M$,\mathcal{O}_{\text{M}\times\text{M}})\neq0.$ This fact contradicts the
definition of a CY manifold. Thus the assumption that $r_{n}$ is not an
isomorphism leads to a contradiction. Lemma \ref{lp3} is proved.
$\blacksquare$

\textbf{The end of the proof of Proposition }\ref{RT}: \textbf{ }Theorem
\ref{igor11} implies that there exists a canonical isomorphism%
\[
H_{\Delta}^{n}(\text{M}\times\text{M},\mathcal{O}_{\text{M}\times\text{M}%
})\approxeq Ext_{\mathcal{O}_{\text{M}\times\text{M}}}^{n}\left(
\mathcal{O}_{\Delta},\mathcal{O}_{\text{M}\times\text{M}}\right)  .
\]
The Grothendieck's duality implies that%
\[
Ext_{\mathcal{O}_{\text{M}\times\text{M}}}^{n}\left(  \mathcal{O}%
_{\text{M}\times\text{M}},\mathcal{O}_{\text{M}\times\text{M}}\right)
\approxeq H^{n}(\text{M}\times\text{M},\mathcal{O}_{\text{M}\times\text{M}})
\]
are canonically isomorphic. From the properties of the local cohomology as
stated in \cite{SGA2} and Grothendieck duality it follows that we have the
following commutative diagram:%
\begin{equation}%
\begin{array}
[c]{ccc}%
Ext_{\mathcal{O}_{\text{M}\times\text{M}}}^{n}\left(  \mathcal{O}_{\Delta
},\mathcal{O}_{\text{M}\times\text{M}}\right)  & \overset{i_{n}}{\rightarrow}
& Ext_{\mathcal{O}_{\text{M}\times\text{M}}}^{n}\left(  \mathcal{O}%
_{\text{M}\times\text{M}},\mathcal{O}_{\text{M}\times\text{M}}\right) \\
\parallel &  & \parallel\\
H_{\Delta}^{n}(\text{M}\times\text{M},\mathcal{O}_{\text{M}\times\text{M}}) &
\overset{\iota_{n}}{\rightarrow} & H^{n}(\text{M}\times\text{M,}%
\mathcal{O}_{\text{M}\times\text{M}})
\end{array}
. \label{ig6}%
\end{equation}
The commutative diagram $\left(  \ref{ig6}\right)  $ and Lemma \ref{lp3} imply
Proposition \ref{RT}. $\blacksquare$

\begin{corollary}
\label{GSCYb}The Serre class $\lambda(\Delta,\eta)$ of the diagonal exists.
\end{corollary}

\begin{proposition}
\label{GSCYa}The restriction of the Serre class $\lambda(\Delta,\eta)$ on an
affine open set $\mathcal{U\times U}$ in M$\times$M can be identified with the
holomorphic Bochner-Martinelli kernel $\mathcal{K}_{\mathcal{U\times U}%
}^{0,n-1}$ defined by formula $\left(  \ref{Bocha}\right)  .$
\end{proposition}

\textbf{Proof: }The Grothendieck duality yields the commutative diagram
\begin{equation}%
\begin{array}
[c]{ccc}%
Ext_{\mathcal{O}_{\text{M}\times\text{M}}}^{n-1}\left(  I_{\Delta}%
,\mathcal{O}_{\text{M}\times\text{M}}\right)  & \overset{d_{n-1}}{\rightarrow}
& Ext_{\mathcal{O}_{\text{M}\times\text{M}}}^{n}\left(  \mathcal{O}_{\Delta
},\mathcal{O}_{\text{M}\times\text{M}}\right) \\
\parallel &  & \parallel\\
H^{n-1}(\text{M}\times\text{M $-$ }\Delta\text{,}\mathcal{O}_{\text{M}%
\times\text{M$-$}\Delta}) & \overset{\delta_{n-1}}{\rightarrow} & H_{\Delta
}^{n}(\text{M}\times\text{M},\mathcal{O}_{\text{M}\times\text{M}})
\end{array}
. \label{ig6a}%
\end{equation}
Proposition \ref{GR} implies that we can identify the Grothendieck class
$\mu(\Delta)$ with the antiholomorphic Dirac current $\overline{\eta}_{\Delta
}$ by using the ismorphism bewteen $Ext_{\mathcal{O}_{\text{M}\times\text{M}}%
}^{n}\left(  \mathcal{O}_{\Delta},\mathcal{O}_{\text{M}\times\text{M}}\right)
$ and $H_{\Delta}^{n}($M$\times$M$,\mathcal{O}_{\text{M}\times\text{M}})$ in
$\left(  \ref{ig6a}\right)  $. Theorem \ref{igor2001a} and Theorem \ref{igor9}
imply that the restriction $\delta_{n-1}^{-1}(\overline{\eta}_{\Delta})$ on a
Zariski open set $\mathcal{U\times U}$ in M$\times$M can be identified with
the holomorphic Bochner-Martinelli kernel $\mathcal{K}_{\mathcal{U\times U}%
}^{0,n-1}.$From the definition of the Serre class $\lambda(\Delta
,\eta)=d_{n-1}^{-1}(\mu(\Delta))$ and the commutativity of the diagram
$\left(  \ref{ig6a}\right)  $ Proposition \ref{GSCYa} follows directly.
$\blacksquare$

\begin{remark}
\label{Atserre}One can show that the Serre class corresponding to the
Grothendieck class $\mu(\Delta)$ as defined in Remark \ref{At0} does not
exists. In fact it is not difficult to see that
\[
\mu(\Delta)\in H_{\Delta}^{n}(\text{M}\times\text{M,}\Omega_{\text{M}%
\times\text{M}}^{n}).
\]
We remarked that the Grothendieck class $\mu(\Delta)$ can be identified with
the Dirac current $\delta_{\Delta}.$Thus it follows that in the exact
sequence:%
\[
0\rightarrow H^{n-1}(\text{M}\times\text{M}-\Delta,\Omega_{\text{M}%
\times\text{M}}^{n})\overset{\delta_{n-1}}{\rightarrow}H_{\Delta}^{n}%
(\text{M}\times\text{M},\Omega_{\text{M}\times\text{M}}^{n})\overset{\iota
_{n}}{\rightarrow}%
\]%
\[
\overset{\iota_{n}}{\rightarrow}H^{n}(\text{M}\times\text{M},\Omega
_{\text{M}\times\text{M}}^{n})\overset{r_{n}}{\rightarrow}H^{n}(\text{M}%
\times\text{M}-\Delta,\Omega_{\text{M}\times\text{M}}^{n})\rightarrow0
\]
the map
\[
H_{\Delta}^{n}(\text{M}\times\text{M},\Omega_{\text{M}\times\text{M}}%
^{n})\overset{\iota_{n}}{\rightarrow}H^{n}(\text{M}\times\text{M}%
,\Omega_{\text{M}\times\text{M}}^{n})
\]
is non zero since $\iota_{n}(\delta_{\Delta})=[\omega_{\Delta}],$ where
$[\omega_{\Delta}]$ represents the Poincare dual class of the diagonal. Thus
Theorem \ref{igor11} implies that the map:%
\[
H^{n-1}(\text{M}\times\text{M}-\Delta,\Omega_{\text{M}\times\text{M}}%
^{n})\overset{\delta_{n-1}}{\rightarrow}H_{\Delta}^{n}(\text{M}\times
\text{M},\Omega_{\text{M}\times\text{M}}^{n})
\]
is the zero map. This implies that the map%
\[
Ext_{\mathcal{O}_{\text{M}\times\text{M}}}^{n-1}\left(  I_{\Delta}%
,\Omega_{\text{M}\times\text{M}}^{n}\right)  \overset{d_{n-1}}{\rightarrow
}Ext_{\mathcal{O}_{\text{M}\times\text{M}}}^{n}\left(  \mathcal{O}%
_{\text{M}\times\text{M}},\Omega_{\text{M}\times\text{M}}^{n}\right)
\]
is the zero map. Thus the Serre class $\lambda(\Delta)$ of the diagonal does
not exist.
\end{remark}

Next we will show that the Serre class of the pair (Y$,\theta)$, where Y is a
submanifold of codimension $m$ embedded in a $n$ dimensional CY threefold M
and $\theta$ is a holomorphic $n-m$ form on Y is well defined on M. From the
definition of the Serre class of the pair (Y$,\theta)$ in a CY manifold it
follows that we need to check that the coboundary map $\delta_{m-1}$ in
$\left(  \ref{LP4}\right)  $:
\begin{equation}
d_{m-1}:Ext_{\mathcal{O}_{\text{M}}}^{m-1}(I_{\text{Y}},\mathcal{O}_{\text{M}%
})\rightarrow Ext_{\mathcal{O}_{\text{M}}}^{m}(\mathcal{O}_{\text{Y}%
},\mathcal{O}_{\text{M}}) \label{LP4a}%
\end{equation}
is an isomorphism.

\begin{proposition}
\label{LP5}Let X be a compact K\"{a}hler manifold such that
\begin{equation}
H^{k}(\text{X,}\mathcal{O}_{\text{X}})=0\text{ for }0<k<\dim_{\mathbb{C}%
}\text{X}=n\text{.} \label{1}%
\end{equation}
Then $\delta_{m-1}$ defined by $\left(  \ref{LP4a}\right)  $ is an isomorphism.
\end{proposition}

\textbf{Proof: }Proposition \ref{LP5} follows from the long exact sequence:%
\[
...\rightarrow Ext_{\mathcal{O}_{\text{X}}}^{m-1}(\mathcal{O}_{\text{X}%
},\mathcal{O}_{\text{X}})\rightarrow Ext_{\mathcal{O}_{\text{X}}}%
^{m-1}(I_{\text{Y}},\mathcal{O}_{\text{X}})\rightarrow Ext_{\mathcal{O}%
_{\text{X}}}^{m}(\mathcal{O}_{\text{X}}/I_{\text{Y}},\mathcal{O}_{\text{X}%
})\rightarrow Ext_{\mathcal{O}_{\text{X}}}^{m}(\mathcal{O}_{\text{X}%
},\mathcal{O}_{\text{X}})\rightarrow...
\]
associated with the exact sequence:%
\[
0\rightarrow I_{\text{Y}}\rightarrow\mathcal{O}_{\text{X}}\rightarrow
\mathcal{O}_{\text{X}}/I_{\text{Y}}=\mathcal{O}_{\text{Y}}\rightarrow0.
\]
In fact the Grothendieck duality and the condition $\left(  \ref{1}\right)  $
imply%
\[
Ext_{\mathcal{O}_{\text{X}}}^{m-1}(\mathcal{O}_{\text{X}},\mathcal{O}%
_{\text{X}})=H^{m-1}(\text{X,}\mathcal{O}_{\text{X}})=Ext_{\mathcal{O}%
_{\text{X}}}^{m}(\mathcal{O}_{\text{X}},\mathcal{O}_{\text{X}})=H^{m}%
(\text{X,}\mathcal{O}_{\text{X}})=0
\]
for $0<m<n.$ Proposition \ref{LP5} is proved. $\blacksquare$

Thus Proposition \ref{LP5} guarantees that the Serre class of the pair
(Y$,\theta)$ is well defined on a CY manifold M since for CY manifolds
condition $\left(  \ref{1}\right)  $ holds. In particular we see that if
$(\Sigma,\theta)$ is a pair of a Riemann surface embedded in a CY threefold
and $\theta$ is a holomorphic form on $\Sigma,$ then the Serre class
$\lambda(\Sigma,\theta)$ and the Grothendieck class $\mu(\Sigma,\theta)$ are
well defined on CY threefold M.

Suppose that Y$_{1}$ and Y$_{2}$ are two submanifolds of complex dimension
$m_{1}$ and $m_{2}$ in a CY manifold M such that
\[
\dim_{\mathbb{C}}Y_{1}+\dim_{\mathbb{C}}Y_{2}=\dim_{\mathbb{C}}M-1
\]
and Y$_{1}\cap$Y$_{2}=\varnothing.$ We will define the restriction of the
Serre class $\lambda($Y$_{1}$,$\theta_{1})|_{\text{Y}_{2}}$ as an element of
$Ext_{\mathcal{O}_{\text{M}}}^{m_{1}-1}(\mathcal{O}_{\text{Y}_{2}}%
,\mathcal{O}_{\text{Y}_{2}})$ as follows; Any element of
\[
\alpha\in Ext_{\mathcal{O}_{\text{M}}}^{m_{1}-1}(I_{\text{Y}_{1}}%
,\mathcal{O}_{\text{M}})
\]
corresponds to an exact sequence of length $m_{1}-1$\ consisting of
$\mathcal{O}_{\text{M}}$ modules with the first one being isomorphic to
$\mathcal{O}_{\text{M}}$ and the last one to $I_{\text{Y}_{1}}.$ The condition
Y$_{1}\cap$Y$_{2}=\varnothing$ implies that the ideal sheaf $I_{\text{Y}_{1}}$
when restricted on Y$_{2}$ will be the structure sheaf $\mathcal{O}%
_{\text{Y}_{2}}$ and thus we have two restrictions map
\[
\rho_{1}:I_{\text{Y}_{1}}\rightarrow\mathcal{O}_{\text{Y}_{2}}\text{ and }%
\rho_{2}:\mathcal{O}_{\text{M}}\rightarrow\mathcal{O}_{\text{Y}_{2}}.
\]
Since the maps $\rho_{1}$ and $\rho_{2}$ are surjective then the exact
sequence of length $m_{1}-1$ corresponds to the same exacts sequence tensored
with $\mathcal{O}_{\text{Y}_{2}}$ and we will get a natural map%
\[
r_{m_{1}-1}:Ext_{\mathcal{O}_{\text{M}}}^{m_{1}-1}(I_{\text{Y}},\mathcal{O}%
_{\text{M}})\rightarrow Ext_{\mathcal{O}_{\text{M}}}^{m_{1}-1}(\mathcal{O}%
_{\text{Y}_{2}},\mathcal{O}_{\text{Y}_{2}}).
\]
Then we define
\[
\lambda(Y_{1},\theta_{1})|_{\text{Y}_{2}}=r_{m_{1}-1}\left(  \lambda
(Y_{1},\theta_{1})\right)  .
\]

\subsection{A Homological Interpretation of the Holomorphic Linking}

\begin{proposition}
\label{GR2c}Let Y be a smooth subvariety of codimension $m$ in CY manifold M
which is represented as the intersections of M with $m$ hypersurfaces. Then
\[
H_{\text{Y}}^{m}(\text{M,}\mathcal{O}_{\text{M}})\approxeq Ext_{\mathcal{O}%
_{\text{M}}}^{m}(\mathcal{O}_{\text{Y}},\mathcal{O}_{\text{M}}).
\]

\end{proposition}

\textbf{Proof:} The proof of Proposition \ref{GR2c} is based several Lemmas:

\begin{lemma}
\label{101}Let Y be a subvariety in a CY manifold M of dimension $m$. Then we
have
\[
\underline{Ext}_{\mathcal{O}_{\text{M}}}^{j}(I_{\text{Y}}^{k}/I_{\text{Y}%
}^{k+1},\mathcal{O}_{\text{M}})=0
\]
for $j<m.$
\end{lemma}

\textbf{Proof: }Since Y is a smooth variety in M then Y is a local complete
intersection in M. Thus we can apply $\left(  \ref{lcs2a}\right)  $ to
conclude
\[
\underline{Ext}_{\mathcal{O}_{\text{M}}}^{j}(I_{\text{Y}}^{k}/I_{\text{Y}%
}^{k+1},\mathcal{O}_{\text{M}})=0
\]
for $j<m.$ Lemma \ref{101} is proved. $\blacksquare$

\begin{lemma}
\label{102} We have the following isomorphisms
\begin{equation}
\underline{Ext}_{\mathcal{O}_{\text{M}}}^{m}(\mathcal{O}_{\text{M}%
}/I_{\text{Y}}^{k+1},\mathcal{O}_{\text{M}})\approxeq\underline{Ext}%
_{\mathcal{O}_{\text{M}}}^{m}(\mathcal{O}_{\text{M}}/I_{\text{Y}}%
^{k},\mathcal{O}_{\text{M}}) \label{iso}%
\end{equation}
for $k>0.$
\end{lemma}

\textbf{Proof: }From the exact sequence
\[
0\rightarrow I_{\text{Y}}^{k}/I_{\text{Y}}^{k+1}\rightarrow\mathcal{O}%
_{\text{M}}/I_{\text{Y}}^{k+1}\rightarrow\mathcal{O}_{\text{M}}/I_{\text{Y}%
}^{k}\rightarrow0,
\]
the associated long exact sequence:%
\begin{equation}
...\rightarrow\underline{Ext}_{\mathcal{O}_{\text{M}}}^{j}(\mathcal{O}%
_{\text{M}}/I_{\text{Y}}^{k},\mathcal{O}_{\text{M}})\rightarrow\underline
{Ext}_{\mathcal{O}_{\text{M}}}^{j}(\mathcal{O}_{\text{M}}/I_{\text{Y}}%
^{k+1},\mathcal{O}_{\text{M}})\rightarrow\underline{Ext}_{\mathcal{O}%
_{\text{M}}}^{j}(I_{\text{Y}}^{k}/I_{\text{Y}}^{k+1},\mathcal{O}_{\text{M}%
})\rightarrow... \label{iso0}%
\end{equation}
and Lemma \ref{101}, $\left(  \ref{iso}\right)  $ follows directly. Lemma
\ref{102} is proved. $\blacksquare$

\begin{corollary}
\label{102a}$\underline{Ext}_{\mathcal{O}_{\text{M}}}^{j}(\mathcal{O}%
_{\text{M}}/I_{\text{Y}}^{k},\mathcal{O}_{\text{M}})=0$ for $0<j<m$ and
$k\geq0.$
\end{corollary}

\begin{lemma}
\label{100} $\underline{\mathcal{H}}_{\text{Y}}^{i}(\mathcal{O}_{\text{M}%
}\mathcal{)}=0$ for $i\neq m$ and
\[
\underline{\mathcal{H}}_{\text{Y}}^{m}(\mathcal{O}_{\text{M}}\mathcal{)}%
=\mathcal{O}_{\text{M}}/(f_{1},...,f_{n-m}),
\]
where Y is locally defined by $f_{1}=...=f_{n-m}=0.$
\end{lemma}

\textbf{Proof: }The proof of Lemma \ref{100} follows directly from the
definition of the sheaves $\underline{\mathcal{H}}_{\text{Y}}^{i}%
(\mathcal{O}_{\text{M}}\mathcal{)}$ given by $\left(  \ref{Grot}\right)  ,$
the result stated in \cite{GH} that
\[
\underline{Ext}_{\mathcal{O}_{\text{M}}}^{m}(\mathcal{O}_{\text{M}%
}/I_{\text{Y}},\mathcal{O}_{\text{M}})=\mathcal{O}_{\text{M}}/(f_{1}%
,...,f_{n-m}),
\]
Lemma \ref{102} and Corollary \ref{102a}. Lemma \ref{100} is proved.
$\blacksquare$

\begin{lemma}
\label{103}We have $H^{j}($Y,$I_{\text{Y}}^{k}/I_{\text{Y}}^{k+1})=0$ for
$0\leq j<m$ and $k\geq0.$
\end{lemma}

\textbf{Proof: }The conditions that Y has a codimension $m$ and Y can be
represented as the intersections of $m$ hyperplanes with M imply that the dual
of the normal bundle is direct sum of line bundles of the type $\mathcal{O}%
_{\text{M}}(-n_{i})$ where $n_{i}>0.$ This fact combined with the fact
$I_{\text{Y}}^{k}/I_{\text{Y}}^{k+1}$ is the symmetric $k$ power of the
conormal bundle of Y in M and Kodaira vanishing Theorem imply Lemma \ref{103}.
$\blacksquare$

Proposition \ref{GR2c} follows directly from long exact sequence:%
\begin{equation}
...\rightarrow Ext_{\mathcal{O}_{\text{M}}}^{j}(\mathcal{O}_{\text{M}%
}/I_{\text{Y}}^{k},\mathcal{O}_{\text{M}})\rightarrow Ext_{\mathcal{O}%
_{\text{M}}}^{j}(\mathcal{O}_{\text{M}}/I_{\text{Y}}^{k+1},\mathcal{O}%
_{\text{M}})\rightarrow Ext_{\mathcal{O}_{\text{M}}}^{j}(I_{\text{Y}}%
^{k}/I_{\text{Y}}^{k+1},\mathcal{O}_{\text{M}})\rightarrow... \label{lesa}%
\end{equation}
Lemmas \ref{101}, \ref{102}, \ref{100}, \ref{103} and the definition of
$H_{\text{Y}}^{m}($M,$\mathcal{O}_{\text{M}})$. $\blacksquare$

Proposition \ref{GR2c} implies that the Grothendieck class $\mu$(Y,$\theta
)$\ of (Y,$\theta)$ can be interpreted as a local cohomology class in
$H_{\text{Y}}^{m}($M,$\mathcal{O}_{\text{M}}).$

\begin{proposition}
\label{GR2a}Let Y be a submanifold of CY manifold M of codimension $m.$ Let
$\theta$ be a holomorphic $n-m$ form on Y. Then there exists a natural
identification of the Grothendieck class $\mu($Y$,\theta)$ with the Dirac
antiholomorphic current $\overline{\theta}_{\text{Y}}.$
\end{proposition}

\textbf{Proof: }The definition of the Grothendieck class $\mu($Y$,\theta)$ and
Proposition \ref{GR2c} imply that
\[
\mu(\text{Y},\theta)\in Ext_{\mathcal{O}_{\text{M}}}^{m}(\mathcal{O}%
_{\text{Y}},\mathcal{O}_{\text{M}})=H_{\text{Y}}^{m}(\text{M},\mathcal{O}%
_{\text{M}}).
\]
From here it follows that we can repeat the arguments in the proof of
Proposition \ref{GR} to conclude Proposition \ref{GR2a}. $\blacksquare$

\begin{proposition}
\label{GR2b}Let
\[
\lambda(\text{Y},\theta)\in Ext_{\mathcal{O}_{\text{M}}}^{m-1}(I_{\text{Y}%
},\mathcal{O}_{\text{M}})
\]
be the Serre class of the pair $($Y$,\theta)$ where Y is an intersection of M
with d hypersurfaces and $m=$codim$_{\mathbb{C}}$Y. Then we have the following
isomorphism:%
\begin{equation}
Ext_{\mathcal{O}_{\text{M}}}^{m-1}(\mathcal{O}_{\text{M$-$Y}},\mathcal{O}%
_{\text{M$-$Y}})\approxeq H^{m-1}(\text{M$-$Y},\mathcal{O}_{\text{M$-$Y}}).
\label{Gr2b}%
\end{equation}
Let $r$ be the restriction map composed with the isomorphism $\left(
\ref{Gr2b}\right)  $
\[
r:Ext_{\mathcal{O}_{\text{M}}}^{m-1}(I_{\text{Y}},\mathcal{O}_{\text{M}%
})\rightarrow Ext_{\mathcal{O}_{\text{M}}}^{m-1}(\mathcal{O}_{\text{M$-$Y}%
},\mathcal{O}_{\text{M$-$Y}})\approxeq H^{m-1}(\text{M$-$Y},\mathcal{O}%
_{\text{M$-$Y}}).
\]
Then we have the following expression for the Serre class $\lambda($%
Y$,\theta)$:%
\begin{equation}
r\left(  \lambda(\text{Y},\theta)\right)  =\overline{\partial}^{-1}\left(
\overline{\theta}_{\text{Y}}\right)  . \label{Gr2}%
\end{equation}

\end{proposition}

\textbf{Proof}: We will recall that according to the Grothendieck duality and
since M is a CY manifold we have%
\begin{equation}
Ext_{\mathcal{O}_{\text{M}}}^{p}(\mathcal{O}_{\text{M}},\mathcal{O}_{\text{M}%
})=H^{p}(\text{M,}\mathcal{O}_{\text{M}})=0 \label{sof}%
\end{equation}
for $0<p<n.$ From the long exact sequence $\left(  \ref{lesa}\right)  $ for
$k=1$ and $\left(  \ref{sof}\right)  $ we can conclude that we have the
following isomorphism:%
\begin{equation}
Ext_{\mathcal{O}_{\text{M}}}^{m-1}(I_{\text{Y}},\mathcal{O}_{\text{M}%
})\overset{\delta_{m-1}}{\rightarrow}Ext_{\mathcal{O}_{\text{M}}}%
^{m}(\mathcal{O}_{\text{M}}/I_{\text{Y}},\mathcal{O}_{\text{M}}).
\label{isomY}%
\end{equation}
The long exact sequence%
\[
...\rightarrow H_{\text{Y}}^{k}(\text{M,}\mathcal{O}_{\text{M}})\rightarrow
H^{k}(\text{M,}\mathcal{O}_{\text{M}})\rightarrow H^{k}(\text{M$-$%
Y,}\mathcal{O}_{\text{M$-$Y}})\rightarrow
\]
combined with $\left(  \ref{sof}\right)  $ give the following isomorphism:
\begin{equation}
H^{m-1}(\text{M $-$ Y},\mathcal{O}_{\text{M}})\overset{\delta_{m-1}^{^{\prime
}}}{\rightarrow}H_{\text{Y}}^{m}(\text{M,}\mathcal{O}_{\text{M}}).
\label{Gr2a}%
\end{equation}
According to Proposition \ref{GR2c}
\begin{equation}
Ext_{\mathcal{O}_{\text{M}}}^{m}(\mathcal{O}_{\text{M}}/I_{\text{Y}%
},\mathcal{O}_{\text{M}})\approxeq H_{\text{Y}}^{m}(\text{M},\mathcal{O}%
_{\text{M}}). \label{az}%
\end{equation}
\ From $\left(  \ref{az}\right)  ,$ $\left(  \ref{isomY}\right)  ,$ $\left(
\ref{Gr2a}\right)  $ and the commutative diagram
\begin{equation}%
\begin{array}
[c]{ccc}%
Ext_{\mathcal{O}_{\text{M}}}^{m-1}(I_{\text{Y}},\mathcal{O}_{\text{M}}) &
\overset{\delta_{m-1}}{\rightarrow} & Ext_{\mathcal{O}_{\text{M}}}%
^{m}(\mathcal{O}_{\text{M}},\mathcal{O}_{\text{M}})\\
\downarrow &  & \downarrow\\
H^{m-1}(\text{M $-$ Y},\mathcal{O}_{\text{M}}) & \overset{\delta
_{m-1}^{^{\prime}}}{\rightarrow} & H_{\text{Y}}^{m}(\text{M,}\mathcal{O}%
_{\text{M}})
\end{array}
\label{az0}%
\end{equation}
we can conclude that the map%
\[
r:Ext_{\mathcal{O}_{\text{M}}}^{m-1}(I_{\text{Y}},\mathcal{O}_{\text{M}%
})\rightarrow Ext_{\mathcal{O}_{\text{M}}}^{m-1}(\mathcal{O}_{\text{M$-$Y}%
},\mathcal{O}_{\text{M$-$Y}})\approxeq H^{m-1}(\text{M$-$Y},\mathcal{O}%
_{\text{M$-$Y}})
\]
induced from the restriction map $I_{\text{Y}}$ to M$-$Y is an isomorphism.
Formula $\left(  \ref{Gr2b}\right)  $ is proved. The definition of the
Grothendieck class implies that%
\[
\mu(\text{Y,}\theta)\in Ext_{\mathcal{O}_{\text{M}}}^{m}(\mathcal{O}%
_{\text{M}}/I_{\text{Y}},\mathcal{O}_{\text{M}})\approxeq H_{\text{Y}}%
^{m}(\text{M},\mathcal{O}_{\text{M}}).
\]
The definition of the Serre class $\lambda($Y$,\theta)$ of a pair $($%
Y$,\theta)$ and $\left(  \ref{az0}\right)  $ imply that
\[
\lambda(\text{Y},\theta)\in Ext_{\mathcal{O}_{\text{M}}}^{m-1}(I_{\text{Y}%
},\mathcal{O}_{\text{M}})\approxeq H^{m-1}(\text{M $-$ Y},\mathcal{O}%
_{\text{M}}).
\]
Now from here and Proposition \ref{Gr2a} we get formula $\left(
\ref{Gr2}\right)  .\blacksquare$

\begin{proposition}
\label{IA1}There exists a canonical pairing%
\[
\left\langle \text{ , }\right\rangle :Ext_{\mathcal{O}_{\text{M}}}%
^{k}(\mathcal{O}_{\text{Y}},\mathcal{O}_{\text{Y}})\times Ext_{\mathcal{O}%
_{\text{M}}}^{n-k}(\mathcal{O}_{\text{Y}},\mathcal{O}_{\text{M}}%
)\rightarrow\mathbb{C}.
\]

\end{proposition}

\textbf{Proof:} Yoneda product defines a pairing%
\begin{equation}
\left\langle \text{ , }\right\rangle :Ext_{\mathcal{O}_{\text{M}}}%
^{k}(\mathcal{O}_{\text{Y}},\mathcal{O}_{\text{Y}})\times Ext_{\mathcal{O}%
_{\text{M}}}^{n-k}(\mathcal{O}_{\text{Y}},\mathcal{O}_{\text{M}})\rightarrow
Ext_{\mathcal{O}_{\text{M}}}^{n}(\mathcal{O}_{\text{Y}},\mathcal{O}_{\text{M}%
}). \label{It1}%
\end{equation}
By Grothendieck duality $Ext_{\mathcal{O}_{\text{M}}}^{n}(\mathcal{O}%
_{\text{Y}},\mathcal{O}_{\text{M}})$ is canonically isomorphic to the dual of
\[
Ext_{\mathcal{O}_{\text{M}}}^{0}(\mathcal{O}_{\text{Y}},\mathcal{O}_{\text{M}%
})=H^{0}(\text{Y},\mathcal{O}_{\text{Y}})
\]
which is canonically isomorphic to $\mathbb{C}.$ Proposition \ref{IA1} is
proved. $\blacksquare$

\begin{theorem}
\label{IA}Since M is a CY threefold we know that
\[
\lambda(\Sigma_{1},\theta_{1})|_{\Sigma_{2}}\in Ext_{\mathcal{O}_{\text{M}}%
}^{1}(\mathcal{O}_{\Sigma_{2}},\mathcal{O}_{\Sigma_{2}})\text{ and }\mu
(\Sigma_{2},\theta_{2})\in Ext_{\mathcal{O}_{\text{M}}}^{2}(\mathcal{O}%
_{\Sigma},\mathcal{O}_{\text{M}}).
\]
Let
\[
\left\langle \lambda(\Sigma_{1},\theta_{1})|_{\Sigma_{2}},\mu(\Sigma
_{2},\theta_{2})\right\rangle
\]
be the pairing defined by Proposition \ref{IA1}. Then we have the following
formula:%
\[
\#\left(  \left(  \Sigma_{1},\theta_{1}\right)  ,\left(  \Sigma_{2},\theta
_{2}\right)  \right)  =\left\langle \lambda(\Sigma_{1},\theta_{1}%
)|_{\Sigma_{2}},\mu(\Sigma_{2},\theta_{2})\right\rangle .
\]

\end{theorem}

\textbf{Proof: }Since
\[
\lambda(\Sigma_{1},\theta_{1})|_{\Sigma_{2}}\in Ext_{\mathcal{O}_{\text{M}}%
}^{1}(\mathcal{O}_{\Sigma_{2}},\mathcal{O}_{\Sigma_{2}})\text{ and }\mu
(\Sigma_{2},\theta_{2})\in Ext_{\mathcal{O}_{\text{M}}}^{2}(\mathcal{O}%
_{\Sigma},\mathcal{O}_{\text{M}})
\]
by Propositions \ref{GR2a}, \ref{GR2b} we can identify $\lambda(\Sigma
_{1},\theta_{1})|_{\Sigma_{2}}$ and $\mu(\Sigma_{2},\theta_{2})$ with
$\overline{\partial}^{-1}\left(  \overline{\theta}_{\Sigma_{1}}\right)  $ and
$\mu\left(  \overline{\theta}_{\Sigma_{2}}\right)  .$ According to
Grothendieck the Yoneda pairing is the same as the Serre pairing which is just
integration. See \cite{SGA2}. The proof of Theorem \ref{IA} follows directly
from\ the analytic formula $\left(  \ref{W3a}\right)  $ for the holomorphic
linking. $\blacksquare$

The condition $\Sigma_{1}\cap\Sigma_{2}=\varnothing$ implies that $\Sigma
_{1}\times\Sigma_{2}\cap\Delta=\varnothing$ thus we get that $I|_{\Sigma
_{1}\times\Sigma_{2}}\approxeq\mathcal{O}_{\Sigma_{1}\times\Sigma_{2}}.$ From
here we conclude for that the Serre class
\begin{equation}
\lambda(\Delta,\eta)|_{\Sigma_{1}\times\Sigma_{2}}\in Ext_{\mathcal{O}%
_{\text{M}\times\text{M}}}^{2}(\mathcal{O}_{\Sigma_{1}\times\Sigma_{2}%
},\mathcal{O}_{\Sigma_{1}\times\Sigma_{2}}). \label{ia}%
\end{equation}
The Grothendieck class $\mu(\Sigma_{1}\times\Sigma_{2},\pi_{1}^{\ast}%
(\theta_{1})\wedge\pi_{2}^{\ast}(\theta_{2}))$ of the pair
\[
(\Sigma_{1}\times\Sigma_{2},\pi_{1}^{\ast}(\theta_{1})\wedge\pi_{2}^{\ast
}(\theta_{2}))
\]
is an element of $Ext_{\mathcal{O}_{\text{M}\times\text{M}}}^{4}%
(\mathcal{O}_{\Sigma_{1}\times\Sigma_{2}},\mathcal{O}_{\text{M}\times\text{M}%
}),$ i.e.
\[
\mu(\Sigma_{1}\times\Sigma_{2},\pi_{1}^{\ast}(\theta_{1})\wedge\pi_{2}^{\ast
}(\theta_{2}))\in Ext_{\mathcal{O}_{\text{M}\times\text{M}}}^{4}%
(\mathcal{O}_{\Sigma_{1}\times\Sigma_{2}},\mathcal{O}_{\text{M}\times\text{M}%
}).
\]
Yoneda pairing%
\[
\left\langle \text{ },\text{ }\right\rangle :Ext_{\mathcal{O}_{\text{M}%
\times\text{M}}}^{2}(\mathcal{O}_{\Sigma_{1}\times\Sigma_{2}},\mathcal{O}%
_{\Sigma_{1}\times\Sigma_{2}})\times Ext_{\mathcal{O}_{\text{M}\times\text{M}%
}}^{4}(\mathcal{O}_{\Sigma_{1}\times\Sigma_{2}},\mathcal{O}_{\text{M}%
\times\text{M}})\rightarrow
\]%
\begin{equation}
\rightarrow Ext_{\mathcal{O}_{\text{M}\times\text{M}}}^{6}(\mathcal{O}%
_{\Sigma_{1}\times\Sigma_{2}},\mathcal{O}_{\text{M}\times\text{M}}).
\label{ia1}%
\end{equation}
defines a non degenerate pairing since by Grothendieck duality
\[
Ext_{\mathcal{O}_{\text{M}\times\text{M}}}^{6}(\mathcal{O}_{\Sigma_{1}%
\times\Sigma_{2}},\mathcal{O}_{\text{M}\times\text{M}})\approxeq H^{0}%
(\Sigma_{1}\times\Sigma_{2},\mathcal{O}_{\Sigma_{1}\times\Sigma_{2}%
})=\mathbb{C}.
\]
Next we will give a homological algebra interpretation of Theorem \ref{igor14}.

\begin{theorem}
\label{Deny}The following formula holds%
\begin{equation}
\#\left(  \left(  \Sigma_{1},\theta_{1}\right)  ,\left(  \Sigma_{2},\theta
_{2}\right)  \right)  =\left\langle \lambda(\Delta,\eta)|_{\Sigma_{1}%
\times\Sigma_{2}},\mu(\Sigma_{1}\times\Sigma_{2},\pi_{1}^{\ast}(\theta
_{1})\wedge\pi_{2}^{\ast}(\theta_{2}))\right\rangle . \label{16}%
\end{equation}

\end{theorem}

\textbf{Proof: }According to Corollary \ref{GSCYa} the restriction of
$\lambda(\Delta,\eta)$ on the product $\mathcal{U\times U}$ of \ an affine set
$\mathcal{U}$ in M can be identified with the holomorphic Bochner-Martinelli
kernel $\mathcal{K}_{\mathcal{U}\times\mathcal{U}}^{0,2}.$ Thus we can
identify the restriction of the Serre class $\lambda(\Delta,\eta)$ on $\left(
\Sigma_{1}\times\Sigma_{2}\right)  \cap\mathcal{U\times U}$ with the
restriction of the holomorphic Bochner-Martinelli kernel on $\left(
\Sigma_{1}\times\Sigma_{2}\right)  \cap\mathcal{U\times U}.$ Proposition
\ref{GR2a} implies that the Grothendieck class $\mu(\Sigma_{1}\times\Sigma
_{2},\pi_{1}^{\ast}(\theta_{1})\wedge\pi_{2}^{\ast}(\theta_{2}))$ of the pair
$(\Sigma_{1}\times\Sigma_{2},\pi_{1}^{\ast}(\theta_{1})\wedge\pi_{2}^{\ast
}(\theta_{2}))$ can be identified with with the antiholomorphic Dirac current
$\left(  \overline{\pi_{1}^{\ast}(\theta_{1})\wedge\pi_{2}^{\ast}(\theta_{2}%
)}\right)  _{\Sigma_{1}\times\Sigma_{2}}.$ We pointed out above that the
Yoneda pairing is the same as the Serre pairing which is just integration.
Formula $\left(  \ref{16}\right)  $ follows directly from the above described
identifications of the restrictions of the Serre class and the Grothendieck
class, the interpretation of the Yoneda pairing as integration and formula
$\left(  \ref{Igo}\right)  $ in Theorem \ref{b1}. Theorem \ref{Deny} is
proved. $\blacksquare$

One can generalize formulas $\left(  \ref{Igo}\right)  $ and $\left(
\ref{BK96}\right)  $ and Theorems \ref{IA}, \ref{Deny} for the pairs
(Y$_{1},\theta_{1})$ and (Y$_{2},\theta_{2})$ of submanifolds Y$_{1}$ and
Y$_{2}$ in a CY manifold M of any dimension such that
\[
\dim_{\mathbb{C}}\text{Y}_{1}+\dim_{\mathbb{C}}\text{Y}_{2}=\dim_{\mathbb{C}%
}\text{M}-1
\]
and Y$_{1}\cap$Y$_{2}=\varnothing$ to define their holomorphic linking.

\section{\label{At}Generalizations of Holomorphic Linking}

In this section we would like to establish an explicit connection between our
holomorphic linking and the one studied by Atiyah. Atiyah arrived at his
formula for holomorphic linking by considering the twistor transform of the
Green function of the Laplacian. Thus he only defined linking of spheres while
our definition makes sense for curves with genus greater than zero. Also the
twistor space is never a CY space required in our approach. Thus in order to
include the Atiyah holomorphic linking in our picture, we need to extend
further our construction.

\subsection{The Holomorphic Linking of Riemann Surfaces with Marked Points}

We would like to generalize the holomorphic linking of two Riemann surfaces in
a CY manifold to the case of Riemann surfaces with punctures. Let%
\[
(\Sigma_{i};p_{1},..,p_{m_{i}};\theta_{i}),\text{ }i=1,2
\]
be two Riemann surfaces with m$_{i}$ marked points on them embedded in a CY
threefold M and let $\theta_{i}$ be meromorphic forms on $\Sigma_{i}$ with
poles of order at most one at $p_{1},...,p_{m_{i}}.$ We will assume as before
that
\[
\Sigma_{1}\cap\Sigma_{2}=\emptyset.
\]
We can define the holomorphic linking of
\[
\left(  \Sigma_{i};p_{1},...,p_{m_{i}};\theta_{i}\right)
\]
in any of the equivalent ways considered above with appropriate modifications.
For example, the definition of the holomorphic linking via the Green kernel is
defined as follows:
\[
\#((\Sigma_{1};p_{1},...,p_{m_{1}};\theta_{1}),(\Sigma_{2};p_{1},...,p_{m_{2}%
};\theta_{2}))=
\]%
\[
\underset{\varepsilon\rightarrow0}{\lim}\int_{\Sigma_{1}-D_{1,\varepsilon
}\times\Sigma_{2}-D_{2,\varepsilon}}\left(  \overline{\partial}^{-1}\left(
\overline{\eta}_{\Delta}\right)  |_{\Sigma_{1}\times\Sigma_{2}}\right)
\wedge\pi_{1}^{\ast}(\theta_{1})\wedge\pi_{2}^{\ast}(\theta_{2}),
\]
where $D_{i,\varepsilon}$ are the union of disks with radius $\varepsilon$
around the marked points $p_{1},...,p_{m_{i}}$ with fixed local coordinates.
We denote the limit as before by
\[
\int_{\Sigma_{1}\times\Sigma_{2}}\left(  \overline{\partial}^{-1}\left(
\overline{\eta}_{\Delta}\right)  |_{\Sigma_{1}\times\Sigma_{2}}\right)
\wedge\pi_{1}^{\ast}(\theta_{1})\wedge\pi_{2}^{\ast}(\theta_{2}).
\]
Repeating the arguments in the proof of Theorem \ref{igor19}, one can also
derive the geometric formula for the holomorphic linking of the punctured
spheres
\[
\#((\Sigma_{1};p_{1},...,p_{m_{1}};\theta_{1}),(\Sigma_{2};p_{1},...,p_{m_{2}%
};\theta_{2}))=\underset{x\in S_{1}\cap\Sigma_{2}}{\sum}\frac{\omega
_{1}(x)\wedge\theta_{2}(x)}{\eta(x)},
\]
where $\omega_{1}$ is a meromorphic 2-form on the complex surface $S_{1}$
containing $\Sigma_{1}$, whose residue on $\Sigma_{1}$ is $\theta_{1}.$

One can also generalize the notion of holomorphic linking to the case when M
is not necessary CY manifold. In the general case a holomorphic form $\eta$
does not exist and we have to fix a meromorphic three form. In certain cases,
(for example, when M is a Fano variety) we can uniquely fix a meromorphic
3-form $\eta$ by prescribing its poles along a given hypersurface H$\subset$M
whose homology class represents the first Chern class of M.

Another generalization of the holomorphic linking of Riemann surfaces in an
arbitrary complex manifolds can be obtained by considering still holomorphic
forms $\theta_{1},$ $\theta_{2}$ and $\eta$ but coupled with a line bundle.
The simplest example of this sort is the linking of two spheres in a three
dimensional projective space. In this example one can consider unique up to a
scale holomorphic forms $\theta_{1},$ $\theta_{2}$ and $\eta$ with values in
natural line bundles. Then the holomorphic linking becomes a well defined
number depending only on the embeddings of the Riemann surfaces.

In fact the two generalizations of the holomorphic linking of Riemann surfaces
in CY manifolds are closely related when the holomorphic forms are replaced by
holomorphic forms with values in the line bundles on $\Sigma_{1},$ $\Sigma
_{2}$ defined by the fixed points. Thus the we can replace the holomorphic
form $\eta$ on M with a meromorphic form with a simple pole along some divisor
and $\theta_{i}$ can be replaced by meromorphic forms with values in the
trivial bundle. In particular the example of the linking of two spheres in the
$\mathbb{CP}^{3}$ with natural line bundles correspond to the holomorphic
linking of two spheres with two marked points embedded in $\mathbb{CP}^{3}$
and the poles of meromorphic three form is determined by hypersurface
consisting of four hyperplanes in general position. In the next subsection we
will consider this example in a more general context of twistor spaces first
studied by Atiyah in \cite{A}.

\subsection{Relations with Atiyah's Results on Linking of Rational Curves in
Twistor Space}

In \cite{A} Atiyah discovered a version of holomorphic linking while studying
the twistor transform of the Green function of the Laplacian. In order to
express his result in an invariant form he used the Serre classes of the
pairs
\[
\Delta\subset\text{M}\times\text{M}%
\]
and
\[
\Sigma_{i}\subset\text{M},
\]
$i=1,2,$ where M (denoted by Z in his paper) is the twistor space of a compact
four-dimensional manifold N with a self dual metric, and $\Sigma_{i}$ are
rational curves that appear as the praimages of the points of N in M in the
twistor transform. If N is a spin-manifold, there is a natural holomorphic
line bundle $L$ on M such that
\[
L^{-4}\approxeq\Omega_{\text{M}}^{3}.
\]
We denote by $\mathcal{O}_{\text{M}}(n)$ the sheaf of sections of $L^{n}.$ If
N is not a spin-manifold, then $L$ does not exist but $L^{2}$ always exists
yielding $\mathcal{O}_{\text{M}}(n)$ for even $n.$ In particular, there is a
canonical 3-form $\eta^{0}$ with coefficients in $\mathcal{O}_{\text{M}}(4)$,
which plays a similar role for the twistor space as the holomorphic volume
form $\eta$ for the Calabi-Yau space. Since $\Sigma_{i}$ are spheres there are
also canonical one forms $\theta_{i}^{0}$ with coefficients in $\mathcal{O}%
_{\Sigma_{i}}(2).$ Thus the forms
\[
\theta_{i}^{0}\in H^{0}(\Sigma_{i},\mathcal{O}_{\Sigma_{i}})
\]
and
\[
\eta^{0}\in H^{0}(\text{M},\mathcal{O}_{\text{M}}).
\]

If we apply directly twice the Grothendieck duality to the pairs
\[
(\Delta,\eta^{0})\text{ and }(\Sigma_{1}\times\Sigma_{2},\pi_{1}^{\ast}%
(\theta_{1}^{0})\wedge\pi_{1}^{\ast}(\theta_{1}^{0}))
\]
as we did in Section 7.2\textbf{ }we can define the analogues of Grothendieck
classes
\[
\mu(\Delta,\eta^{0})\text{ and }\mu(\Sigma_{1}\times\Sigma_{2},\pi_{1}^{\ast
}(\theta_{1}^{0})\wedge\pi_{1}^{\ast}(\theta_{1}^{0}))
\]
and Serre classes
\[
\lambda(\Delta,\pi_{1}^{\ast}(\eta^{0})|_{\Delta}),\text{ }\lambda(\Sigma
_{1}\times\Sigma_{2},\pi_{1}^{\ast}(\theta_{1}^{0})\wedge\pi_{1}^{\ast}%
(\theta_{1}^{0})).
\]
Since
\[
\mathcal{O}_{\text{M}\times\text{M}}(-2,-2))\approxeq\Omega_{\text{M}%
\times\text{M}}^{6}%
\]
it is easy to see that
\[
\lambda(\Delta,\eta^{0})\in Ext_{\mathcal{O}_{\text{M}\times\text{M}}}%
^{2}(I_{\Delta},\mathcal{O}_{\text{M}\times\text{M}}(-2,-2)).\text{ }%
\]
and
\[
\mu(\Sigma_{1}\times\Sigma_{2},\pi_{1}^{\ast}(\theta_{1}^{0})\wedge\pi
_{1}^{\ast}(\theta_{1}^{0}))\in Ext_{\mathcal{O}_{\text{M}\times\text{M}}}%
^{4}(\mathcal{O}_{\Sigma_{1}\times\Sigma_{2}},\mathcal{O}_{\text{M}%
\times\text{M}}(-2,-2)).\text{ }%
\]
In the same manner we can define the Grothendieck $\mu((\Sigma_{i},\theta
_{i}^{0})$ and Serre classes $\lambda(\Sigma_{i},\theta_{i}^{0})$\ of the pair
$(\Sigma_{i},\theta_{i}^{0})$ and the definitions imply
\[
\lambda(\Sigma_{i},\theta_{i}^{0})\in Ext_{\mathcal{O}_{\text{M}}}%
^{1}(I_{\Sigma_{i}},\mathcal{O}_{\text{M}}(-2))
\]
and%
\[
\mu((\Sigma_{i},\theta_{i}^{0})\in Ext_{\mathcal{O}_{\text{M}}}^{2}%
(\mathcal{O}_{\Sigma_{i}},\mathcal{O}_{\text{M}}(-2)).
\]

By repeating the arguments from Proposition \ref{IA1} we can see that the
Grothendieck duality defines the following two non degenerate pairings:%
\begin{equation}
Ext_{\mathcal{O}_{\text{M}}}^{1}(\mathcal{O}_{\Sigma_{2}},\Omega_{\Sigma_{2}%
}^{1})\times Ext_{\mathcal{O}_{\text{M}}}^{2}(\Omega_{\Sigma_{2}}%
^{1},\mathcal{O}_{\text{M}}\mathcal{(}-4))\rightarrow\mathbb{C} \label{At1}%
\end{equation}
and
\begin{equation}
Ext_{\mathcal{O}_{\text{M}}}^{2}(\mathcal{O}_{\Sigma_{1}\times\Sigma_{2}%
},\Omega_{\Sigma_{1}\times\Sigma_{2}}^{2})\times Ext_{\mathcal{O}%
_{\text{M}\times\text{M}}}^{4}(\Omega_{\Sigma_{1}\times\Sigma_{2}}%
^{2},\mathcal{O}_{\text{M}}\mathcal{(}-4,-4))\rightarrow\mathbb{C}.
\label{At2}%
\end{equation}
Let
\[
\lambda(\Sigma_{1},\theta_{1})|_{\Sigma_{2}}\in Ext_{\mathcal{O}_{\text{M}}%
}^{1}(\mathcal{O}_{\Sigma_{2}},\Omega_{\Sigma_{2}}^{1})\text{ and }%
\lambda\left(  \Delta,\eta^{0}\right)  |_{\Sigma_{1}\times\Sigma_{1}}\in
Ext_{\mathcal{O}_{\text{M}}}^{2}(\mathcal{O}_{\Sigma_{1}\times\Sigma_{2}%
},\Omega_{\Sigma_{1}\times\Sigma_{2}}^{2})
\]
be the restrictions of the Serre classes of the pairs $(\Sigma_{1},\theta
_{1})$ and $\left(  \Delta,\pi_{1}^{\ast}(\eta^{0})|_{\Delta}\right)  $ on
$\Sigma_{2}.$ Then we can define the holomorphic linking of the pair
$(\Sigma_{i},\theta_{i}^{0})$ by using $\left(  \ref{At1}\right)  $ and
$\left(  \ref{At2}\right)  $ in two different ways according to Theorems
\ref{IA} and \ref{Deny}:%
\begin{equation}
At(\Sigma_{1},\Sigma_{2})=\left\langle \lambda(\Sigma_{1},\theta_{1}%
^{0})|_{\Sigma_{2}},\mu(\Sigma_{2},\theta_{2}^{0})\right\rangle \label{At3}%
\end{equation}
and
\begin{equation}
At(\Sigma_{1},\Sigma_{2})=\left\langle \lambda(\Delta,\eta^{0})|_{\Sigma_{2}%
},\mu(\Sigma_{1}\times\Sigma_{2},\pi_{1}^{\ast}(\theta_{1}^{0})\wedge\pi
_{2}^{\ast}(\theta_{2}^{0})\right\rangle . \label{At4}%
\end{equation}
Then one can prove that that formulas $\left(  \ref{At3}\right)  $ and
$\left(  \ref{At4}\right)  $ are the same up to a constant with Atiyah's
formulas for the linking of $\Sigma_{1}$and $\Sigma_{2}$ expressed as the
value of the Green function $G(x,y)$ of the conformally invariant Laplacian
$\square$ on a self dual compact four manifold $N$ of the rational curves
$\Sigma_{x}$ and $\Sigma_{y}$ in N corresponding to $x,y\in$N.

We expect that the various forms of the holomorphic linking considered in this
paper have appropriate analogues for the Atiyah linking number. The
possibility of the path integral presentation has been noted by Gerasimov
\cite{G}. The analytic expression for the Atiyah linking number should
coincide with the Penrose integral formula for the Green function of the
Laplacian. Finally the geometric formula should have the following form
similar to $\left(  \ref{IGOR}\right)  $:%
\[
At(\Sigma_{1},\Sigma_{2})=%
%TCIMACRO{\dsum \limits_{x\in\text{S}_{1}\cap\Sigma_{2}}}%
%BeginExpansion
{\displaystyle\sum\limits_{x\in\text{S}_{1}\cap\Sigma_{2}}}
%EndExpansion
\text{Res}_{x}\left(  \frac{\eta_{1}^{0}}{\eta^{0}}|_{\Sigma_{2}}\right)
\theta_{2}^{0}%
\]
where $\eta_{1}^{0}$ is a twisted holomorphic 3-form on M whose double residue
is the 1-form $\theta_{2}^{0}$ on $\Sigma_{1}.$

In order to relate the Atiyah holomorphic linking with ours, we will pick out
a pair of marked points for each rational curve $\Sigma_{i}$ in M$.$ In the
homogeneous coordinates $(z_{0}:z_{1})$, the marked points can be represented
by linear functionals, say p$_{1}(z),$ p$_{2}(z)$ for $\Sigma_{1}$ and
p$_{3}(z),$ p$_{4}(z)$ for $\Sigma_{2}.$ Then the meromorphic forms with
simple poles at the marked points will have the forms
\[
\theta_{1}=\frac{\theta_{1}^{0}}{\text{p}_{1}(z)\text{p}_{2}(z)}\text{ and
}\theta_{2}=\frac{\theta_{2}^{0}}{\text{p}_{3}(z)\text{p}_{4}(z)}.\text{ }%
\]
Suppose that the twistor space is a projective variety and there is a section
of the canonical bundle on M that yields a hypersurface H$\subset$M that
passes through the four marked points and M$-$H is an open CY manifold. Let
$\eta$ be a meromorphic form M with a singularities along the hypersurface H.
Then we define the holomorphic linking of a rational curves with two marked
points as discussed in \textbf{Section 8.1.} We conjecture that the Atiyah
linking coincides with ours, namely,
\[
At(\Sigma_{1},\Sigma_{2})=\#((\Sigma_{1},\theta_{1}),(\Sigma_{2},\theta
_{2})).
\]

We will illustrate this relation in the case of the twistor space of the four
dimensional sphere $S^{4}.$ In this case M=$\mathbb{P}^{3}$ and its canonical
bundle is isomorphic to $\mathcal{O}_{\text{M}}(-4).$ Thus $H^{0}($%
M,$\Omega_{\text{M}}^{3}(4))$ is one dimensional and its generator can be
written in the projective coordinates as follows
\[
\eta^{0}:=\sum_{i=0}^{3}(-1)^{i}z_{i}dz_{0}\wedge...\wedge d\hat{z}_{i}%
\wedge...\wedge dz_{3}.
\]
Let us consider two non intersecting lines $\Sigma_{1}$ and $\Sigma_{2}$ in
$\mathbb{P}^{3}$ given by pairs of linear functionals
\[
l_{1}(z)=l_{2}(z)=0
\]
and
\[
l_{3}(z)=l_{4}(z)=0.
\]
We denote by $S_{i}$ a hyperplane in M given by the equation $l_{i}(z)=0.$ We
also consider
\[
\theta_{1}^{0}:=\frac{\eta^{0}}{l_{1}(z)l_{2}(z)},
\]
the meromorphic 3-form on M whose double residue is the 1-form $\eta_{1}^{0}$
on $\Sigma_{1}.$ So we obtain%
\begin{equation}
At(\Sigma_{1},\Sigma_{2})=\underset{x\in S_{1}\cap\Sigma_{2}}{\sum}res\left(
\frac{\eta_{1}^{0}}{\eta^{0}}\theta_{2}^{0}\right)  =%
%TCIMACRO{\dsum \limits_{x\in S_{1}\cap\Sigma_{2}}}%
%BeginExpansion
{\displaystyle\sum\limits_{x\in S_{1}\cap\Sigma_{2}}}
%EndExpansion
res\left(  \frac{z_{0}dz_{1}-z_{1}dz_{0}}{l_{1}(z)l_{2}(z)}\right)  .
\label{At6}%
\end{equation}
On the other hand, we can compute the holomorphic linking of two
nonintersecting lines $\Sigma_{1}$and $\Sigma_{2}$ with two marked points also
using the geometric formula. Let p$_{i}(z)=0,$ $i=1,2,2,4$ defines four
hyperplanes H$_{i},$ which are in general position and also intersect the
hyperplanes defined by $l_{i}(z)=0,$ $i=1,2,3,4$ transversely. The pairs of
marked points on $\Sigma_{1}$ and $\Sigma_{2}$ are precisely their
intersection with the hyperplanes H$_{1},$H$_{2}$ and H$_{3},$ H$_{4},$
respectively. We can consider M $-$ $\left\{  \cup_{i=1}^{4}H_{i}\right\}  $
as an open CY manifold with a holomorphic form
\[
\eta:=\frac{\eta^{0}}{p_{1}(z)p_{2}(z)p_{3}(z)p_{4}(z)}.
\]
We will also define a meromorphic 3-form on M $-$ $\left\{  \cup_{i=1}%
^{4}H_{i}\right\}  $%
\[
\eta_{1}:=\frac{\eta^{0}}{p_{1}(z)p_{2}(z)l_{1}(z)l_{2}(z)}%
\]
whose double residue is the 1-form $\theta_{1}$ on $\Sigma_{1}.$ Then the
geometric formula for the holomorphic linking yields
\[
\#((\Sigma_{1},\theta_{1}),(\Sigma_{2},\theta_{2}))=
\]%
\[
\underset{S_{1}\cap\Sigma_{2}}{\sum}res\left(  \frac{\eta_{1}}{\eta}\theta
_{2}\right)  =\underset{S_{1}\cap\Sigma_{2}}{\sum}res\left(  \frac{\frac
{\eta_{1}^{0}}{\text{p}_{1}(z)\text{p}_{2}(z)}}{\frac{\eta^{0}}{\text{p}%
_{1}(z)\text{p}_{2}(z)\text{p}_{3}(z)\text{p}_{4}(z)}}\frac{\theta_{2}^{0}%
}{\text{p}_{3}(z)\text{p}_{4}(z)}\right)  ,
\]
is the Atiyah holomorphic linking $At(\Sigma_{1},\Sigma_{2}).$

We would like to conclude the discussion of the relation between the Atiyah
holomorphic linking and the one arising from the complex CSW theory with a
remark concerning the non abelian case. Atiyah noted in his paper \cite{A}
that the conformal Laplacian has a natural covariant analogue for a unitary
vector bundle $E$ with a connection over X and if the connection is self-dual
the $E$ lifts to a holomorphic bundle $\widetilde{E}$ on M. In this case the
Green function of the conformally invariant Laplacian can be expressed using a
non abelian version of an Atiyah holomorphic linking. In this case, the Serre
class of the diagonal $\Delta\subset$M$\times$M should be extended to
endomorphisms of $\widetilde{E}$ and yields a non abelian version of the
holomorphic linking of two rational curves in M.

A variant of Atiyah construction for arbitrary pairs of nonintersecting curves
in a CY manifold will also yield a non abelian generalization of the
holomorphic linking studied in this paper. We conjecture that it can also be
derived from the complex CSW theory with a gauge group GL($n,\mathbb{C})$
where $n$ is the rank of $\widetilde{E}.$

\end{document}